\documentclass[12pt,amsfonts,amssymb,leqno]{article}
\newfont\got{eufm10}

\usepackage{amssymb}
\usepackage[all]{xy}

\newtheorem{proposition}{Proposition}[section]
\newtheorem{thm}[proposition]{Theorem}
\newtheorem{cor}[proposition]{Corollary}
\newtheorem{lemma}[proposition]{Lemma}
\newtheorem{defn}[proposition]{Definition}

\newcommand{\zerohandgrenade}{0^{\: \mid^{\! \! \! \bullet}}}

\newcommand{\forces}{\ \ |\!\!\!|\!- \ }

\newcounter{secnum}

\begin{document}
\setcounter{section}{-1}

\begin{center}

\vspace*{2cm}

{\Large \bf The core model for almost linear iterations}

\end{center}

\begin{center}
\renewcommand{\thefootnote}{\fnsymbol{footnote}}
\renewcommand{\thefootnote}{arabic{footnote}}
\renewcommand{\thefootnote}{\fnsymbol{footnote}}
{\large Ralf-Dieter Schindler}
\end{center}
\begin{center} 
{\footnotesize
{\it Institut f\"ur Formale Logik, Universit\"at Wien, 1090 Wien, Austria}}
\end{center}

\begin{center}
{\tt 
rds@logic.univie.ac.at}
\end{center}

\begin{abstract}
\noindent
We introduce $\zerohandgrenade$ (``zero hand-grenade'') as a sharp for an inner model
with a proper class of strong cardinals. 
We 
prove the existence of the core model 
$K$ in the theory ``${\sf ZFC} \ + \ \zerohandgrenade$ doesn't
exist.'' Combined with work of Woodin, Steel, and earlier work of the author,
this provides the last step for determining the exact
consistency strength of the assumption in the statement of the
$12^{{\rm th}}$ Delfino problem (cf.
\cite{KeMaSt}). 
\end{abstract}


\section{Introduction.}

Core models were constructed in the papers \cite{dodd-jensen}, \cite{koepke},
\cite{MOZ}, \cite{mitchell} and \cite{mitchell2}, \cite{jensen} (see also
\cite{thesis}), \cite{CMIP}, and \cite{CMMW}.
We refer the reader to \cite{jensen2}, \cite{mitchell3}, and \cite{st-loewe}
for less painful introductions into core model theory.

A core model is intended to be an inner model of set theory (that is,
a transitive
class-sized model of ${\sf ZFC}$) which meets two requirements:

{\bf F}${}_1 \ \ $ It is close to $V$ ($=$ the universe of all sets), and

{\bf F}${}_2 \ \ $ it can be analyzed in great detail.

\noindent Both requirements should be formulated more precisely, of course. However,
as
``core model'' is no formal concept, we can't expect a thorough
{\em general} definition.
Let's try to give some hints.

As for {\bf F}${}_1$, a core model, call it $K$,
should reflect the large cardinal situation of $V$
(for example, the large cardinals that exist in $V$ should be found
in $K$ as well up to a certain size), it should satisfy certain forms of covering
(for example, $K$ should compute successors of singular cardinals correctly),
and it should be absolute for set-forcings (i.e., the definition of $K$
should determine the
same object in any set-forcing extension of $V$).
As for {\bf F}${}_2$, a core model should be a fine structural inner model (a
class-sized premouse, technically speaking), which satisfies certain forms of
condensation (which in turn typically follow from a certain amount of iterability of
$K$, and make $K$ amenable for combinatorial studies).

For the purposes of this introduction we'll say that $K$ exists if there is an inner
model satisfying appropriate versions of {\bf F}${}_1$ and {\bf F}${}_2$ above.
The key strategy for constructing core models has not changed through time.
As a matter of fact, in order to have any chance to build $K$ satisfying
{\bf F}${}_1$ and {\bf F}${}_2$ one has to work in a theory $${\sf ZFC} \ + \ \lnot \ \Psi
{\rm , }$$ where $\lnot \ \Psi$ denotes an anti large cardinal assumption (for example,
``there is no inner model with a certain large cardinal''). The following table
lists the achievements of \cite{dodd-jensen}, \cite{koepke},
\cite{MOZ}, \cite{mitchell} and \cite{mitchell2}, \cite{jensen}, and \cite{CMIP}, and
of their forerunner, G\"odel.

\begin{center}
\begin{tabular}{l|r}
  Author(s) & $\lnot \ \Psi$ \\ \hline
  G\"odel & $\lnot \ 0^\#$ \\
  Dodd + Jensen & $\lnot \ 0^\dagger$ \\
  Koepke & $\lnot \ 0^{long}$ \\
  Jensen & $\lnot \ 0^{sword}$ \\
  Mitchell &  $\lnot \ o(\kappa) = \kappa^{++}$ \\
  Jensen & $\lnot \ 0^\P$ \\
  Steel & $\lnot \ M_1^\#$ \\
\end{tabular}
\end{center}

Except for the last entry this means that the person(s) listed on the left
hand side has (have)
shown $K$ to exist in the theory ${\sf ZFC} \ + \ \lnot \ \Psi$. (The definitions of the
respective anti
large cardinal assertions can be found in the above cited papers.)
It must however be mentioned that it was Jensen who proved 
the Covering Lemma for $L$ and who also developed the fine structure theory 
for $L$.

What Steel does in
\cite{CMIP} is to prove the existence of $K$ in the theory $${\sf ZFC} \ + \ \lnot \
M_1^\#
\ + \ {\rm there \ is \ a \ measurable \ cardinal,}$$ where by $\lnot \ M_1^\#$
we mean that there is no sharp
for an inner model with a Woodin cardinal.
For certain applications of $K$ this is a deficit (see the discussion in \S 0 of
\cite{FSIT}): Core model theory is typically applied as follows. 
If $\Phi$ is any
statement whatsoever, then one can try to lead $${\sf ZFC} + \ \lnot \ \Psi \ + \ \Phi$$ to
a contradiction by using $${\sf ZFC} + \ \lnot \ \Psi \ \Rightarrow \ K {\rm \ exists}$$
and showing that $\Phi$ contradicts the existence of $K$ ``below $\Psi$.''
If this is the case, one has arrived at finding a lower bound in terms of the
consistency strength of $\Phi$, that is, at proving
$${\sf ZFC} \ + \ \Phi \ \Rightarrow \ \Psi.$$
Consequently, from \cite{CMIP} we can (often, not always) only hope to get a theorem
of the form
$${\sf ZFC} \ + \ \Phi \ + \ {\rm there \ is \ a \ measurable \ cardinal }
\ \Rightarrow$$ $${\rm there \ is \ a \ sharp \ for \ an \ inner \ model \ with
\ a \ Woodin \ cardinal.}$$
(It should be noted that there is an important variation of this use of core model
theory. Instead of trying to lead $${\sf ZFC} + \ \lnot \ \Psi \ + \ \Phi$$ to an
outright
contradiction one can try to show that ${\sf ZFC} + \Phi$ implies that such-and-such
large cardinals exist in the $K$
``below $\Psi$.'')

At the time of writing
it is not known how to develop the theory of $K$ just assuming ``${\sf ZFC}
+$ there is no inner model with a Woodin cardinal.''
The present paper
solves the puzzle of constructing the core model for many strong cardinals.
Specifically, we improve
Jensen's notes \cite{jensen} by establishing that $K$ can
be shown to exist in the theory
\begin{center}
${\sf ZFC} +$ there is no sharp for an inner model with

a proper class of strong cardinals,
\end{center}
or $${\sf ZFC} +
\zerohandgrenade {\rm \ (zero \ hand-grenade)
\ does \ not \ exist,}$$
as we shall say.
It has turned out that the means by which \cite{jensen} can be improved in this
direction are
far from being straightforward generalizations of the means provided by \cite{jensen}.
We also want to emphasize that the work done here yields applications
which could not be obtained before (cf. section 9).

The main new idea here is that if $\zerohandgrenade$ does not exist
then normal iteration trees are simple enough
("almost linear") so that
$K^c$ can (still) be built by just requiring new extenders in its
recursive definition to be
countably complete (rather than having certificates).
This will allow a proof of weak covering for $K^c$, and
finally a proof of the existence of $K$ in the theory ${\sf ZFC} + \lnot \ \zerohandgrenade$.
There are serious obstacles to proving the existence of $K$ in a theory weaker than
${\sf ZFC} + \lnot \ \zerohandgrenade$.

We believe that the present paper raises (and answers) questions and develops
techniques which are interesting
from the point of view of
inner model theory, and which should be useful also outside the
realm of this paper.

The paper is organized as follows.
Section 1 discusses fine structure. In particular, we shall present our notation and
state what we'll have to assume familiarity with.
Section 2 introduces $\zerohandgrenade$, and proves that any normal iteration tree
of a
premouse which is below $\zerohandgrenade$ is ``almost linear.'' Section 3 builds
$K^c$, the countably complete core model below $\zerohandgrenade$. $K^c$ is a
preliminary version of $K$; it is not known how to develop the theory of $K$ without
going through $K^c$ and proving weak covering for it.
Section 4 proves that $K^c$ is maximal in the mouse order, and it shows a ``goodness''
property of $K^c$.
Sections 5 and 6 establish technical tools for proving the ``weak covering lemma for
$K^c$'' in section 7: we have to study bicephali which are somewhat more liberal than
any bicephalus studied so far, and we have to show a ``maximality'' property of $K^c$.
Section 8
finally isolates $K$.

Section 9
contains an application, and motivation, for this work. We'll show that the results of
the present paper provide the last step for determining the exact consistency strength
of the assumption in the statement of the $12^{{\rm th}}$ Delfino problem
(cf. \cite{KeMaSt}).
Section 9 also contains
a brief description of that problem.

\bigskip
{\em Historical note.} The work described in sections 2, 3, 4, 8, and 9 was done in
Berkeley in the fall of 1997 (parts of section 8 are due to John Steel). The results in
sections 5, 6, and 7 were obtained in Vienna in the fall and winter of 1999/2000.

\section{Preliminaries.}

We assume familiarity to a certain extent with Jensen's classical fine structure
theory
(cf. \cite{classic}, or
\cite[Chapter 1]{zeman}), with inner model theory as presented in \cite{NFS} or
\cite{FSIT}, and with core model theory as developed in
\cite{MOZ}, \cite{jensen}, \cite{CMIP}, or \cite{zeman}.
The policy of the present paper is that we have elaborated only on new ideas. In most
cases we won't give proofs for things which can be found elsewhere.

This paper builds upon the concept of a ``(Friedman-Jensen)
premouse'' proposed in \cite{NFS}.
Let us first briefly recall essentials from 
\cite[\S \S 1 and 4]{NFS}.
We explicitly warn the reader that this part of the current section of our present
paper
does not give any complete definitions. The reader may find any missing details in the
monograph \cite{zeman}.

An extender is a partial map $F \colon {\cal P}(\kappa) \rightarrow {\cal P}(\lambda)$
for some ordinals $\kappa < \lambda$ (cf. \cite[\S 1 p.2]{NFS})
rather than a system of hyper-measures.
(We shall sometimes abuse the notation by writing $F$ for the partial map
${\cal P}([\kappa]^{<\omega}) \rightarrow {\cal P}([\lambda]^{<\omega})$
induced by $F$ via soft coding [cf. \cite[\S 1 p.2]{NFS}];
and we'll often write $i_F$
for an ultrapower map ${\cal M} \rightarrow {\cal N}$ induced by $F$.)
Premice will be $J$-structures constructed from
certain well-behaved extender sequences.

A pre-premouse (cf. \cite[\S 4 p.1]{NFS})
is an acceptable $J$-structure of the form ${\cal M} =
( J_\alpha[{\vec E}] ; \in , {\vec E} , E_\alpha )$,  where $${\vec E} =
\{ ( \nu , \xi , X ) \ \colon \ \xi < \nu < \alpha \
{\rm and } \ \xi \in E_\nu(X) \}$$ codes a sequence of extenders,
with
the following two properties:

(a) If $E_\nu \not= \emptyset$ for $\nu \leq \alpha$, then
$E_\nu$ is an extender
whose domain  is ${\cal P}(\kappa) \cap J_\nu[{\vec E}]$ for some $\kappa < \nu$,
$(J_\nu[{\vec E} {\upharpoonright} \nu];\in,{\vec E} {\upharpoonright} \nu)$ is the
($\Sigma_0$-) ultrapower of
$(J_{\kappa^+}[{\vec E} {\upharpoonright} {\kappa^+}];\in,{\vec E} {\upharpoonright}
{\kappa^+})$ where $\kappa^+$ is calculated in $J_\nu[{\vec E} {\upharpoonright} \nu]$,
and if
$$i_{E_\nu} \ \colon \  ( J_\nu[{\vec E} {\upharpoonright} \nu] ;
\in , {\vec E} {\upharpoonright} \nu) \rightarrow_{E_\nu}
N$$
is the ($\Sigma_0$-) ultrapower map then
${\vec E}^N {\upharpoonright} \nu = {\vec E} {\upharpoonright} \nu$ and
$E^N_\nu = \emptyset$, and

(b) Proper initial segments of ${\cal M}$ are sound.

Condition (a) is often referred to as ``coherence.'' The ultrapower
formed in (a) is according to the upward extension of embeddings
technique using $E_\nu$ as a fragment of
the $i_{E_\nu}$ to be formed.
We always suppose the well-founded part of a model to be transitive.
(In our formulation of (a), we suppose that $\nu+1$ is a subset
of the well-founded part.)
The concept of "soundness" in (b) is according to Jensen's fine
structure (see below).

For pre-premice ${\cal M} =
( J_\alpha[{\vec E}] ; \in , {\vec E} , E_\alpha )$ as above and for $\nu \leq
\alpha$ we adopt the notation of
\cite[Definition 5.0.4]{FSIT} and write
${\cal J}^{\cal M}_\nu$ to mean $
( J_\nu[{\vec E} {\upharpoonright} \nu] ; \in , {\vec E}  {\upharpoonright} \nu, E_\nu )$.
Sometimes we shall confuse ${\cal J}^{\cal M}_\nu$ with its underlying universe,
$J_\nu[{\vec E} {\upharpoonright} \nu]$.

Let ${\cal M}$ be a pre-premouse as above. Let $F = E_\nu \not=
\emptyset$ be an extender with critical point $\kappa$
for some $\kappa < \nu$ and
$\nu \leq {\cal M} \cap {\rm OR}$. For any $\lambda \leq F(\kappa)$ we define $F|\lambda$ by
setting $(X,y) \in F|\lambda$ if and only if 
$\exists Y ( (X,Y) \in F \wedge \sigma(y)=Y )$ where
$\sigma$ the inverse of the collapse of the $\Sigma_1$ hull of $\lambda$
inside ${\rm Ult}({\cal J}^{\cal M}_\nu;F)$.
Let $C^M_\nu$ be the set of all $\lambda \in (\kappa,F(\kappa))$ such
that $F|\lambda$ is its own trivial completion, i.e., such that
$$\alpha < \lambda \wedge f \in {}^\kappa\kappa \cap {\cal J}^{\cal M}_\nu
\Rightarrow i_F(f)(\alpha) < \lambda.$$ (For such $\lambda$ we'll have that
$(F|\lambda)(X) = F(X) \cap \lambda$ for $X  \in {\rm dom}(F)$.)
The pre-premouse ${\cal M}$ is now called a premouse (cf. 
\cite[\S 4 p.2]{NFS})
provided we always have
$$\lambda \in C^M_\nu \Rightarrow F|\lambda \in {\cal J}^{\cal M}_\nu.$$
This clause is called the initial segment condition (cf. \cite[I]{CAR},
which gives the correction to the initial segment condition 
proposed in \cite{NFS}).

If ${\cal M}$ is a premouse with $C^{\cal M}_\nu \not= \emptyset$ for some $\nu \leq
{\cal M} \cap {\rm OR}$ then ${\cal M}$ is not ``below superstrong,'' in that ${\cal M}$ has
then an extender on its sequence witnessing that ${\cal J}^{\cal M}_\nu \models$
``there is a superstrong cardinal.'' That is, below superstrongs does the initial
segment condition collapse to the requirement that we have
$C^{\cal M}_\nu = \emptyset$ for all $\nu \leq
{\cal M} \cap {\rm OR}$.

\bigskip
Inner model theory has to iterate premice. Fine iterations in \cite{NFS} 
are based on
Jensen's smooth $\Sigma^*$-theory.
Jensen has developed a machinery for taking fine ultrapowers
$\pi \colon {\cal M} \rightarrow_F^\star {\cal N}$ in such a way that
$\pi$ will be what he calls $\Sigma^{(n)}_0$-elementary for all $n<\omega$
with $\rho_n({\cal M}) > {\rm c.p.}(\pi)$.
As presented in \cite{FSIT},
in order to develop the theory of $K$ one would only need
$\Sigma^{(n)}_0$-elementarity here for all those $n<\omega$ such that ${\cal M}$ is
$n$-sound. When following this latter route and
defining the concept of a fine ultrapower \`a la \cite{FSIT},
one can moreover in fact stick to
the more traditional master codes, that is to
``coding ${\cal M}$ onto $\rho_n({\cal M})$, taking a
$\Sigma_0$-ultrapower of the coded structure, and then decoding'' (\cite{FSIT} p. 40;
see also the discussion in the introduction to
\cite{NFS}). Being pragmatic, it is this latter approach which we shall follow here.
It has a couple of advantages:
it suffices for our purposes, and it
will simplify our iterability proof (cf.
\ref{iterability} below) and will make it accessible for people ignorant
of Jensen's $\Sigma^*$-theory.

\bigskip
Let ${\cal M}$ be an acceptable $J$-structure. We shall 
write (cf. \cite[Chapter 1]{zeman}):

\bigskip
$\bullet \ \ $ $\rho_n({\cal M})$ $ \ \ $ for the $n^{{\rm th}}$ projectum of ${\cal M}$,

$\bullet \ \ $ $P^n_{\cal M}$ $ \ \ $ for the set of good parameters (i.e., for the set
of parameters witnessing $\rho_n({\cal M})$ is the $n^{{\rm th}}$ projectum),

$\bullet \ \ $ $p_{{\cal M},n}$ $ \ \ $ for the $n^{{\rm th}}$
standard parameter of ${\cal M}$
(i.e., the least element of $P^n_{\cal M} \cap [{\rm OR}]^{< \omega}$
where ``least'' refers to the
canonical well-order of $[{\rm OR}]^{<\omega}$),

$\bullet \ \ $ $W^\nu_{\cal M}$ $ \ \ $ for the 
witness for $\nu \in p_{{\cal M},n}$,

$\bullet \ \ $ $A^{n,p}_{\cal M}$ $ \ \ $ for the $n^{{\rm th}}$
standard code determined by $p$,

$\bullet \ \ $ ${\cal M}^{n,p}$ $ \ \ $ for the $n^{{\rm th}}$ reduct determined by $p$,

$\bullet \ \ $ ${\cal M}^n$ $ \ \ $ for ${\cal M}^{n,p_{{\cal M},n}}$,

$\bullet \ \ $ $h^{n,p}_{\cal M}$ ($h^n_{\cal M}$) $ \ \ $
for the canonical $\Sigma_1$ Skolem function of
${\cal M}^{n,p}$ (of
${\cal M}^n$),

$\bullet \ \ $ $R^n_{\cal M}$ $ \ \ $ for the 
set of very good parameters (i.e.,
for the set of $p$ such that ${\cal M}^{n-1,p}$ is generated by
$h^{n-1,p}_{\cal M}$ from $\rho_n({\cal M}) \cup \{ p \}$), and

$\bullet \ \ $ ${\frak C}_n({\cal M})$  $ \ \ $ for 
the $n$-core of ${\cal M}$ (i.e., for that
${\bar {\cal M}}$ such that ${\bar {\cal M}}^k$ is the transitive collapse of
what is generated by
$h^k_{\cal M}$ from $\rho_{k+1}({\cal M}) \cup \{ p_{{\cal
M}^k,1} \}$ for all $k < n$).

\bigskip
By definition, ${\cal M}$ is $n$-sound if and only if $P^n_{\cal M} = R^n_{\cal M}$
(we thereby follow \cite{zeman} and
emphasize that this is in contrast to \cite[Defition 2.8.3]{FSIT}).
As a matter of fact,
${\cal M}$ is $n$-sound if and only if $p_{{\cal M},k} \in R^k_{\cal M}$ for all $k \leq n$
if and only if ${\cal M} = {\frak C}_n({\cal M})$.
Moreover, if ${\cal M}$ is $n$-sound then $p_{{\cal M},k} = p_{{\cal M},n}
{\upharpoonright} k$ for all $k \leq n$.

Let ${\cal M}$ and ${\cal N}$ be acceptable $J$-structures, and let $n \leq \omega$.
We call $\pi \colon {\cal M} \rightarrow {\cal N}$ an $n$-embedding (cf.
\cite[Definition 2.8.4]{FSIT}) if and only if

\bigskip
$\bullet \ \ $ ${\cal M}$ and ${\cal N}$ are both $n$-sound,

$\bullet \ \ $ $\pi(p_{{\cal M},n}) = p_{{\cal N},n}$,

$\bullet \ \ $ $\pi(\rho_k({\cal M})) = \rho_k({\cal N})$ for all $k < n$
(by convention, $\pi({\cal M} \cap {\rm OR}) = {\cal N} \cap {\rm OR}$),

$\bullet \ \ $ $\pi$"$\rho_n({\cal M})$ is cofinal in $\rho_n({\cal N})$, and

$\bullet \ \ $ $\pi {\upharpoonright} {\cal M}^k \colon
{\cal M}^k \rightarrow_{\Sigma_1} {\cal N}^k$ for all $k \leq n$.

\bigskip
The significance of $n$-embeddings is that they
are generated by taking $n$-ultra-powers.
Let ${\cal M}$ be an $n$-sound premouse where $n \leq \omega$.
Suppose that $F$ is an extender which measures all the subsets of its
critical point which are in ${\cal M}^{n}$, and let $$i_F \colon {\cal M}^n
\rightarrow_F {\bar {\cal N}}$$
be the $\Sigma_0$-ultrapower map. There will be at most one
transitive structure ${\cal N}$ together with a $\Sigma_1$-elementary embedding
$\pi \colon {\cal M} \rightarrow {\cal N}$ such that
$\pi(p_{{\cal M},n})
\in R^n_{\cal M}$, ${\bar {\cal N}} = {\cal N}^{n,\pi(p_{{\cal M},n})}$, and
$\pi {\upharpoonright} {\cal M}^n = i_F$. This follows
from the upward extension of embeddings lemma (cf. 
\cite[Chapter 1]{zeman}).
By forming a term model, one can always produce a (not necessarily well-founded)
unique (up to isomorphism)
candidate for such an ${\cal N}$. We shall denote this situation by
$$\pi \colon {\cal M} \rightarrow_F^n {\cal N} = {\rm Ult}_n({\cal M};F) {\rm , }$$
and call it the $n$-ultrapower of ${\cal M}$ by $F$. As a matter of fact, if $F$ is
``close to'' ${\cal M}$ (cf. \cite[Definition 4.4.1]{FSIT}) then we'll
have that $\rho_{n+1}({\cal M}) = \rho_{n+1}({\cal N})$.

Now if $\pi \colon {\cal M} \rightarrow_F^n {\cal N}$ and ${\cal N}$ is transitive
then $\pi$ is an $n$-embedding if and only if
$\pi(p_{{\cal M},n}) = p_{{\cal N},n}$ and
$\pi(\rho_k({\cal M})) = \rho_k({\cal N})$ for all $k < n$.
This in turn will follow from a certain solidity and universality of $p_{{\cal
M},n}$.
We call
${\cal M}$ $n$-solid just in case that for all $\nu \in p_{{\cal M},n}$ we have that
$W^\nu_{\cal M} \in {\cal M}$.
We call ${\cal M}$ $n$-universal if and only if $${\cal P}(\rho_n({\cal M})) \cap {\cal M}
= \{
h^n_{\cal M}(i,({\vec \alpha},p_{{\cal M},n})) \cap \rho_n({\cal M}) \ \colon \
i < \omega \wedge {\vec \alpha} \in \rho_n({\cal M}) \}.$$
We have that $\pi(p_{{\cal M},n}) = p_{{\cal N},n}$ and
$\pi(\rho_k({\cal M})) = \rho_k({\cal N})$ for all $k < n$ hold true in the above
setting if ${\cal M}$ is $n$-solid and $n$-universal.
One of the key theorems of fine structure theory (cf. \cite[\S 8]{CMIP}, or
\cite[\S 7]{NFS}) will tell us that if ${\cal M}$ is $n$-sound and
$n$-iterable then ${\cal M}$ is $(n+1)$-solid and $(n+1)$-universal.
Here, ``$n$-iterability'' has to be explained.

Before turning to that it remains to introduce a weakened version of $n$-embeddings
which will come up in the proofs of \ref{normal-iterability}
and \ref{it-of-bicephalus}.
Let ${\cal M}$ and ${\cal N}$ be acceptable $J$-structures, and let $n \leq \omega$.
We call $\pi \colon {\cal M} \rightarrow {\cal N}$ a weak $n$-embedding (cf.
\cite[p. 52 ff.]{FSIT}) if and only if

\bigskip
$\bullet \ \ $ ${\cal M}$ and ${\cal N}$ are both $n$-sound,

$\bullet \ \ $ $\pi(p_{{\cal M},n}) = p_{{\cal N},n}$,

$\bullet \ \ $ $\pi(\rho_k({\cal M})) = \rho_k({\cal N})$ for all $k < n$
(by convention, $\pi({\cal M} \cap {\rm OR}) = {\cal N} \cap {\rm OR}$),

$\bullet \ \ $ ${\rm sup}(\pi$"$\rho_n({\cal M})) \leq \rho_n({\cal N})$,

$\bullet \ \ $ $\pi {\upharpoonright} {\cal M}^k \colon
{\cal M}^k \rightarrow_{\Sigma_1} {\cal N}^k$ for all $k < n$, and

$\bullet \ \ $ there is a set $X \subset {\cal M}$ with $\{ p_{n,{\cal M}} ,
\rho_n({\cal M}) \} \subset X$, $X \cap \rho_n({\cal M})$ is cofinal in
$\rho_n({\cal M})$, and $\pi {\upharpoonright} {\cal M}^n \colon
{\cal M}^n \rightarrow {\cal N}^n$
is $\Sigma_1$-elementary on parameters from $X$.

\bigskip
There is a generalization of taking an ultrapower by an extender, namely,
taking an ultrapower by a ``long extender'' (cf. for example 
\cite[\S 2.5]{MiSchSt}).
If ${\cal M}$ is a premouse, if $\nu$ is a regular cardinal of ${\cal M}$, and if
$\pi \colon {\cal J}^{\cal M}_\nu \rightarrow {\cal N}$ is an embedding
then we shall write
$${\rm Ult}_n({\cal M};\pi)$$ for the $n$-ultrapower of ${\cal M}$ by the long extender
derived from $\pi$ (see \cite[\S 2.5]{MiSchSt}),
also called the lift up of ${\cal N}$ by $\pi$.
We'll use this technique in sections 6
and 7.

\bigskip
Let us now discuss iterations of premice.
Premice are iterated by forming iteration trees. When ${\cal M}$ is a premouse
then a typical iteration tree on ${\cal M}$ will have the form
$${\cal T} = (({\cal M}^{\cal T}_\alpha,\pi^{\cal T}_{\alpha \beta} \colon \alpha
T
\beta < \theta),(E_\alpha^{\cal T} \colon \alpha+1 < \theta),T)$$ where $T$ is the
tree structure, ${\cal M}^{\cal T}_\alpha$ are the models, $\pi^{\cal T}_{\alpha
\beta}$ are the embeddings, and $E_\alpha^{\cal T}$ are the extenders used. When
dealing with iteration trees we shall use the notation from \cite{FSIT}.
We call an iteration tree ${\cal T}$ as above normal (cf. 
\cite[\S 4 p. 4]{NFS})
if and only if the following three
requirements are met.

\bigskip
$\bullet \ \ $ The index of $E_\beta^{\cal T}$ is (strictly) greater than the index
of $E_\alpha^{\cal T}$ for all $\alpha < \beta < \theta$,

$\bullet \ \ $ $T$-pred$(\alpha+1) =$ the least $\beta \leq \alpha$ such that
$${\rm c.p.}(E_\alpha^{\cal T}) < E_\beta^{\cal T}({\rm c.p.}(E_\beta^{\cal T})) {\rm , \ and }$$

$\bullet \ \ $ setting $\alpha^\star = T$-pred$(\alpha+1)$,
${\rm dom}(\pi^{\cal T}_{\alpha^\star \alpha+1}) = {\cal J}^{{\cal
M}^{\cal T}_{\alpha^\star}}_\eta$ where $\eta$ is largest such that
$E_\alpha^{\cal T}$ measures all the subsets of its critical point which are in
${\cal J}^{{\cal
M}^{\cal T}_{\alpha^\star}}_\eta$, for all $\alpha+1 < \theta$.

\bigskip
We shall build our fine iteration trees by forming $n$-ultrapowers.
That is, if $${\cal T} =
(({\cal M}^{\cal T}_\alpha,\pi^{\cal T}_{\alpha \beta} \colon \alpha T
\beta < \theta),(E_\alpha^{\cal T} \colon \alpha+1 < \theta),T)$$
is an iteration tree we'll have that for all $\alpha+1 < \theta$, setting
$\alpha^\star = T$-pred$(\alpha+1)$,
$$\pi^{\cal T}_{\alpha^\star \alpha} \colon {\rm dom}(\pi^{\cal T}_{\alpha^\star \alpha})
\rightarrow^n_{E_\alpha^{\cal T}} {\cal M}^{\cal T}_{\alpha+1}$$ for some $n \leq
\omega$. In this case, we shall denote $n$ by ${\rm deg}^{\cal T}(\alpha+1)$.
Notice that if ${\rm deg}^{\cal T}(\alpha+1) = n$ then
${\rm dom}(\pi^{\cal T}_{\alpha^\star \alpha})$ has to be $n$-sound, and
${\cal M}^{\cal T}_{\alpha+1}$ will be $n$-sound by design.
We call
an iteration tree ${\cal T}$ on a premouse ${\cal M}$
$n$-bounded if ${\rm deg}^{\cal T}(\alpha+1) \leq n$ for all $\alpha+1 < \theta$ such that
${\cal D}^{\cal T} \cap (0,\alpha+1]_T = \emptyset$, i.e., there is no drop along
$[0,\alpha+1]_T$.

We call ${\cal T}$ a putative iteration tree if ${\cal T}$ is an iteration tree except
for the fact that if ${\cal T}$ has a last model, ${\cal M}^{\cal T}_\infty$, then
${\cal M}^{\cal T}_\infty$ does not have to be well-founded.

Let ${\cal M}$ be a premouse. For our purposes, we may think of ${\cal M}$ as being
normally $n$-iterable if all putative
normal $n$-bounded iteration trees on ${\cal M}$ of successor length have the property
that ${\cal M}^{\cal T}_\infty$, the last model of ${\cal T}$, is transitive, and
${\cal D}^{\cal T} \cap (0,\infty]_T$ is finite.\footnote{The
informed reader might miss an assertion about the existence of well-founded branches
for trees of limit length; the reason for not including such a thing is that
by virtue of \ref{unique-branch} we'll only have to deal with trees with exactly one
cofinal branch.
In general we'd have to say that there exists an iteration strategy for ${\cal
M}$, i.e., a certain partial function picking cofinal well-founded branches through
certain trees on ${\cal M}$ of limit length.}
In a comparison process we'll typically use maximal trees.
For $n \leq \omega$
we call ${\cal T}$ $n$-maximal (cf. \cite[Definition 6.1.2]{CMIP})
if the following two requirements are met.

\bigskip
$\bullet \ \ $ if ${\cal D}^{\cal T} \cap (0,\alpha+1]_T = \emptyset$ then
${\rm deg}^{\cal T}(\alpha+1)$ is the largest $k \leq n$ such that ${\rm c.p.}(E_\alpha^{\cal T})
< \rho_k({\rm dom}(\pi_{{\alpha^*} \alpha+1}^{\cal T}))$ where $\alpha^* =
T$-pred$(\alpha+1)$, and

\bigskip
$\bullet \ \ $ if ${\cal D}^{\cal T} \cap (0,\alpha+1]_T \not= \emptyset$ then
${\rm deg}^{\cal T}(\alpha+1)$ is the largest $k \leq \omega$
such that ${\rm c.p.}(E_\alpha^{\cal T})
< \rho_k({\rm dom}(\pi_{{\alpha^*} \alpha+1}^{\cal T}))$ where $\alpha^* =
T$-pred$(\alpha+1)$.

\bigskip
That
${\cal M}$ is $n$-iterable (see \cite[Definition 5.1.4]{FSIT})
will say that every stack $${\cal T}_0 {}^\frown
{\cal T}_1 {}^\frown {\cal T}_2 {}^\frown ...$$ of normal iteration trees
is well-behaved, where ${\cal T}_{i+1}$ is a tree on an initial segment of
the last
model of ${\cal T}_i$. More precisely, ${\cal M}$ is $n$-iterable
provided the following holds true: {\em if} $0 < k \leq \omega$ and $({\cal T}_i
\colon i<k)$ is a sequence of {\it normal} iteration trees such that

\bigskip
$\bullet \ \ $ ${\cal M}^{{\cal T}_0}_0 = {\cal M}$,

$\bullet \ \ $ ${\cal T}_0$ is $n$-bounded,

$\bullet \ \ $ ${\cal T}_i$ has successor length $\theta_i$ for all $i<k$,

$\bullet \ \ $ $\exists \gamma_i \
{\cal M}^{{\cal T}_{i+1}}_0 = {\cal J}_{\gamma_i}^{{\cal M}^{{\cal T}_i}_{\theta_i}}$
for all
$i+1 < k$, and

$\bullet \ \ $ ${\cal T}_{i+1}$
is $n(i)$-bounded, where $n(i)$ is maximal such that
${\cal M}^{{\cal T}_{i+1}}_0$ is $n(i)$-sound and
$\forall j \leq i \ [ \gamma_j = {\cal M}^{{\cal T}_j}_{\theta_j} \cap
{\rm OR} \wedge {\cal D}^{{\cal T}_j} \cap (0,\theta_j]_{T_j} =
\emptyset ] \Rightarrow n(i) \leq n$, for all $i+1 < k$,

\bigskip
{\em then} either $k<\omega$ and ${\cal M}^{{\cal T}_{k-1}}_{\theta_{k-1}}$ is
normally ${\bar n}$-iterable,
where ${\bar n}$ is maximal such that
${\cal M}^{{\cal T}_{k-1}}_0$ is ${\bar n}$-sound and
$\forall j < k-1 \
[ \gamma_j = {\cal M}^{{\cal T}_j}_{\theta_j} \cap
{\rm OR} \wedge
{\cal D}^{{\cal T}_j} \cap (0,\theta_j] =
\emptyset ] \Rightarrow {\bar n} \leq n$, or
else $k=\omega$ and ${\cal D}^{{\cal T}_i} \cap
(0,\theta_i]_{T_i} = \emptyset$ for all but finitely many $i<\omega$ and the direct
limit of the ${\cal M}^{{\cal T}_i}_0$'s together with the obvious maps
(for sufficiently large $i$) is well-founded.

\bigskip
The following lemma, which we might call the ``strong initial segment condition''
for Friedman-Jensen mice,
is an immediate consequence of the
condensation lemma \cite[\S 8 Lemma 4]{NFS}. It is just a slight strengthening of
\cite[\S 8 Corollary 4.2]{NFS}.

\begin{lemma}\label{strong-ISC}
Let ${\cal M} = (J_\alpha[{\vec E}];\in,{\vec E},E_\alpha)$
be a $0$-iterable premouse,
where $E_\alpha \not= \emptyset$.
Suppose that for no $\mu \leq {\cal M} \cap {\rm OR}$ do we have ${\cal J}^{\cal M}_\mu
\models {\sf ZFC} +$ ``there is a Woodin cardinal.''
Set $\kappa = {\rm c.p.}(E_\alpha)$, and let $\xi \in
(\kappa,\rho_1({\cal M}))$. Then one of (a) or (b) below holds:

(a) There is some $\eta < \alpha$ such that $E_\alpha | \xi = E_\eta^{\cal M}$.

(b) $\xi$ is not a cardinal in ${\cal M}$, and
there are $\mu < \eta < {\cal M} \cap {\rm OR}$ and $k<\omega$
such that $E_\alpha | \xi$ is the top
extender of ${\rm Ult}_k({\cal J}^{\cal M}_\eta;E_\mu^{\cal M})$ and $(\eta,k)$ is
$<_{lex}$-least with $\eta \geq \mu$ and
$\rho_{k+1}({\cal J}_\eta^{\cal M}) \leq {\rm c.p.}(E_\mu^{\cal M})$.

Moreover, in any case there is some $E_{\tilde \nu}^{\cal M}$ with
$\xi < {\tilde \nu} < \xi^{+{\cal M}}$ and ${\rm c.p.}(E_{\tilde \nu}^{\cal M}) = \kappa$.
\end{lemma}

{\sc Proof.} As $E_\alpha | \xi = E_\alpha | (\xi \cup \kappa^{+{\cal M}})$ for $\xi >
\kappa$,
let us assume without loss of generality that $\xi \geq \kappa^{+{\cal M}}$. Consider
$$\sigma \colon {\bar {\cal M}} \cong H^{\cal M}(\xi) \prec_{\Sigma_1} {\cal M}
{\rm , }$$ where $H^{\cal M}(\xi)$ denotes the
hull generated from $\xi$ by $h^0_{\cal M}$, and ${\bar {\cal M}}$ is transitive.
Let $\nu \leq {\bar {\cal M}} \cap {\rm OR}$ be maximal with $\sigma
{\upharpoonright} \nu = {\rm id}$.
We have that $\rho_1({\bar {\cal M}}) \leq \xi$, and ${\bar {\cal M}}$
is trivially $1$-sound above $\nu \geq \xi$ (i.e., ${\bar {\cal M}}$
is the hull generated from  $\nu \cup \{ p_{{\bar {\cal M}},1} \}$
by $h^0_{\bar {\cal M}}$).
Because $\rho_1({\bar {\cal M}}) \leq \xi < \rho_1({\cal M})$, by
the condensation lemma \cite[\S 8 Lemma 4]{NFS} there remain just two
possibilities:

(i) ${\bar {\cal M}} = {\cal J}^{\cal M}_\eta$ for some $\eta < {\cal M} \cap {\rm OR}$, or

(ii) ${\bar {\cal M}} = {\rm Ult}_k({\cal J}^{\cal M}_\eta;E_\mu^{\cal M})$ where
$$\rho_\omega({\cal J}_\eta^{\cal M}) <
\nu = ({\rm c.p.}(E_\mu^{\cal M}))^{+{\cal J}_\mu^{\cal M}} =
({\rm c.p.}(E_\mu^{\cal M}))^{+{\cal J}_\eta^{\cal M}} < \mu \leq \eta < {\cal M}
\cap {\rm OR} {\rm , }$$ $E_\mu^{\cal M}$ is generated by its critical point, and
$k<\omega$ is such that $\rho_{k+1}({\cal J}_\eta^{\cal M}) < \nu \leq
\rho_{k}({\cal J}_\eta^{\cal M})$.

Now let $F$ denote the top extender of ${\bar {\cal M}}$. It is easy to see that
$F | \xi = E_\alpha | \xi$. Moreover, all the generators of $F$ are below $\xi$, so
that in fact $F = F | \xi = E_\alpha | \xi$.

Let us first assume that (i) holds.
Then $F = E^{\cal M}_\eta$, so that (a) in the statement of \ref{strong-ISC}
holds.

Let us then assume that (ii) holds.
Then $F$ has to be the top extender of ${\rm Ult}_k({\cal J}^{\cal M}_\eta;E_\mu^{\cal M})$.
It is also easy to see that $\eta > \mu$.
Set ${\bar \mu} = {\rm c.p.}(E_\mu^{\cal M})$. As
$\rho_\omega({\cal J}_\eta^{\cal M}) \leq {\bar \mu}$ and $\kappa^{+{\bar {\cal M}}} =
\kappa^{+{\cal M}}$, we must have ${\bar \mu} > \kappa$. Hence by (ii),
${\cal J}^{\cal M}_\eta$ has to have a top extender with critical point $\kappa$,
namely $E^{\cal M}_\eta$. Notice that $\eta > \nu \geq \xi$.

Now suppose that
$\xi$ is a cardinal in ${\cal
M}$. We'll then have that $\rho_\omega({\cal J}_\eta^{\cal M}) \leq {\bar \mu}$
and $\eta < {\cal M} \cap {\rm OR}$ imply that $\xi \leq
{\bar \mu}$. But taking the ultrapower ${\rm Ult}_k({\cal J}^{\cal M}_\eta;E_\mu^{\cal M})$
is supposed to give ${\bar {\cal M}}$ and it
adds ${\bar \mu}$ as a generator of its top extender, $F$. However,
all the generators of $F$ are below $\xi$. Contradiction!

We have shown that if (ii) holds then (b) in the statement of \ref{strong-ISC}
holds.

Finally, the careful reader will have observed that
we also have established the ``moreover'' clause in the statement of \ref{strong-ISC}:
notice that we'll always have that ${\bar {\cal M}} \in {\cal M}$, and that
${\bar {\cal M}}$ has the same cardinality as $\xi$ inside ${\cal M}$.

\bigskip
\hfill $\square$ (\ref{strong-ISC})

\bigskip
\ref{strong-ISC} should be compared with the initial 
segment condition of \cite[Definition 1.0.4 (5)]{FSIT}
(cf. also \cite[Definition 2.4]{ssz}). However, as
Friedman-Jensen extenders allow
more room between the sup
of their generators and their index, one has to add the hypothesis
that $\xi < \rho_1({\cal M})$
in \ref{strong-ISC}.

\begin{defn}\label{defn-tau-strong}
Let ${\cal M}$ be a premouse, and let $\kappa <
\tau \leq {\cal M} \cap {\rm OR}$.
Then $\kappa$ is said to be $<\tau$-strong in ${\cal M}$ if for all $\alpha < \tau$
there is some extender $F \in {\cal M}$ with ${\rm dom}(F) = {\cal P}(\kappa) \cap {\cal M}$,
${\rm c.p.}(F) = \kappa$,
and ${\cal J}^{\cal M}_\alpha
\in {\rm Ult}_0({\cal M};F)$.
Furthermore,
$\kappa$ is said to be $<\tau$-strong in ${\cal M}$ as witnessed by ${\vec E}^{\cal
M}$ if for all $\alpha < \tau$
there is some $E^{\cal M}_\beta \not= \emptyset$ with ${\rm c.p.}(E^{\cal M}_\beta) =
\kappa$ and $\alpha < \beta < \tau$.
\end{defn}

We may now phrase an immediate corollary to \ref{strong-ISC} as follows.

\begin{cor}\label{strong-ISC+}
Let ${\cal M} = (J_\alpha[{\vec E}];\in,{\vec E},E_\alpha)$
be a $0$-iterable premouse,
where $E_\alpha \not= \emptyset$.
Suppose that for no $\mu \leq {\cal M} \cap {\rm OR}$ do we have ${\cal J}^{\cal M}_\mu
\models {\sf ZFC} +$ ``there is a Woodin cardinal.''
Set $\kappa = {\rm c.p.}(E_\alpha)$, and let $\tau \in (\kappa^{+{\cal M}},\rho_1({\cal
M})]$ be a cardinal in ${\cal M}$. Then $\kappa$ is $<\tau$-strong in ${\cal M}$
as witnessed by ${\vec E}^{\cal
M}$.
\end{cor}

%
%
%

\section{Almost linear iterations and $\zerohandgrenade$.}

In this section we show that for premice which do not encompass
a measurable limit of strong
cardinals (i.e., which are ``below $\zerohandgrenade$,'' as we shall say)
almost linear iterability -- in a sense to be made precise --
suffices for comparison.
It appears that T. Dodd already knew that such small mice can be linearily compared,
at least the mice of his time (see \cite{dodd}).

\begin{defn}\label{almost-linear}
An iteration tree ${\cal T}$ is called almost 
linear\footnote{\cite[Definition 1.10]{john-philip} introduces a different concept of
almost linearity. If a Mitchell-Steel premouse ${\cal M}$ is below $0^\P$ then
every iteration tree on ${\cal M}$ is almost linear in the 
sense of \cite[Definition 1.10]{john-philip}.}
if the following
hold true.

(a) For all $i+1 < {\rm lh}({\cal T})$
we have that ${\cal T}$-${\rm pred}(i+1) \in [0,i]_T$, and

(b) any $i < {\rm lh}({\cal T})$ only has finitely
many immediate ${\cal T}$-successors.
\end{defn}

\[
\xymatrix{
  & \omega &   & 4 &   &   &   & \\
\omega+1  & & &   & 3 \ar[lu] &   & 2 & \\
  &  &   & 5 \ar@{.>}[lluu] \ar[lllu] &  & 1 \ar[lu] \ar[ru] & & \\
  &  &   &   & 0 \ar[lu] \ar[ru] &   &   & \\
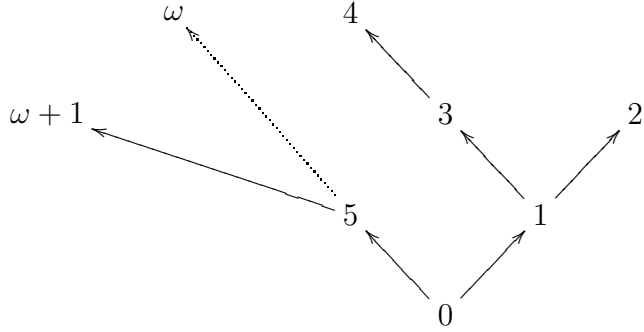
  &  &   & {}\save[]*{\txt{{\sc Figure 1.} An almost linear iteration tree}} \restore & & & & \\
}
\]

\bigskip
We shall not have to worry about finding
cofinal branches for almost linear iteration trees of limit length.

\begin{lemma}\label{unique-branch}
Let ${\cal T}$ be an almost linear iteration
tree. Let $\lambda \leq {\rm lh}({\cal T})$ be a limit ordinal (or
$\lambda = {\rm OR}$). Then ${\cal T} {\upharpoonright} \lambda$ has a unique
cofinal branch $b$, and if $\lambda
< {\rm lh}({\cal T})$ then $b = [0,\lambda)_T$.
\end{lemma}

{\sc Proof.} Obviously, the second part follows from the first.
We get $b$ by the following simple recursion.
If $i \in b$, $i+1 < \lambda$, then
the immediate ${\cal T}$-successor of $i$ in $b$ is the maximal (as an
ordinal) $i' > i$ being immediate ${\cal T}$-successor of $i$. And
if ${\bar \lambda} < \lambda$ is a limit ordinal and cofinally many
$i < {\bar \lambda}$ are in $b$ then ${\bar \lambda}$ is in $b$, too.
\ref{almost-linear} is just what is needed to see that this works.

\bigskip
\hfill $\square$ (\ref{unique-branch})

\begin{defn}\label{0-handgrenade}
Let ${\cal M}$ be a premouse. ${\cal M}$ is said to be
below $\zerohandgrenade$ (pronounced "zero hand-grenade")
if 
there is no $\kappa$ which is the critical point of an extender
$E^{\cal M}_\nu \not= \emptyset$
with the property that
$$\{ \mu < \kappa \ \colon \ \mu {\rm \ is \ } 
<\kappa{\rm -strong \ in \ } {\cal
M} \} {\rm \ is \ unbounded \ in \ } \kappa.$$
\end{defn}

We admit that ${}^\bullet \! \! \! \! \mid$ doesn't seem to
look like a hand-grenade.
The mice of inner model theory don't resemble the mice in our backyard
as well. ``Hand-grenade'' is just another math term 
in the tradition of daggers,
swords, and pistols.

We want to emphasize that if the premouse ${\cal M}$ is below $\zerohandgrenade$ then
${\cal M}$ is of course below superstrong as well, so that we'll have that $C^{\cal
M}_\nu = \emptyset$ for all $\nu \leq {\cal M} \cap {\rm OR}$.

The following lemma shows the benefit of life without hand-grenades.

\begin{lemma}\label{normal-linear}
Let ${\cal M}$ be a premouse which is below $\zerohandgrenade$.
Then any normal iteration of
${\cal M}$ is almost linear.
\end{lemma}

{\sc Proof.} Let us fix ${\cal M}$. We first aim to prove:

\bigskip
{\bf Claim 1.} For all normal ${\cal T}$ on ${\cal M}$ of successor length
$\vartheta+1$ do we have the following. Let $i+1 \in (0,\vartheta]_T$,
and set $\kappa = {\rm c.p.}(E^{\cal T}_i)$, and $\lambda = E^{\cal T}_i(\kappa)$. 
Then there is no
$G = E^{{\cal M}^{\cal T}_\vartheta}_\nu \not= \emptyset$ with
${\rm c.p.}(G) \in [\kappa,\lambda)$ and $\nu \geq \lambda$.

\bigskip
{\sc Proof.}
Let ${\cal T}$ be a normal iteration tree on ${\cal M}$ of length $\vartheta+1$,
and let us
consider $$\pi^{\cal T}_{i^* i+1} \ \colon {\cal J}^{{\cal M}^{\cal
T}_{i^*}}_{\eta} \rightarrow_F {\cal M}^{\cal T}_{i+1} {\rm , }$$
where $F = E^{\cal T}_i$, $i^* =
T$-${\rm pred}(i+1)$, ${\cal J}^{{\cal M}^{\cal
T}_{i^*}}_{\eta} = {\rm dom}(\pi^{\cal T}_{i^* i+1})$, and $i+1 T \vartheta$. We put
$\kappa = {\rm c.p.}(F)$, and $\lambda = F(\kappa)$. Suppose that there is some
$G = E^{{\cal M}^{\cal T}_\vartheta}_\nu \not= \emptyset$ with
${\rm c.p.}(G) \in [\kappa,\lambda)$ and $\nu \geq \lambda$.
Let
$\mu = {\rm c.p.}(G)$.

\bigskip
{\bf Subclaim.} $\rho_\omega({\cal J}^{{\cal M}^{\cal T}_\vartheta}_\nu) \geq \lambda$.

\bigskip
{\sc Proof.}
Of course, $\lambda$ is a cardinal of all ${\cal M}^{\cal T}_k$
with $k \geq i+1$. This trivially implies
$\rho_\omega({\cal J}^{{\cal M}^{\cal T}_\vartheta}_\nu) \geq \lambda$ if ${\cal
D}^{\cal T} \cap (i+1,\vartheta]_T \not= \emptyset$ or if $G$ is not the top extender
of ${\cal M}^{\cal T}_\vartheta$.

Let us suppose that ${\cal
D}^{\cal T} \cap (i+1,\vartheta]_T = \emptyset$ and that $G$ is the top extender
of ${\cal M}^{\cal T}_\vartheta$. By the normality of
${\cal T}$ we will have that $\pi_{i+1 \vartheta}^{\cal T} {\upharpoonright} \lambda =
{\rm id}$. But $\pi_{i+1 \vartheta}^{\cal T} \colon {\cal M}_{i+1}^{\cal T}
\rightarrow {\cal M}_{\vartheta}^{\cal T}$
is sufficiently elementary to yield
that then ${\cal M}_{i+1}^{\cal T}$
must have a top extender with critical point $\mu = {\rm c.p.}(G)$. In particular,
$\mu \in {\rm ran}(\pi^{\cal T}_{i^* i+1})$. However, $\mu \in [\kappa,\lambda)$ clearly
implies
$\mu \notin {\rm ran}(\pi^{\cal T}_{i^* i+1})$. Contradiction!

\bigskip
\hfill $\square$ (Subclaim)

\bigskip
Now by the Subclaim we get that
$G | \xi \in {\cal J}^{{\cal M}^{\cal T}_\vartheta}_\nu$ for all $\xi \in
(\mu,\lambda)$. Thus $$( \ G | \xi \ \colon \ \xi \in
(\mu,\lambda))$$ witnesses that
$\mu$ is $<\lambda$-strong in ${\cal M}^{\cal T}_\vartheta$.
But as
${\cal J}^{{\cal M}^{\cal T}_\vartheta}_{\lambda} =
{\cal J}^{{\cal M}^{\cal T}_{i+1}}_{\lambda}$ and $\lambda$ is a cardinal in both of
these models,
this says that
$\mu$ is $<\lambda$-strong in ${\cal M}^{\cal T}_{i+1}$.

Let $\alpha < \kappa$ be arbitrary. We have seen that
$${\cal M}_{i+1}^{\cal T}
\models {\rm \ `` }\exists {\tilde \mu} \in (\alpha,\lambda) \ ({\tilde \mu}
{\rm \ is \ } <\lambda {\rm -strong}){\rm ."}$$ As $\pi_{i^* i+1}^{\cal T}$ is
sufficiently elementary, we
can deduce that $${\cal M}_{i^*}^{\cal T}
\models {\rm \ `` }\exists {\tilde \mu} \in (\alpha,\kappa) \ ({\tilde \mu}
{\rm \ is \ } <\kappa {\rm -strong}){\rm ."}$$
But ${\cal J}_\kappa^{{\cal M}_{i^*}^{\cal T}} =
{\cal J}_\kappa^{{\cal M}_{i}^{\cal T}}$ and $\alpha < \kappa$ was arbitrary, so that
$$\{ {\tilde \mu} < \kappa \ \colon \ {\tilde \mu}
{\rm \ is \ } <\kappa{\rm -strong \ in \ } {\cal
M}^{\cal T}_i \} {\rm \ is \ unbounded \ in \ } \kappa.$$
Because $F = E^{\cal T}_i$ has critical point
$\kappa$, this shows
that ${\cal M}^{\cal T}_i$ is not below $\zerohandgrenade$.
But this implies that
${\cal M}$ was not below $\zerohandgrenade$ to begin with. Contradiction!

\bigskip
\hfill $\square$ (Claim 1)

\bigskip
{\bf Claim 2.} For all normal ${\cal T}$ on ${\cal M}$ of double
successor length
$\vartheta+2$ do we have the following. Let $i < j+1 < \vartheta+1$
be such that $j+1 \in (i,\vartheta]_T$, and
$i = T$-pred$(j+1) = T$-pred$(\vartheta+1)$. Then ${\rm c.p.}(E^{\cal T}_j) >
{\rm c.p.}(E^{\cal T}_\vartheta)$.

\bigskip
{\sc Proof.} This easily follows from Claim 1.

\bigskip
\hfill $\square$ (Claim 2)

\bigskip
Let us finally prove that every normal ${\cal T}$ on ${\cal M}$ is almost linear.
Suppose not, and let ${\cal T}$
be a normal iteration tree on ${\cal M}$ of length $\vartheta$
such that ${\cal T} {\upharpoonright} \theta$ is almost linear for all $\theta <
\vartheta$, whereas ${\cal T}$ is not almost linear.

\bigskip
{\it Case 1.} $\vartheta$ is a limit ordinal.

\bigskip
Then there is some $i < \vartheta$ such that $i$ has infinitely many immediate
$T$-successors, say $i_k +1$ for all $k<\omega$. We may then apply Claim 2
successively to the almost linear trees ${\cal T} {\upharpoonright} i_k +2$ to get an
infinite descending sequence of ordinals, 
namely $({\rm c.p.}(E^{\cal T}_{i_k}) \colon k \in \omega)$. Contradiction!

\bigskip
{\it Case 2.} $\vartheta$ is a successor ordinal.

\bigskip
By \ref{unique-branch} 
we then know that ${\rm lh}({\cal T})$ must be a double successor,
say ${\rm lh}({\cal T}) =
\vartheta+2$, and
${\cal T} {\upharpoonright} \vartheta+1$ is almost linear.
For $i < \vartheta+1$ let $\kappa_i = {\rm c.p.}(E_i^{\cal T})$, and
$\lambda_i = E_i^{\cal T}(\kappa_i)$. Because ${\cal T}$ is normal, we must have
$i<j<\vartheta+1 \Rightarrow \lambda_i<\lambda_j$.

Set $j = {\cal T}$-${\rm pred}(\vartheta+1)$. Note that $j \notin
[0,\vartheta]_T$ (and so in particular $j < \vartheta$)
by our assumption on ${\cal T}$.
Let then $k < j$ be maximal such that $k \in
[0,j)_T$
as well
as $k \in [0,\vartheta)_T$.

By the normality of ${\cal T}$, $j$ is least such that $\kappa_\vartheta
< \lambda_j$. In particular,
\begin{eqnarray}
\forall i \ ( \ j \leq i < \vartheta+1 \ \Rightarrow
\kappa_\vartheta < \lambda_j \leq \lambda_i \ ).
\end{eqnarray}
Let $i+1$ be minimal in $(k,\vartheta]_T$, so that $E^{\cal T}_i$ is
the first extender used on $[k,\vartheta)_{\cal
T}$.
As ${\cal T} {\upharpoonright} \vartheta+1$ is
almost linear, we must have that $i+1 > j$, i.e., $i \geq j$. Hence by (1),
\begin{eqnarray}
\kappa_\vartheta < \lambda_j \leq \lambda_i \leq \lambda_\vartheta.
\end{eqnarray}

Recall that $k<j$. As $k = {\cal T}$-${\rm pred}(i+1)$ we must have $\kappa_i < \lambda_k$.
We thus have that
\begin{eqnarray}
\kappa_{i} < \kappa_\vartheta {\rm , }
\end{eqnarray}
because
otherwise $k < j$ would be such that $\kappa_\vartheta \leq \kappa_{i} < \lambda_k$;
but $j = {\cal T}$-${\rm pred}(\vartheta+1)$ is least with $\kappa_\vartheta < \lambda_j$.

But (2) and (3) contradict Claim 1!

\bigskip
\hfill $\square$ (\ref{normal-linear})

\bigskip
For future reference, let us state an immediate consequence of the proof of
Claim 1 in the proof of \ref{normal-linear}.

\begin{lemma}\label{not-strong}
Let ${\cal M}$ be a premouse which is below $\zerohandgrenade$, and let ${\cal T}$ be
a normal tree on ${\cal M}$ with last model ${\cal N} = {\cal M}^{\cal T}_\infty$.
Let $i+1 < {\rm lh}({\cal T})$, $F = E_i^{\cal T}$, $\kappa = {\rm c.p.}(F)$, and $\lambda =
F(\kappa)$. Then for no ${\tilde \lambda} \geq \lambda$ and ${\tilde \kappa} \in
[\kappa,\lambda)$ do we have that ${\tilde \kappa}$ is $<{\tilde \lambda}$-strong in
${\cal N}$.
\end{lemma}

\bigskip
Using \ref{strong-ISC+} we could have shown the following. Let ${\cal M}$ be a
$0$-iterable premouse which is such that
for no $E^{\cal M}_\nu \not= \emptyset$
with critical point $\kappa$ do we
have $$\{ \mu < \kappa \ \colon \ \mu {\rm \ is \ } <\kappa{\rm -strong \ in \ } {\cal
M} {\rm \ as \ witnessed \ by \ } {\vec E}^{\cal M}
\} {\rm \ is \ unbounded \ in \ } \kappa.$$ Then any normal iteration tree on ${\cal
M}$ is almost linear. The advantage of \ref{normal-linear} is that it can be applied
to premice of which we don't (yet) know that they are iterable.

\bigskip
The insight which gives \ref{normal-linear} also shows that normal iterations of
phalanxes below $\zerohandgrenade$ are of a simple form.

\begin{defn}\label{phalanx-below-0h}
Let ${\vec {\cal P}} = (({\cal P}_i \colon i<\alpha+1),(\mu_i \colon i<\alpha))$ be a
phalanx. ${\vec {\cal P}}$ is said to be below $\zerohandgrenade$ if every ${\cal
P}_i$ for $i<\alpha+1$ is below $\zerohandgrenade$.
\end{defn}

\begin{lemma}\label{phalanx-iterations}
Let the phalanx ${\vec {\cal P}} =
(({\cal P}_i \colon i<\alpha+1),(\mu_i \colon i<\alpha))$
be below $\zerohandgrenade$, and let ${\cal T}$ be a normal iteration of
${\vec {\cal P}}$. Then there are $i_n < i_{n-1} < ... < i_1 < i_0 = \alpha$
and $\alpha = \beta_0 < \beta_1 < ... < \beta_{n-1} < \beta_n < \beta_{n+1}
= {\rm lh}({\cal T})$
suchthat for all $\gamma < {\rm lh}({\cal T})$ and for all $k \leq n$ do we have that
$$\gamma \in [\beta_k,\beta_{k+1}) \ \Leftrightarrow i_k T \gamma.$$
\end{lemma}

The lemma says, among other things,
that we can write ${\cal T}$ as
$${\cal T} = {\cal T}_0 {}^\frown {\cal T}_1 {}^\frown ... {}^\frown {\cal T}_n$$
where
${\cal T}_0$ is an iteration of ${\cal P}_{i_0} = {\cal P}_\alpha$ (and may be
trivial) and
${\cal T}_{k+1}$ is an iteration of ${\cal P}_{i_{k+1}}$ (except for the fact
that its first extender is taken from the last model of ${\cal T}_{k}$).
However, the proof of \ref{phalanx-iterations} is straightforward in the light of the
proof of \ref{normal-linear} and may be left to the reader.

Lemma \ref{normal-linear} is
optimal in the following sense. 

\begin{lemma}\label{new}
Let ${\cal M}$ be a premouse with top
extender $F$ such that for $\kappa = {\rm c.p.}(F)$,
there are arbitrary large $\mu < \kappa$ such that $\mu$ is $< \kappa$ strong in
${\cal M}$ as witnessed by ${\vec E}^{\cal M}$.
(In particular, ${\cal M}$ is not below $\zerohandgrenade$.)
Then assuming that ${\cal M}$ is sufficiently iterable there is a
normal iteration tree ${\cal T}$ on ${\cal M}$ of length $4$ which is not almost
linear.
\end{lemma}

{\sc Proof.} 
Let $\lambda$ be the largest cardinal of ${\cal M}$. Let $\nu_0 <
\lambda$ be
such that $E^{\cal M}_{\nu_0} \not= \emptyset$ is total on ${\cal
M}$ and has critical point
$\kappa_0 > \kappa$. Let $\pi_{0 1}^{\cal T} \ \colon {\cal M}
\rightarrow_{E^{\cal M}_{\nu_0}} {\cal M}_1$. Let $F'$ be the top
extender of ${\cal M}_1$, and let
$\pi_{0 2}^{\cal T} \ \colon {\cal M}
\rightarrow_{F'} {\cal M}_2$.
Let $F''$ be the top extender of ${\cal M}_2$.

Note that $F'(\kappa) = \pi^{\cal T}_{0 1}(\lambda)$
is the critical point of $F''$, and that
$F'(\kappa) \geq \lambda > \nu_0$, so that we may pick some $\mu \in
(\nu_0,F'(\kappa))$ which is
$< F'(\kappa)$-strong in ${\cal M}_1$ as witnessed by ${\vec E}^{{\cal M}_1}$.
Using $F''$, $\mu$ is also
$< F''(F'(\kappa))$-strong in ${\cal M}_2$ as witnessed by ${\vec E}^{{\cal
M}_2}$.
We may thus pick $\nu_2 > {\cal M}_1 \cap {\rm OR}$
such that $E^{{\cal M}_2}_{\nu_2} \not= \emptyset$ is total on ${\cal
M}_2$ and has critical point $\mu$. Let $\pi^{\cal T}_{1 3} \ \colon
{\cal M}_1 \rightarrow_{E^{{\cal M}_2}_{\nu_2}} {\cal M}_3$.

It is easy to see that we have constructed a normal iteration tree on
${\cal M}$ of length $4$ which is not almost linear.

\bigskip
\hfill $\square$ (\ref{new})

\[
\xymatrix{
  &   &   & 3 \\
 2 &  & 1 \ar[ru] &  \\
  & 0 \ar[lu] \ar[ru] &   &  \\
  &  {}\save[]*{\txt{{\sc Figure 2.} A normal tree 
generated by $\zerohandgrenade$}} \restore & &  \\
}
\]

\bigskip

One can easily generalize \ref{normal-linear} and 
\ref{new}. 

\begin{defn}
Let ${\cal M}$ be a premouse, and let
$\mu < \kappa \leq {\cal M} \cap {\rm OR}$. We call $\mu$ $<\kappa$-$0$-strong 
in ${\cal
M}$ if $\mu$ is a measurable cardinal in 
${\cal J}^{\cal M}_\kappa$ as witnessed by ${\vec E}^{\cal M}$, i.e,
there is an extender $E^{\cal M}_\nu \not= \emptyset$ with critical point 
$\mu$ and such that $\nu > \mu^{+{\cal J}^{\cal M}_\kappa}$.
For $n<\omega$ we call $\mu$
$<\kappa$-$(n+1)$-strong in ${\cal M}$
if $\mu$ is $<\kappa$-strong in ${\cal
M}$ as witnessed by ${\vec E}^{\cal M}$ and
there are arbitrary large ${\bar \mu} < \mu$ such that ${\bar \mu}$ is
$<\mu$-$n$-strong in ${\cal M}$.
\end{defn}

\begin{defn}
Let ${\cal M}$ be a premouse, and let $n < \omega$. ${\cal M}$ is said to be 
below
$0^{n \: \mid^{\! \! \! \bullet}}$ if there is no $\kappa$
which is the critical point of an extender
$E^{\cal M}_\nu \not= \emptyset$
with the property that 
$$\{ \mu < \kappa \ \colon \ \mu {\rm \ is \ } 
<\kappa{\rm -}n{\rm -strong \ in \ } {\cal M} \} {\rm \ is 
\ unbounded \ in \ } 
\kappa.$$
\end{defn}

Notice that 
${\cal M}$ is below $0^{1 \: \mid^{\! \! \! \bullet}}$ if 
and only if ${\cal M}$ is below $\zerohandgrenade$.
The following lemma, whose easy proof we omit, generalizes \ref{new}.

\begin{lemma}
Let ${\cal M}$ be a premouse which is not below
$0^{n \: \mid^{\! \! \! \bullet}}$. 
Then -- provided all the
involved ultrapowers are
well-founded -- one can build a normal alternating chain on ${\cal M}$ of length $3+n$, i.e.,
a normal iteration tree with tree structure $T$ given by $i T j$ if and only if $i=0$ or $(i \leq
j \wedge i \equiv j \ (mod \ 2))$ for $i \leq j < 3+n$.
\end{lemma}

This is tight in the sense
that there are no iteration trees with so 
much ``jumping from branch to branch'' if
every premouse is below $0^{n \: \mid^{\! \! \! \bullet}}$. We 
leave it to the reader to formulate a precise 
generalization of \ref{normal-linear}.

Now suppose that ${\cal M}$ is a premouse with top extender $F$, which is
such that, setting $\kappa =
{\rm c.p.}(F)$, $${\cal M}  \models
\{ \mu < \kappa \ \colon \ \mu {\rm \ is \ } <\kappa{\rm -strong \
in \ } {\cal M} {\rm \ as \ witnessed \ by \ } {\vec E}^{\cal M}
\} {\rm \ is \ stationary.}$$ Then
for no $n<\omega$ is ${\cal M}$ below $0^{n \: \mid^{\! \! \! \bullet}}$, and
-- again provided all the
involved ultrapowers are
well-founded -- for every $n<\omega$ can one
build a normal alternating chain on ${\cal M}$ of length $3+n$. Recall that (by \cite[Theorem 6.1]{FSIT};
see also \cite[\S 6]{NFS}) to build an alternating chain of
length $\omega$ requires an assumption at the level of a (definably) Woodin cardinal.

\bigskip
We also want to mention an extension of \ref{normal-linear} into another direction,
namely by
revising the definition of ``normal trees.''

\begin{defn}
Let ${\cal M}$ be a premouse. ${\cal M}$ is said to be below $0\ulcorner$
if for no $\kappa$ which is the critical point of some extender
$F = E^{\cal M}_\nu \not= \emptyset$ 
do we have that, setting $\lambda = F(\kappa)$,
there is some $\mu \in (\kappa,\lambda)$ with $${\cal J}^{\cal M}_\lambda \models
{\rm
`` } \mu {\rm \ is \ a \ strong \ cardinal," }$$ and $F$
has a generator (strictly) above
$\mu$.
\end{defn}

If we were to let $T$-${\rm pred}(i+1)$ be the least $j \leq i$ such that
${\rm c.p.}(E_i^{\cal T}) < $ the ${\rm sup}$ of the generators of
$E_j^{\cal T}$ in a D-normal iteration tree ${\cal T}$
then we would get that any D-normal iteration tree of a premouse which is below
$0\ulcorner$ is almost linear.\footnote{$0\ulcorner$ is mentioned 
in the introduction to
\cite{dodd}. Cf. also \cite[\S 3]{ronald}.}
This observation can
be
used to develop the theory of $K$ up to the level of 
$0\ulcorner$ (cf. \cite{d-premice}).

\section{Iterability, and the existence of $K^c$.}

This section introduces $K^c$, a preliminary version of $K$. $K^c$ is constructed
recursively, exactly as in \cite[\S 11]{NFS} 
(see also \cite[p. 6 f.]{CMIP}),
except for the fact that we only require
new extenders to be countably complete (rather than ``certifiable'' as in
\cite[\S 11]{NFS} 
when they are put onto the sequence. Recall that an extender $F$ with $\kappa =
{\rm c.p.}(F)$ is called countably complete if for all $(a_n , X_n \colon n<\omega)$ such
that $a_n \in [F(\kappa)]^{<\omega}$ and $a_n \in F(X_n)$ for all $n<\omega$ there is
some
order preserving
$\tau \colon \bigcup_{n<\omega} \ a_n \rightarrow \kappa$ with $\tau {\rm " } a_n
\in X_n$ for every $n<\omega$.

\bigskip
Let us define premice ${\cal N}_\xi$ and ${\cal M}_\xi$ by induction on $\xi \in {\rm OR}
\cup \{ {\rm OR} \}$
(cf. \cite[\S 10 p. 9 f.]{NFS}).
We let ${\cal N}_0 = (J_\omega;\in,\emptyset)$.
Having defined ${\cal N}_\xi$, we
let the construction break down unless ${\cal N}_\xi$ is reliable, i.e., unless
for all $n \leq \omega$, ${\frak C}_n({\cal N}_\xi)$ is $n$-iterable.
If ${\cal N}_\xi$ is ``reliable'' then we continue by
setting $${\cal
M}_\xi = {\frak C}_\omega({\cal N}_\xi).$$
Now suppose that ${\cal M}_\xi$ has been defined.

\bigskip
{\it Case 1.} ${\cal M}_\xi = (J_\alpha[{\vec E}];\in,{\vec E})$
is passive, and there is a
unique countably complete extender $F$
such that $$(J_\alpha[{\vec E}];\in,{\vec
E},F)$$
is a premouse below $\zerohandgrenade$.
In this case we set ${\cal N}_{\xi+1} = (J_\alpha[{\vec E}];\in,{\vec
E},F)$.

\bigskip
{\it Case 2.} Otherwise. Then we just construct one more step. I.e., if
${\cal M}_\xi = (J_\alpha[{\vec E}];\in,{\vec E})$ then we let ${\cal N}_{\xi+1} =
(J_{\alpha+1}[{\vec E}];\in,{\vec E})$, and if
${\cal M}_\xi = (J_\alpha[{\vec E}];\in,{\vec E},F)$ (with $F \not=
\emptyset$) then we let ${\cal N}_{\xi+1} =
(J_{\alpha+1}[{\vec E}^\frown \langle F \rangle];\in,{\vec
E}^\frown F)$.

\bigskip
Now suppose that ${\cal M}_\xi$ has been defined
for all $\xi < \lambda$ where $\lambda$ is a limit ordinal. We may let
$$\mu = {\rm sup}_{\xi < \lambda} \ {\rm min}_{\xi < \zeta < \lambda} \ \rho_\omega({\cal
M}_\zeta)^{+{\cal M}_\xi} {\rm , }$$ and let ${\cal N}_\lambda$ be that passive
premouse of height $\mu$ such that for all ${\bar \mu} < \mu$,
${\cal J}^{{\cal N}_\lambda}_{\bar \mu}$ is the eventual value of
${\cal J}^{{\cal N}_\xi}_{\bar \mu}$ as $\xi \rightarrow \lambda$ 
(cf. \cite[\S
10 p. 10]{NFS},
where it is shown that ${\cal N}_\lambda$ is well-defined if all ${\cal
M}_\xi$, $\xi < \lambda$, are well-defined).

\bigskip
Notice that, whereas there is no restriction on
${\rm cf}^V({\rm c.p.}(F))$ in Case 1 above, it is automatic that ${\rm cf}^V({\rm c.p.}(F)) > \omega$,
because no countably complete extender can have a critical point with cofinality
$\omega$.

We are now going to prove that the construction never breaks down.
For this
purpose, we need the resurrection maps provided by the following lemma. (We handle
the issue of resurrecting slightly different than both
\cite[\S 12]{FSIT} and \cite[\S 9]{CMIP}.)

\begin{lemma}\label{resurrecting} {\bf and Definition}
Suppose that ${\cal N}_{\xi_0}$ exists for some $\xi_0 \in {\rm OR}$.
Then there are $$(\tau_{(\xi,\eta)} \
\colon \ \xi \leq \xi_0 \wedge \eta \leq {\cal
N}_\xi \cap {\rm OR})$$ and 
$$(\varphi(\xi,\eta) \
\colon \ \xi \leq \xi_0 \wedge \eta \leq {\cal
N}_\xi \cap {\rm OR})$$
such that for every pair $(\xi,\eta)$ with
$\xi \leq \xi_0$ and $\eta \leq {\cal
N}_\xi \cap {\rm OR}$ do we have that $$\tau_{(\xi,\eta)} \colon {\cal J}_\eta^{{\cal
N}_\xi} \rightarrow {\cal N}_{(\varphi_{(\xi,\eta)})} {\rm , }$$
where, for every $n < \omega$,
$$\tau_{(\xi,\eta)} {\upharpoonright} ({\cal J}_\eta^{{\cal N}_\xi})^n
\ \colon \ ({\cal J}_\eta^{{\cal N}_\xi})^n \rightarrow_{\Sigma_0} ({\cal
N}_{(\varphi(\xi,\eta))})^n {\rm ; }$$ the map $\tau_{(\xi,\eta)}$ is called a resurrection map.
Moreover, whenever $\eta < \eta' \leq {\cal
N}_\xi \cap {\rm OR}$ then, setting $\rho = {\rm min} \{ \rho_\omega({\cal J}_{\bar
\eta}^{{\cal N}_\xi}) \colon \eta \leq {\bar \eta} < \eta' \}$,
we have that $\tau_{(\xi,\eta)}$ agrees with $\tau_{(\xi,\eta')}$ up to $\rho$, i.e.,
$${\cal J}^{{\cal N}_{(\varphi(\xi,\eta))}}_{\tau_{(\xi,\eta')}(\rho)} =
{\cal J}^{{\cal N}_{(\varphi(\xi,\eta'))}}_{\tau_{(\xi,\eta')}(\rho)} {\rm , \ and }$$
$$\tau_{(\xi,\eta)} {\upharpoonright} {\cal J}^{{\cal N}_\xi}_\rho =
\tau_{(\xi,\eta')} {\upharpoonright} {\cal J}^{{\cal N}_\xi}_\rho.$$
\end{lemma}

{\sc Proof sketch.} The maps $\tau_{(\xi,\eta)}$ and their target models are chosen
by induction on $\xi$.
Suppose $\tau_{(\xi,\eta)}$ has been defined for all $\eta \leq {\cal
N}_\xi \cap {\rm OR}$. We let $\tau_{(\xi+1,{\cal
N}_{\xi+1} \cap {\rm OR})} = {\rm id}$. If $\eta < {\cal
N}_{\xi+1} \cap {\rm OR}$ then we have that ${\cal J}_\eta^{{\cal N}_{\xi+1}} =
{\cal J}_\eta^{{\cal M}_{\xi}}$. Let $\pi \colon {\cal M}_{\xi} \rightarrow {\cal
N}_\xi$ be the core map. We then let
$$\varphi(\xi+1,\eta) = \varphi(\xi,\pi(\eta)) {\rm , \ and }$$
$$\tau_{(\xi+1,\eta)} = \tau_{(\xi,\pi(\eta))}
\circ \pi {\upharpoonright} {\cal J}_\eta^{{\cal M}_\xi}.$$
It is straightforward to verify the required agreement between appropriate
$\tau_{(\xi+1,\eta)}$ and $\tau_{(\xi+1,\eta')}$.

At limit stages $\lambda$ we exploit the fact that any proper initial segment of
${\cal N}_\lambda$ is an initial segment of some ${\cal N}_\xi$ for $\xi < \lambda$.
We leave the details to the reader.

\bigskip
\hfill $\square$ (\ref{resurrecting})

\bigskip
The following two lemmas, \ref{normal-iterability} and \ref{iterability}, proving
normal and full $n$-iterability of ${\frak C}_n({\cal
N}_\xi)$, are shown
simultaneously by induction on $n < \omega$.

\begin{lemma}\label{normal-iterability}
Suppose that ${\cal N}_\xi$ exists for some $\xi \in {\rm OR}$ and $n <
\omega$. Then ${\frak C}_n({\cal
N}_\xi)$ is normally $n$-iterable.
\end{lemma}

{\sc Proof.} Let us fix some $n < \omega$
and assume that ${\frak C}_m({\cal
N}_\xi)$ is $m$-iterable
for all $m < n$. By \cite[\S 4, and Lemmas 6.1.5 and 8.1]{FSIT},
this will buy us that
we already know that for all $k \leq n$,
any $k$-bounded iteration of ${\frak C}_n({\cal
N}_\xi)$ will move
$p_{{\frak C}_n({\cal
N}_\xi),k}$ as well as all
standard parameters of all non-simple iterates of
${\frak C}_n({\cal
N}_\xi)$ correctly.

By \ref{unique-branch} and \ref{normal-linear},
if the lemma fails then there is a
putative normal $n$-bounded iteration tree ${\cal U}$ of minimal
length $\beta$
on ${\frak C}_n({\cal N}_\xi)$ such that either ${\cal U}$ has a last
ill-founded model, or else ${\cal U}$ has limit length and
${\cal D}^{\cal U} \cap b'$ is infinite
where $b'$ is the unique cofinal branch through $U$.
A standard argument yields a (fully elementary) embedding
$$\sigma \ \colon {\cal M} \rightarrow {\frak C}_n({\cal N}_\xi)$$ such that ${\cal
M}$ is countable and there is a putative normal $n$-bounded iteration tree ${\cal
T}$ of minimal length $\alpha$ on ${\cal M}$ such that $\alpha <
\omega_1$ and either ${\cal T}$ has a last ill-founded model,
or else ${\cal T}$ has limit length and ${\cal D}^{\cal T} \cap b$ is infinite
where $b$ is the unique cofinal branch through $T$. Let
$${\cal T} = (({\cal M}_i^{\cal T}, \pi^{\cal T}_{i j} \colon
i T j < \alpha+1),(E_i^{\cal T} \colon i < \alpha),T).$$
For $i<\alpha$ we set
$\kappa_i = {\rm c.p.}(E_i^{\cal T})$, $\lambda_i = E_i^{\cal T}(\kappa_i)$, and
$\eta_i = {\rm dom}(\pi^{\cal T}_{i^\star i+1}) \cap {\rm OR}$ where $i^\star = T$-pred$(i+1)$.
We set $n(0) = n$, and
for $i>0$ we let $n(i)$ be such that
$$\exists j \ ( \ j <_T i \wedge \forall k+1 \in (j,i]_T \ {\cal M}_{k+1}^{\cal T} =
{\rm Ult}_{n(i)}({\cal J}^{{\cal M}^{\cal T}_{T-pred(k+1)}}_{\eta_k};E^{\cal T}_k) \ ).$$
(That is, $n(i) = {\rm deg}^{\cal T}(i)$ for successor ordinals $i$.)
Notice that for all $i<\alpha+1$ do we have that
${\cal M}^{\cal T}_i$ is $n(i)$-sound.

We shall now pick for all $i
\leq \alpha$ some $\xi(i) \leq \xi$ together with
a weak $n(i)$-embedding $${\tilde \sigma}_i \
\colon \ {\cal M}^{\cal T}_i \rightarrow {\frak C}_{n(i)}({\cal N}_{\xi(i)})
{\rm . }$$
The maps ${\tilde \sigma}_i$ may be obtained by composing some $\sigma_i$ with
(the inverse of) a core
map. We shall therefore just 
recursively pick $\xi(i)$ and $\sigma_i$.
We'll inductively maintain that the following three requirements are met.

\bigskip
{\bf R 1}${}_i$ $\ \ \forall \ k \in [0,i]_T$
do we have that
$\sigma_k \colon {\cal M}^{\cal T}_k
\rightarrow {\cal
N}_{\xi(k)}$ is the extension of
$$\sigma_k {\upharpoonright} ({\cal M}^{\cal T}_k)^{n(k)} \ \colon \
({\cal M}^{\cal T}_k)^{n(k)}\rightarrow_{\Sigma_0}
({\cal
N}_{\xi(k)})^{n(k)}$$
given by the downward extension of embeddings lemma.

\bigskip
{\bf R 2}${}_i$ $\ \ \forall \ k+1 \in (0,i]_T$, setting
$k^* = {\cal T}$-${\rm pred}(k+1)$, we have that
$\sigma_i$ agrees with
$\tau_{(\xi(k^*),\sigma_{k^*}(\eta_k))} \circ \sigma_{k^*}$
up to $\kappa_k$, i.e., $${\cal J}^{{\cal N}_{\xi(i)}}_{\sigma_i(\kappa_k)} =
{\cal J}^{{\cal N}_{\varphi(\xi(k^*),\sigma_{k^*}(\eta_k))}}_{\sigma_i(\kappa_k)}
{\rm , \ and }$$
$$\sigma_i {\upharpoonright} {\cal J}_{\kappa_k}^{{\cal M}_i^{\cal T}} =
\tau_{(\xi(k^*),\sigma_{k^*}(\eta_k))} \circ
\sigma_{k^*} {\upharpoonright} {\cal J}_{\kappa_k}^{{\cal M}_{k^*}^{\cal T}}.$$

\bigskip
{\bf R 3}${}_i$ $\ \ $ if $k T j \in (0,i]_T$ and ${\cal D}^{\cal T} \cap [k,j]_T =
\emptyset$ then $\sigma_j \circ \pi^{\cal T}_{k j } = \sigma_k$.

\bigskip
To commence, we set $\xi(0) = \xi$ and $\sigma_0 = \pi \circ \sigma$ where $\pi \colon
{\frak C}_n({\cal N}_\xi) \rightarrow {\cal N}_\xi$ is the core map. It is
easy to see that {\bf R
1}${}_0$ holds. Moreover, {\bf R 2}${}_0$ and {\bf R 3}${}_0$ are vacuously true.

Now suppose we have defined $\sigma_j \colon {\cal M}^{\cal T}_j
\rightarrow {\cal N}_{\xi(j)}$ for all $j \leq i$
in such a way that
{\bf A 1}${}_i$, {\bf A 2}${}_i$, and {\bf A 3}${}_i$ hold. Suppose that $i+1 <
\alpha$, so that the $n(i+1)$-ultrapower $$\pi^{\cal T}_{i^* i+1} \ \colon \
{\cal J}^{{\cal M}^{\cal T}_{i^*}}_{\eta_i} \rightarrow_F {\cal
M}^{\cal T}_{i+1}$$ exists with $F = E^{\cal T}_i$
and $i^* = {\cal T}$-${\rm pred}(i+1)$. It is convenient to split
the construction into two cases.

\bigskip
{\it Case 1.} $i^* < i$.

\bigskip
In this case, $\lambda_{i^*}$ is a cardinal in ${\cal M}^{\cal
T}_i$, so that by $\kappa_i < \lambda_{i^*}$
we have that
$F$ is a total extender on ${\cal M}^{\cal T}_i$.
Hence $G = \sigma_i(F)$ is countably complete, being
total on ${\cal N}_{\xi(i)}$. So we may pick
$\rho \colon \sigma_i(\lambda_i) \cap
{\rm ran}(\sigma_i) \rightarrow \sigma_i(\kappa_i)$ order preserving
such that $$a \in G(X) \ \Rightarrow \
\rho {\rm " } a \in X$$
for all appropriate $a$, $X \in {\rm ran}(\sigma_i)$.
Let $\xi(i+1) = \varphi(\xi(i^*),\sigma_{i^*}(\eta_i))$,
i.e.,  ${\cal M}_{\xi(i+1)}$
is the target model of the
resurrection map
$$\tau_{(\xi(i^*),\sigma_{i^*}(\eta_i))} \ \colon {\cal J}^{{\cal
M}_{\xi(i^*)}}_{\sigma_{i^*}(\eta_i)} \rightarrow {\cal
N}_{\xi(i+1)} {\rm . }$$

By \ref{normal-linear}, $i^* \in [0,i)_{\cal
T}$.
Let $k+1$ be minimal in $(i^*,i]_{\cal T}$. In particular, ${\cal
T}$-${\rm pred}(k+1) = i^*$.
Now {\bf R 2}${}_i$ gives us that
$\sigma_i$ agrees with
$\tau_{\xi(i^*),\sigma_{i^*}(\eta_k)} \circ \sigma_{i^*}$
up to $\kappa_k$. But by the proof of \ref{normal-linear}
we have
$\kappa_i < \kappa_k$, so that
$\eta_i \geq \eta_k$ and $$\forall \ \xi
\in [\eta_k , \eta_i) \ (\rho_\omega({\cal J}^{{\cal M}_{i^*}^{\cal
T}}_\xi) > \kappa_i) {\rm , \ and \ thus }$$
$$\forall \ \xi
\in [\sigma_{i^*}(\eta_k) , \sigma_{i^*}(\eta_i))
\ (\rho_\omega({\cal J}^{{\cal N}_{\xi(i^*)}}_\xi) > \sigma_{i^*}(\kappa_i))
{\rm , }$$
which by \ref{resurrecting}
implies that $\tau_{(\xi(i^*),\sigma_{i^*}(\eta_i))}$ agrees with
$\tau_{(\xi(i^*),\sigma_{i^*}(\eta_k))}$ up to
$\sigma_{i^*}(\kappa_i^+)$ (where $\kappa_i^+$ is calculated in
${\cal J}^{{\cal M}^{\cal T}_{i^*}}_{\eta_i}$).
These agreements combined easily give that $\sigma_i$ agrees with
$\tau_{(\xi(i^*),\sigma_{i^*}(\eta_i))} \circ \sigma_{i^*}$ up to
$\kappa_i^+$. Let us write $\tau = \tau_{(\xi(i^*),\sigma_{i^*}(\eta_i))}$.

Let $$p = p_{{\cal J}^{{\cal M}^{\cal T}_{i^*}}_{\eta_i},n(i+1)}.$$
By {\bf R 1}${}_i$ we will have that
$$\sigma_{i^*}(p) = p_{{{\cal J}^{{\cal
N}_{\xi(i^*)}}_{\sigma_{i^*}(\eta_i)}},n(i+1)}.$$ Moreover, by
\ref{resurrecting}
we then get that
$$\tau \circ \sigma_{i^*}(p) = p_{{{\cal N}_{\xi(i+1)}},n(i+1)} {\rm , \ and }$$
$$\tau {\upharpoonright} ({\cal J}^{{\cal
N}_{\xi(i^*)}}_{\sigma_{i^*}(\eta_i)})^{n(i+1),\sigma_{i^*}(p)}
\colon ({\cal J}^{{\cal
N}_{\xi(i^*)}}_{\sigma_{i^*}(\eta_i)})^{n(i+1),\sigma_{i^*}(p)}
\rightarrow_{\Sigma_0}
({\cal N}_{\xi(i+1)})^{n(i+1),\tau \circ \sigma_{i^*}(p)}.$$
%

We may now define $${\bar \sigma} \ \colon
\ ({\cal M}^{\cal T}_{i+1})^{n(i+1),\pi^{\cal T}_{i^* i+1}(p)} \rightarrow_{\Sigma_0}
({\cal N}_{\xi(i+1)})^{n(i+1),\tau \circ \sigma_{i^*}(p)} =
({\cal N}_{\xi(i+1)})^{n(i+1)}$$ by
setting $${\bar \sigma}([a,f]) = \tau \circ
\sigma_{i^*}(f)(\rho {\rm " }
\sigma_i(a)).$$

To show that this is well-defined and $\Sigma_0$-elementary we may
reason as follows. Let $\Phi$ be a $\Sigma_0$ formula. Then
$$({\cal M}^{\cal T}_{i+1})^{n(i+1),\pi^{\cal T}_{i^* i+1}(p)}
\ \models \ \Phi([a_1,f_1],...,[a_k,f_k]) \ \ \Leftrightarrow$$
$$(a_1,...,a_k) \in F( \{ (u_1,...,u_k) \ \colon \
({\cal J}^{{\cal M}^{\cal T}_{i^*}}_{\eta_i})^{n(i+1),p} \ \models \
\Phi(f_1(u_1),...,f_k(u_k)) \} ) \ \ \Leftrightarrow$$
$$(\sigma_i(a_1),...,\sigma_i(a_k)) \in G( \sigma_i(
\{ (u_1,...,u_k) \ \colon \
({\cal J}^{{\cal M}^{\cal T}_{i^*}}_{\eta_i})^{n(i+1),p} \ \models \
\Phi(f_1(u_1),...,f_k(u_k)) \} ) ) {\rm , }$$
which, by the amount of agreement of $\sigma_i$ with $\tau \circ
\sigma_{i^*}$, holds if and only if

\bigskip
$(\sigma_i(a_1),...,\sigma_i(a_k)) \in$ $$G( \tau \circ \sigma_{i^*}(
\{ (u_1,...,u_k) \ \colon \
({\cal J}^{{\cal M}^{\cal T}_{i^*}}_{\eta_i})^{n(i+1),p} \ \models \
\Phi(f_1(u_1),...,f_k(u_k)) \} ) ) \ \ \Leftrightarrow$$

$(\sigma_i(a_1),...,\sigma_i(a_k)) \in$ $$G( \{ (u_1,...,u_k) \ \colon \
({\cal N}_{\xi(i+1)})^{n(i+1),\tau \circ \sigma_{i^*}(p)} \ \models \
\Phi(\tau \circ \sigma_{i^*}(f_1)(u_1),...,\tau \circ 
\sigma_{i^*}(f_k)(u_k)) \} ) \ \ \Leftrightarrow$$
$$({\cal N}_{\xi(i+1)})^{n(i+1),\tau \circ \sigma_{i^*}(p)}
\ \models \ \Phi(\tau \circ
\sigma_{i^*}(f_1)(\rho {\rm " } \sigma_i(a_1),...,\tau \circ
\sigma_{i^*}(f_k)(\rho {\rm " } \sigma_i(a_k))).$$
We let $\sigma_{i+1}$ be the extension of ${\bar \sigma}$ given by the
downward extension of embeddings lemma.
By the remark in the first paragraph of this proof, we'll have that
$$\pi^{\cal T}_{i^* i+1}(p) = p_{{\cal M}_{i+1}^{\cal T},n(i+1)}.$$ Also,
by the definition of $\sigma_{i+1}$,
$$\sigma_{i+1} ( \pi^{\cal T}_{i^* i+1}(p)) =
\tau \circ \sigma_{i^*}(p) = p_{{\cal N}_{\xi(i+1)},n(i+1)} {\rm , }$$
and hence $$\sigma_{i+1}(p_{{\cal M}_{i+1}^{\cal T},n(i+1)}) =
p_{{\cal N}_{\xi(i+1)},n(i+1)}.$$
Hence we
have established {\bf R 1}${}_{i+1}$.

Let us verify {\bf R 2}${}_{i+1}$.
It is clear by construction that $\sigma_{i+1}$ agrees
with $\tau \circ \sigma_{i^*}$ up to $\kappa_i$. So let $k+1 \in
(0,i+1]_T$ be such that $k < i$. Then $\lambda_k \leq \kappa_i$
is a cardinal in ${\cal M}^{\cal T}_{i^*}$, so that
$$\forall \ \xi \in [\lambda_k,{\rm OR} \cap {\cal M}^{\cal T}_{i^*}) \ \
\rho_\omega({\cal J}^{{\cal M}^{\cal T}_{i^*}}_\xi)
\geq \lambda_k {\rm , thus
}$$
$$\forall \ \xi \in [\sigma_{i^*}(\lambda_k),{\rm OR} \cap {\cal N}_{\xi(i^*)}) \ \
\rho_\omega({\cal J}^{{\cal N}_{\xi(i^*)}}_\xi)
\geq \sigma_{i^*}(\lambda_k) {\rm ,
}$$
which implies that $$\tau_{(\xi(i^*),\sigma_{i^*}(\eta_i))}
{\upharpoonright} \sigma_{i^*}(\lambda_k) = {\rm id}.$$
Combining this with {\bf R 2}${}_{i^*}$ we get that
$\tau_{(\xi(k^*),\sigma_{k^*}(\eta_k))}
\circ \sigma_{k^*}$ agrees with
$\tau_{(\xi({i^*}),\sigma_{i^*}(\eta_i))} \circ \sigma_{i^*}$
up to $\kappa_k$, which in turn agrees
with $\sigma_{i+1}$ up to $\kappa_k$ by
the construction of $\sigma_{i+1}$.
This proves {\bf R 2}${}_{i+1}$.

It is now straightforward to verify {\bf R 3}${}_{i+1}$.

\bigskip
{\it Case 2.} $i^* = i$.

\bigskip
In this case, $F$ may be partial on ${\cal M}^{\cal T}_i$.
Let $\xi(i+1) = \varphi(\xi(i),\sigma_i(\eta_i))$, i.e.,
${\cal M}_{\xi(i+1)}$ is the target model of the resurrection map
$$\tau_{(\xi(i),\sigma_i(\eta_i))} \ \colon
{\cal J}^{{\cal M}_{\xi(i)}}_{{\sigma_i}(\eta_i)} \rightarrow
{\cal N}_{\xi(i+1)} {\rm , }$$ and let us write $\tau =
\tau_{(\xi(i),\sigma_i(\eta_i))}$.
We then have that $G = \tau \circ \sigma_i(F)$ is
countably complete, being total on ${\cal
N}_{\xi(i+1)}$. So we may pick an orderpreserving $\rho \colon \tau \circ
\sigma_i(\kappa_i) \cap {\rm ran}(\tau \circ \sigma_i) \rightarrow
\tau \circ \sigma_i(\kappa_i)$ such that
$$a \in G(X) \ \Rightarrow \ \rho {\rm " } a \in X$$ for all appropriate
$a$, $X \in {\rm ran}(\tau \circ \sigma_i)$.

Let $$p = p_{{{\cal J}^{{\cal M}^{\cal T}_{i}}_{\eta_i}},n(i+1)}.$$
%

We may define $${\bar \sigma} \ \colon
\ ({\cal M}^{\cal T}_{i+1})^{n(i+1),\pi^{\cal T}_{i i+1}(p)} \rightarrow_{\Sigma_0}
({\cal N}_{\xi(i+1)})^{n(i+1),\sigma_{i}(p)}$$ by
setting $${\bar \sigma}([a,f]) = \tau \circ
\sigma_{i}(f)(\rho {\rm " }
\tau \circ \sigma_i(a)).$$

To show that this is well-defined and $\Sigma_0$-elementary we may
reason as follows. Let $\Phi$ be a $\Sigma_0$ formula. Then
$$({\cal M}^{\cal T}_{i+1})^{n(i+1),\pi^{\cal T}_{i i+1}(p)}
\ \models \ \Phi([a_1,f_1],...,[a_k,f_k]) \ \ \Leftrightarrow$$
$$(a_1,...,a_k) \in F( \{ (u_1,...,u_k) \ \colon \
({\cal J}^{{\cal M}^{\cal T}_{i}}_{\eta_i})^{n(i+1),p} \ \models \
\Phi(f_1(u_1),...f_k(u_k)) \} ) \ \ \Leftrightarrow$$

\bigskip
$(\tau \circ \sigma_i(a_1),...,\tau \circ \sigma_i(a_k)) 
\in$ $$G( \tau \circ \sigma_i(
\{ (u_1,...,u_k) \ \colon \
({\cal J}^{{\cal M}^{\cal T}_{i}}_{\eta_i})^{n(i+1),p} \ \models \
\Phi(f_1(u_1),...,f_k(u_k)) \} ) ) \ \ \Leftrightarrow$$

$(\tau \circ \sigma_i(a_1),...,\tau \circ \sigma_i(a_k)) \in$ 
$$G( \{ (u_1,...,u_k) \ \colon \
({\cal N}_{\xi(i+1)})^{n(i+1),\tau \circ \sigma_{i}(p)} \ \models \
\Phi(\tau \circ 
\sigma_{i}(f_1)(u_1),...,\tau \circ 
\sigma_{i}(f_k)(u_k)) \} ) \ \ \Leftrightarrow$$
$$({\cal N}_{\xi(i+1)})^{n(i+1),\tau \circ \sigma_{i}(p)}
\ \models \ \Phi(\tau \circ
\sigma_{i}(f_1)(\rho {\rm " } \tau \circ \sigma_i(a_1)),...,\tau \circ
\sigma_{i}(f_k)(\rho {\rm " } \tau \circ \sigma_i(a_k))).$$
We may now let $\sigma_{i+1}$ be the extension of ${\bar \sigma}$ given by the
downward extension of embeddings lemma.
We can then argue exactly as in Case 1 to see
that we have established {\bf R 1}${}_{i+1}$.

As for {\bf R 2}${}_{i+1}$ and {\bf R 3}${}_{i+1}$, the proofs are
similar as in the previous case.

Now suppose that we have defined $\sigma_i \colon {\cal M}^{\cal T}_i \rightarrow
{\cal N}_{\xi(i)}$ for all $i < \lambda$ where $\lambda < \alpha$ is a limit ordinal.
By our minimality assumption on $\alpha$ we have that ${\cal D}^{\cal T} \cap
(0,\lambda]_T$ is finite. Hence we may use 
that for all $i<\lambda$ the statement {\bf R 3}${}_{i}$ holds
to define $${\sigma_\lambda}(x) = \sigma_i \circ (\pi^{\cal T}_{i \lambda})^{-1}(x)$$
where $i <_T \lambda$ is large enough. It is straightforward to verify that with this
definition we have {\bf R 1}${}_{\lambda}$, {\bf R 2}${}_{\lambda}$, and
{\bf R 3}${}_{\lambda}$.

Now if ${\cal T}$ has a last model ${\cal M}^{\cal T}_{\alpha-1}$ then it follows from
{\bf R 1}${}_{\alpha-1}$ that ${\cal M}^{\cal T}_{\alpha-1}$ can't be ill-founded.
Also, if ${\cal T}$ has limit length and ${\cal D}^{\cal T} \cap b$ is infinite where
$b$ is the unique cofinal branch through $T$, then the indices of the target models of
$\sigma_i$ for $i \in b$ yield an infinite descending sequence of ordinals. We have
reached a contradiction!

\bigskip
\hfill $\square$ (\ref{normal-iterability})

\bigskip
We can now exploit the previous argument a bit further and arrive at the following.

\begin{lemma}\label{iterability}
Suppose that ${\frak C}_n({\cal N}_\xi)$ exists for some $\xi < \infty$ and $n <
\omega$. Then ${\frak C}_n({\cal
N}_\xi)$ is $n$-iterable.
\end{lemma}

{\sc Proof sketch.} This time, we have to deal with a sequence of iteration trees
$${\cal T}_0 {}^\frown {\cal T}_1 {}^\frown {\cal T}_2 {}^\frown ...$$ of length $\leq
\omega$, where ${\cal T}_0$ is on ${\cal M}$, and ${\cal T}_{i+1}$ is on the last
model of ${\cal T}_i$. We may now repeatedly apply
the proof of \ref{normal-iterability} to see that
the last model of any ${\cal T}_i$ can be embedded into some ${\cal N}_\xi$. This will
give the desired conclusion.

\bigskip
\hfill $\square$ (\ref{iterability})

\bigskip
It is now easy to see that ${\cal M}_{{\rm OR}} = {\cal N}_{{\rm OR}}$ is a model of height ${\rm OR}$.

\begin{defn}
We write $K^c$ for ${\cal M}_{{\rm OR}} = {\cal N}_{{\rm OR}}$. $K^c$ is called the countably
complete core model below $\zerohandgrenade$.
\end{defn}

The method of the proof of \ref{iterability} also yields the following.

\begin{lemma}\label{unique-next-ext}
Let ${\cal M}_\xi = (J_\alpha[{\vec E}];\in,{\vec E})$ be passive. Then there is at
most one countably complete extender $F$ such that $(J_\alpha[{\vec E}];\in,{\vec
E},F)$ is a premouse.
\end{lemma}

{\sc Proof sketch.} Suppose $F$, $F'$ are both countably complete and such that
$(J_\alpha[{\vec E}];\in,{\vec
E},F)$ as well as $(J_\alpha[{\vec E}];\in,{\vec
E},F')$ is a premouse. We may then form the prebicephalus
$${\cal B} = (J_\alpha[{\vec E}];\in,{\vec
E},F,F').$$ (Cf. \cite[\S 9]{FSIT} for a paradigmatic theory of bicephali.)
The proof of \ref{normal-iterability} shows that
${\cal B}$ is $0$-iterable in the
obvious sense. By coiterating ${\cal B}$ against itself we may then as usual conclude
that in fact $F = F'$.

\bigskip
\hfill $\square$ (\ref{unique-next-ext})

\bigskip
Lemma \ref{unique-next-ext}
may also be derived from our more general theory of bicephali
which we shall develop in section 5 below.

Notice that in Case 1 of the recursive construction of ${\cal N}_{\xi+1}$ we had $F$
required to be unique. \ref{unique-next-ext} now says that this requirement poses no
restriction at all.

\bigskip
For the purposes of isolating $K$ we shall in fact need a variant of $K^c$.
Let $\Gamma \not= \emptyset$ be a class of regular cardinals, $\omega \notin \Gamma$.
We then define ${\cal N}_\xi$ and ${\cal M}_\xi$ exactly as before except that in Case
1 of the definition of ${\cal N}_{\xi+1}$ we also require that ${\rm cf}^V({\rm c.p.}(F)) \in
\Gamma$ (here, $F$ is the top extender of ${\cal N}_{\xi+1}$).
Let us write ${\cal N}_\xi^\Gamma$ and ${\cal M}_\xi^\Gamma$ instead of
${\cal N}_\xi$ and ${\cal M}_\xi$. The proofs of the current section show that
they all exist, and that ${\frak C}_n({\cal N}_\xi^\Gamma)$ is always $n$-iterable.

\begin{defn}\label{Kc-Gamma}
Let $\Gamma \not= \emptyset$ be a class of regular cardinals, $\omega \notin \Gamma$.
We write $K^c_\Gamma$ for ${\cal M}_{{\rm OR}}^\Gamma = {\cal N}_{{\rm OR}}^\Gamma$.
$K^c$ is called the (countably
complete) $\Gamma$-full core model below $\zerohandgrenade$.
\end{defn}

\section{Universality of $K^c$.}

\begin{defn}\label{defn-full}
Let $W$ be a weasel. $W$ is called full provided
the following holds.
Let $\alpha < \beta$ be any cardinals of $W$.
There are no ${\tilde \tau} > \beta$ and a
countably complete ${\tilde F} \colon {\cal P}(\alpha) \cap W \rightarrow
{\cal P}(\beta) \cap {\cal J}^W_{\tilde \tau}$ such that
$({\cal J}^W_{\tilde \tau},{\tilde F})$ is a pre-premouse,
${\rm c.p.}({\tilde F}) = \alpha$, and
${\tilde F}(\alpha) = \beta$.
\end{defn}

We shall from now on denote by ``$\lnot \ \zerohandgrenade$'' the assumption that
every premouse is below $\zerohandgrenade$.

\begin{lemma}\label{weak-max} $( \ \lnot \ \zerohandgrenade \ )$
$K^c$ is full.
\end{lemma}

{\sc Proof.} Suppose otherwise.
Let $\alpha < \beta$ be cardinals of $K^c$,
let ${\tilde \tau} > \beta$ and let 
${\tilde F} \colon {\cal P}(\alpha) \cap K^c \rightarrow
{\cal P}(\beta) \cap {\cal J}_{\tilde \tau}^{K^c}$
be countably complete and
such that
$({\cal J}_{\tilde \tau}^{K^c},{\tilde F})$ is a pre-premouse,
${\rm c.p.}({\tilde F}) = \alpha$ and
${\tilde F}(\alpha) = \beta$.

Let $\lambda \leq \beta$ be least such
that ${\tilde F}(f)(\xi)
< \lambda$ whenever $f \in {}^\alpha\alpha \cap K^c$
and $\xi < \lambda$. Set $F = {\tilde F} | \lambda$,
and $\tau = \alpha^{+K^c}$.

There is $\pi \colon {\rm Ult}({\cal J}^{K^c}_\tau,F) \rightarrow {\cal
J}^{K^c}_{\tilde \tau}$ defined by $\pi(F(f)(a)) = {\tilde F}(f)(a)$,
and of course $\pi {\upharpoonright} \lambda = {\rm id}$ and ${\tilde F} =
\pi \circ F$. By the choice of
$\lambda$, we also easily get $F(\alpha) = \lambda$.

As $\lambda$ is a limit cardinal in ${\rm Ult}({\cal J}^{K^c}_\tau,F)$,
by
$\pi
{\upharpoonright} \lambda = {\rm id}$ we then get that in fact $\lambda$ is a
(limit) cardinal in $K^c$. Also, $\lambda$ is the largest cardinal in
${\rm Ult}({\cal J}^{K^c}_\tau,F)$, so that applying the condensation lemma
\cite[Lemma 4]{NFS}
to cofinally many restrictions of $\pi$ gives that
actually ${\rm Ult}({\cal J}^{K^c}_\tau,F) = {\cal J}^{K^c}_\gamma$ for some
$\gamma$.

Notice that $({\cal J}^{K^c}_\gamma,F)$ is a premouse. (The initial
segment condition for $F$ is vacuously true by the choice of
$\lambda$.) Moreover, we have that ${\cal J}^{K^c}_\lambda = {\cal
M}_\eta$ for some $\eta$, as $\lambda$ is a cardinal in $K^c$.
It is also straightforward to check that because $\lambda$ is the
largest cardinal of ${\cal J}^{K^c}_\gamma$ there is $\eta' >
\eta$ with ${\cal J}^{K^c}_\gamma = {\cal M}_{\eta'}$.

We thus have ${\cal M}_{\eta'+1} = ({\cal J}^{K^c}_\gamma,F)$
by \ref{unique-next-ext}, so
that in particular $\rho_\omega({\cal M}_{\eta'+1}) <
\lambda$. But
the fact that $\lambda$ is a cardinal in $K^c$ implies that
$\rho_\omega({\cal M}_{\tilde \eta})
\geq \lambda$ for all ${\tilde
\eta} \geq \eta$. Contradiction!

\bigskip
\hfill $\square$ (\ref{weak-max})

\begin{defn}\label{defn-universal}
Let $W$ be a weasel. $W$ is called universal if $W$ is iterable,
and whenever ${\cal T}$, ${\cal U}$ are iteration
trees arising from the coiteration of $W$
with some (set- or class-sized) premouse ${\cal M}$ such that ${\rm lh}({\cal T}) =
{\rm lh}({\cal U}) = {\rm OR}+1$, then ${\cal M}$ is a weasel, ${\cal D}^{\cal U} \cap (0,{\rm OR}]_U =
\emptyset$, $\pi^{\cal U}_{0 \infty} {\rm " } {\rm OR} \subset {\rm OR}$, and
${\cal M}^{\cal U}_\infty \trianglelefteq {\cal M}^{\cal T}_\infty$.
\end{defn}

\begin{lemma}\label{W-is-universal}
$( \ \lnot \ \zerohandgrenade \ )$ Every
full weasel is universal.
\end{lemma}

{\sc Proof.} This is shown by varying an argument which is due to Jensen, cf.
the Addendum to \cite[\S 3 Theorem 5]{NFS}.

Fix a full weasel $W$, and suppose that $W$ is not universal.
Let ${\cal N}$ witness this, i.e.,
${\cal N}$ is a (set or proper class sized) premouse, and
if $({\cal T},{\cal U})$ denotes the coiteration of $W$ with ${\cal
N}$ then both ${\cal T}$ and ${\cal U}$ have
length ${\rm OR} +1$ and there is a club $C \subset {\rm OR}$
of strong limit cardinals
such that $D^{\cal U} \cap [{\rm min}(C),{\rm OR}]_U = \emptyset$, 
$\pi^{\cal U}_{\alpha \infty}
{\upharpoonright} \alpha = {\rm id}$, and
$\pi^{\cal U}_{\alpha \beta}(\alpha) = \beta$
for all $\alpha \leq \beta \in C$, and $D^{\cal T} \cap
[0,{\rm OR}]_T = \emptyset$ and $\pi^{\cal T}_{0 \infty} {\rm " } \beta
\subset \beta$ for all $\beta \in C$.

\bigskip
{\it Case 1.} There is some $\beta \in C$ such that $\pi^{\cal T}_{0 \beta}(\beta) >
\beta$.

\bigskip
Fix such $\beta \in C$. As $\pi^{\cal T}_{0 \beta} {\rm " } \beta
\subset \beta$ but $\pi^{\cal T}_{0 \beta}(\beta) >
\beta$, we must have that $\mu = {\rm cf}^W(\beta) < \beta$ is measurable in $W$ and that
$\mu$ is used in the iteration giving $\pi^{\cal T}_{0 \beta}$. Let $\alpha$ be least
in $[0,\beta)_T$ such that ${\tilde \mu} = \pi^{\cal T}_{0 \alpha}(\mu)$ is the
critical point of $\pi^{\cal T}_{\alpha \beta}$. Then $${\rm cf}^{{\cal M}_\alpha^{\cal
T}}(\beta) = {\tilde \mu}.$$
Pick $f \in {\cal M}_\alpha^{\cal
T}$, $f \colon {\tilde \mu} \rightarrow \beta$ cofinal. Then $\pi_{\alpha {\rm OR}}^{\cal
T}(f) {\upharpoonright} {\tilde \mu} \colon {\tilde \mu} \rightarrow \beta$ is cofinal,
too. (For $\xi < {\tilde \mu}$ do we have $\pi_{\alpha \infty}^{\cal
T}(f)(\xi) = \pi_{\alpha \infty}^{\cal
T}(f(\xi)) < \beta$, as $\pi^{\cal T}_{0 \infty} {\rm " } \beta
\subset \beta$; and of course $\pi_{\alpha \infty}^{\cal
T}(f(\xi)) \geq f(\xi)$.) But $\pi_{\alpha \infty}^{\cal
T}(f) {\upharpoonright} {\tilde \mu} \in {\cal M}_{\infty}^{\cal T}$, and thus $\beta$ is
singular in ${\cal M}_{\infty}^{\cal T}$.

On the other hand, $\beta$ is of course inaccessible in ${\cal M}_{\infty}^{\cal U}$,
and ${\cal M}_{\infty}^{\cal U} \trianglerighteq {\cal M}_{\infty}^{\cal T}$. This
is a contradiction!

\bigskip
We may hence assume that:

\bigskip
{\it Case 2.} All $\beta \in C$ are such that $\pi^{\cal T}_{0 \beta}(\beta) =
\beta$.

\bigskip
{\bf Claim 1.} There are a club class $D' \subset {\rm OR}$ and a commutative
system $(\sigma_{\alpha \beta} \ \colon \ \alpha \leq \beta \in D')$ such
that
$\sigma_{\alpha \beta} \ \colon \ {\cal J}^W_{\alpha^{+W}}
\rightarrow_{\Sigma_0} {\cal J}^W_{\beta^{+W}}$ cofinally with
$\sigma_{\alpha \beta} {\upharpoonright} \alpha = {\rm id}$ and
$\sigma_{\alpha \beta}(\alpha) = \beta$ for all $\alpha \leq \beta \in D'$.

\bigskip
{\sc Proof.}
Let us write $W_\alpha = {\cal J}^W_{\alpha^{+W}}$ and
$W_\alpha^\star = \pi^{\cal T}_{0 \alpha}(W_\alpha)$
for $\alpha \in C$. Notice that $$W_\alpha^\star =
{\cal J}^{{\cal M}^{\cal
U}_\alpha}_{\alpha^+} {\rm \ where \ } \alpha^+ =
\alpha^{+{\cal M}^{\cal U}_\alpha} {\rm , }$$
so that we have that $$\pi^{\cal U}_{\alpha \beta} {\upharpoonright}
W_\alpha^\star \ \colon \ W_\alpha^\star
\rightarrow_{\Sigma_0} W_\beta^\star {\rm \
cofinally, }$$ for $\alpha \leq \beta \in C$. But
we clearly also have that
$$\pi^{\cal T}_{0 \beta} {\upharpoonright} W_\beta \ \colon \
W_\beta \rightarrow_{\Sigma_0} W_\beta^\star {\rm \ cofinally}$$
for $\beta \in C$.
We aim to show that typically ${\rm ran}(\pi^{\cal U}_{\alpha \beta}
{\upharpoonright}
W_\alpha^\star \circ \pi^{\cal T}_{0 \alpha}
{\upharpoonright} W_\alpha) \subset {\rm ran}(\pi^{\cal T}_{0 \beta}
{\upharpoonright} W_\beta)$, so that $(\pi^{\cal T}_{0 \beta}
{\upharpoonright} W_\beta)^{-1} \circ \pi^{\cal U}_{\alpha \beta}
{\upharpoonright}
W_\alpha^\star \circ \pi^{\cal T}_{0 \alpha}
{\upharpoonright} W_\alpha$ makes sense.

Let $\tau_\alpha = {\rm OR} \cap W_\alpha$, $\alpha \in C$.
The point is that ${\rm cf}(\tau_\alpha) = {\rm cf}(\tau_\beta)$ for all
$\alpha$, $\beta \in C$. Let $\gamma = {\rm cf}(\tau_\alpha)$, $\alpha
\in C$, and pick $X_\alpha = \{ \xi_0^\alpha < \xi_1^\alpha <
... < \xi_i^\alpha < ... \colon i < \gamma \}
\subset \tau_\alpha$ unbounded in
$\tau_\alpha$ and of order type $\gamma$ for all $\alpha \in C$.
We may define a regressive function $\delta \ \colon \ C \cap \{
\alpha \colon {\rm cf}(\alpha) > \gamma \} \rightarrow {\rm OR}$ by letting
$\delta(\alpha)$ be the least $\delta < \alpha$ such that
$$\pi^{\cal T}_{0 \alpha} {\rm " } X_\alpha \subset
{\rm ran}(\pi^{\cal U}_{\delta \alpha}).$$ We may further define a regressive
function $h \ \colon \ C \cap \{
\alpha \colon {\rm cf}(\alpha) > \gamma \} \rightarrow {\rm OR}$ by letting
$$h(\alpha) = (\pi^{\cal U}_{\delta(\alpha) \alpha})^{-1} {\rm " }
\pi^{\cal T}_{0 \alpha} {\rm " } X_\alpha.$$
By Fodor, there are an ordinal $\delta_0$, some $Y \subset {\rm OR}$, and
some unbounded $D \subset C$ with $\delta(\alpha) =
\delta_0$ and $h(\alpha) = Y$, all $\alpha \in D$.
We then have that
$$\pi^{\cal U}_{\alpha \beta} \circ \pi^{\cal T}_{0
\alpha}(\xi^\alpha_i) = \pi^{\cal T}_{0 \beta}(\xi^\alpha_i)$$ for all
$\alpha \leq \beta \in D$ and $i < \gamma$.

%
We now verify that $$\pi^{\cal U}_{\alpha \beta}
{\upharpoonright}
W_\alpha^\star \circ \pi^{\cal T}_{0 \alpha}
{\upharpoonright} W_\alpha \subset \pi^{\cal T}_{0 \beta}
{\upharpoonright} W_\beta$$ for all $\alpha \leq \beta \in D$.

Well, let $\alpha \leq \beta \in D$ and $x \in W_\alpha$. Then $x \in
{\cal J}^W_{\xi_i^\alpha}$ for some $i < \gamma$. Let $f \in W_\alpha$
be the least (in the order of constructibility) surjection $f \colon
\alpha \rightarrow {\cal J}^W_{\xi^\alpha_i}$. So $x = f(\nu)$ , some
$\nu < \alpha$. We have that $\pi^{\cal U}_{\alpha \beta} \circ
\pi^{\cal T}_{0 \alpha}(f) =$ the least surjection $$g \colon \beta
\rightarrow {\cal J}^{W_\beta^\star}_{\pi^{\cal U}_{\alpha \beta}
\circ
\pi^{\cal T}_{0 \alpha}(\xi^\alpha_i)} =
{\cal J}^{W_\beta^\star}_{\pi^{\cal T}_{0 \beta}(\xi^\beta_i)} {\rm ,
}$$
so that in particular $\pi^{\cal U}_{\alpha \beta} \circ
\pi^{\cal T}_{0 \alpha}(f) \in {\rm ran}(\pi^{\cal T}_{0 \beta})$, say
$\pi^{\cal U}_{\alpha \beta} \circ
\pi^{\cal T}_{0 \alpha}(f) = \pi^{\cal T}_{0 \beta}(g')$.

We then get $\pi^{\cal U}_{\alpha \beta} \circ
\pi^{\cal T}_{0 \alpha}(x) = \pi^{\cal U}_{\alpha \beta} \circ
\pi^{\cal T}_{0 \alpha}(f(\nu)) =
\pi^{\cal T}_{0 \beta}(g')(\pi^{\cal T}_{0 \alpha}(\nu))$, as $\pi^{\cal T}_{0
\alpha}{\rm " }\alpha \subset \alpha$ and $\pi^{\cal U}_{\alpha \beta}
{\upharpoonright} \alpha = {\rm id}$,
$= \pi^{\cal T}_{0 \beta}(g')(\pi^{\cal T}_{0 \beta}(\nu))$, as $\pi^{\cal T}_{\alpha
\beta} {\upharpoonright} \alpha = {\rm id}$,
$= \pi^{\cal T}_{0 \beta}(g'(\nu))$.

Now let $$D' = D \cup \{ \alpha \ \colon \ 
\exists \lambda \in {\rm Lim} \ \alpha = {\rm sup}(D \cap \lambda) \}.$$ 
$D'$ is club. 
In 
order to finish the proof of Claim 1
it suffices to verify the following.

\bigskip
{\bf Subclaim.} If $\{ \beta < \alpha \ \colon \ \beta \in D \}$ is unbounded in
$\alpha$, and if $({\tilde W},\sigma_\beta \ \colon \ \beta \in D \cap \alpha)$ is the
direct limit of the system $$(W_\beta,((\pi^{\cal T}_{0 \beta})^{-1} \circ \pi^{\cal
U}_{\gamma \beta} \circ \pi^{\cal T}_{0 \gamma} \ \colon \ \gamma \leq \beta \in D
\cap \alpha)$$ then we have that $\sigma_\beta \colon W_\beta \rightarrow W_\alpha$ is
cofinal with $\sigma_\beta {\upharpoonright} \beta = {\rm id}$ and $\sigma_\beta(\beta) =
\alpha$.

\bigskip
{\sc Proof.} Let $\delta = {\rm min}(D \setminus \alpha+1)$. There is $\sigma \colon {\tilde
W} \rightarrow W_\delta$ defined by $$x \mapsto (\pi^{\cal T}_{0 \delta})^{-1} \circ
\pi^{\cal U}_{\beta \delta} \circ \pi^{\cal T}_{0 \beta}({\bar x}) {\rm \ where \ }
{\bar x} = \sigma_\beta^{-1}(x).$$ It is clear that $\sigma {\upharpoonright} \alpha =
{\rm id}$, and $\sigma(\alpha) = \delta$. By \cite[\S 8 Lemma 4]{NFS} 
we hence have that
(${\tilde W}$ is transitive and) ${\tilde W} \trianglelefteq W_\alpha$.

Now consider the two maps $\pi_{0 \alpha}^{\cal T} \colon W_\alpha \rightarrow
W_\alpha^\star$ and ${\tilde \sigma} \colon {\tilde W} \rightarrow W_\alpha^\star$,
where $${\tilde \sigma}(x) = \pi_{\beta \alpha}^{\cal U} \circ \pi_{0 \beta}^{\cal
T}({\bar x}) {\rm \ for \ } {\bar x} = \sigma^{-1}_\beta(x).$$
We may define $\tau \colon W_\alpha^\star \rightarrow \pi^{\cal T}_{0 \alpha}({\tilde
W})$ by setting $$\tau({\tilde \sigma}(f)(\xi)) = \pi^{\cal T}_{0 \alpha}(f)(\xi)$$
for $\xi<\alpha$ and appropriate $f$ with ${\rm dom}(f)$ bounded below $\alpha$. $\tau$ is
easily seen to be well-defined, and in fact surjective.

But then $\pi_{0 \alpha}^{\cal T}({\tilde W}) = W_\alpha^\star$, and so
$W_\alpha = {\tilde W}$.

\bigskip
\hfill $\square$ (Subclaim)

\hfill $\square$ (Claim 1)

\bigskip
Now fix $D'$ and maps $\sigma_{\alpha \beta}$ as in Claim 1.  
For each $\alpha \in D'$, let $\alpha' = {\rm min} (D' \setminus \alpha+1)$,
and let $F_\alpha$
the extender derived from $\sigma_{\alpha \alpha'}$.

\bigskip
{\bf Claim 2.} There is some $\beta \in D'$ such that $F_\beta$ is
countably complete.

\bigskip
{\sc Proof.} Let $D'' = \{ \alpha \in D' \ \colon \ D' \cap \alpha$ is unbounded in
$\alpha \}$. Notice that $D''$ is club. Suppose that 
$F_\alpha$ is not countably complete for any 
$\alpha \in D''$. Pick for any $\alpha \in D''$ sequences
$(a^\alpha_n \colon n<\omega)$ and $(X^\alpha_n \colon n<\omega)$ 
witnessing this. This means that 
$a^\alpha_n
\in F_\alpha(X^\alpha_n)$ for all natural numbers $n$, 
but there is no orderpreserving function $\tau \colon
\bigcup_{n<\omega} a^\alpha_n \rightarrow \alpha = {\rm c.p.}(F_\alpha)$ with
the property that for all natural numbers $n$, we have
$\tau {\rm " } a^\alpha_n \in X^\alpha_n$.

Let $\alpha \in D''$. Let $g^\alpha \colon \bigcup_{n<\omega} a^\alpha_n \rightarrow
{\rm otp}(\bigcup_{n<\omega} a^\alpha_n) < \omega_1$ denote the transitive collapse, and
let $$g(\alpha) = ({\rm otp}(\bigcup_{n<\omega} a^\alpha_n),(g^\alpha {\rm " } a_n^\alpha
\colon n<\omega)) {\rm , }$$ i.e., $g(\alpha)$ tells us how the $a_n^\alpha$ sit
inside $\bigcup_{n<\omega} a^\alpha_n$. Let $\beta(\alpha)$ be the least $\beta \in
D'$ such that 
for all natural numbers $n$ do we have that
$X_n^\alpha \in {\rm ran}(\pi_{\beta \alpha})$. Notice that
$\beta(\alpha) < \alpha$, as ${\rm cf}(\alpha) > \omega$.

We may now apply Fodor's lemma to the function $F$ with ${\rm dom}(F) = D''$ defined by
$$F(\alpha) = (\beta(\alpha),(\pi_{\beta(\alpha) \alpha}^{-1}(X_n^\alpha) \colon
n<\omega),g(\alpha))$$ to get some unbounded $E \subset D''$ on which $F$ is constant.
Let $\alpha < \gamma \in E$. As $g(\gamma) = g(\alpha)$, we have that $\tau =
(g^\alpha)^{-1} \circ g^\gamma \colon \bigcup_{n<\omega} a^\gamma_n \rightarrow
\bigcup_{n<\omega} a^\alpha_n$ is order preserving with $\tau {\rm
" } a^\gamma_n = a^\alpha_n$ for all natural numbers $n$. But we now get, by $F(\gamma) = F(\alpha)$ together with
the commutativity of the maps $\pi_{\beta \beta'}$, that
for all $n<\omega$,
$$\tau {\rm
" } a^\gamma_n = a^\alpha_n \in F_\alpha(X_n^\alpha) = F_\alpha(\pi_{\alpha
\gamma}^{-1}(X_n^\gamma)) = X_n^\gamma \cap [F_\alpha(\alpha)]^{< \omega}
{\rm , }$$ which contradicts the choice of
$(a_n^\gamma,X_n^\gamma,n<\omega)$.

\bigskip
\hfill $\square$ (Claim 2)

\bigskip
But Claim 2 now plainly contradicts
the fact that $W$ is supposed to be full.

\bigskip
\hfill $\square$ (\ref{W-is-universal})

\bigskip
The following is an immediate corollary to \ref{weak-max} and \ref{W-is-universal}.

\begin{cor}\label{Kc-is-universal}
$( \ \lnot \ \zerohandgrenade \ )$ $K^c$ is universal.
\end{cor}

Let $\Gamma \not= \emptyset$ be a class of regular cardinals, $\omega \notin
\Gamma$. In the verification of Claim 2 during the proof of
\ref{W-is-universal} we might have
considered the stationary class $D'' \cap \Gamma$ rather than $D''$, thus deducing
that there is some $\beta \in D'$ such that $F_\beta$ is countably complete and
${\rm cf}^V({\rm c.p.}(F_\beta)) \in \Gamma$. This leads to the following.

\begin{cor}\label{KcGamma-is-universal}
$( \ \lnot \ \zerohandgrenade \ )$
Let $\Gamma \not= \emptyset$ be a class of regular cardinals, $\omega \notin
\Gamma$. Then $K^c_\Gamma$ is universal.
\end{cor}

It is however not true
that an initial segment of $K^c$
whose height is a cardinal in $V$ greater than $\aleph_1$
is universal for coiterable premice of at
most the same 
height.\footnote{Compare with \cite[Theorem 3.4]{maximality-paper}.}  
Anticipating the theory of $K$,
here is an example. Suppose that $K$ has a measurable cardinal, and let $\mu$ be the
least one. Let $\lambda > \mu$ be, in $K$, a singular cardinal of cofinality $\mu$.
Suppose that $V = K^{Col(\omega,\mu)}$. Then ${\cal J}^{K^c}_{\lambda^+}$ does not
win the coiteration against ${\cal J}^{K}_{\lambda^+}$, because if $F$ is an extender on $K$ with critical
point $\mu$ then ${\rm Ult}_0({\cal J}^{K}_{\lambda^+};F)$ will have height $> \lambda^+$.
We leave the further (easy) details to the reader.

\begin{lemma}\label{goodness} (Goodness of $K^c$.)
Let $\Gamma \not= \emptyset$ be a class of regular cardinals, $\omega \notin \Gamma$.
Let $\kappa$ be a cardinal of $K^c_\Gamma$, and let ${\cal P}$ be an
iterate
of $K^c_\Gamma$ above $\kappa$, i.e., there is an iteration tree
${\cal T} = (({\cal M}^{\cal T}_\alpha,\pi^{\cal T}_{\alpha \beta} \colon 
\alpha T \beta \leq \theta+1),
(E^{\cal T}_\alpha \ \colon \ \alpha < \theta),T)$ on $K^c_\Gamma$
such that ${\rm c.p.}(E^{\cal T}_\alpha) \geq \kappa$ for 
every $\alpha < \theta$
and ${\cal P} = {\cal M}^{\cal T}_\theta$.
Let $F = E^{\cal P}_\alpha \not= \emptyset$ be such that $\alpha 
> \kappa$, and $\mu =
{\rm c.p.}(F) < \kappa$ (notice that
$\alpha > \mu^{+K^c_\Gamma}$ and $F$ is total on
${\cal P}$).

Then $F$ is countably complete.
\end{lemma}

{\sc Proof.} This is shown by revisiting the arguments for
\ref{normal-iterability} and \ref{iterability}.
Suppose without loss of generality that $K^c_\Gamma = K^c$.
Let $((a_n , X_n) \colon n < \omega)$ be given such
that $$a_n \in F(X_n) {\rm \ \ for \ all \ \ } n < \omega.$$
Let $\lambda^+$, a successor cardinal in $V$,
be large enough such that ${\cal T}$ and 
$F$ are both hereditarily smaller than $\lambda^+$. We may then construe
${\cal T}$ as an iteration tree on ${\cal J}_{\lambda^+}^{K^c}$, and we shall have that $F$  
is on the final model of ${\cal T}$ so construed.
Let $\xi$ be such that ${\cal J}^{K^c}_{\lambda^+} =
{\cal M}_\xi$.

Running the proofs of \ref{normal-iterability} and \ref{iterability},
we may then pick the map
$$\sigma \ \colon \ {\cal M} \rightarrow {\cal M}_\xi {\rm , }$$
such that in fact $\sigma$ is the restriction of an uncollapse
(also
called $\sigma$) of a countable submodel of a large enough
initial segment of $V$ containing all objects of current interest; in particular, we
want that both $\{ a_n \ \colon \ n < \omega \}$ and
$\{ X_n \ \colon \ n < \omega \}$ are contained in
${\rm ran}(\sigma)$.

In the end, we get an embedding
$$\sigma' \ \colon {\cal M}' \rightarrow {\cal M}_{\xi'}$$ for some
$\xi' \leq \xi$,
where ${\cal M}'$ is the final model of $\sigma^{-1}({\cal T})$.
Moreover, as ${\cal T}$ is above $\kappa$, $\sigma^{-1}({\cal T})$ is
above $\sigma^{-1}(\kappa)$, which implies that
$$\sigma' {\upharpoonright} \sigma^{-1}(\kappa) =
\sigma {\upharpoonright} \sigma^{-1}(\kappa).$$

As ${\rm c.p.}(F) < \kappa$, $\sigma' \circ \sigma^{-1}(F)$ is
countably complete. We may thus pick an order-preserving $\rho$ such that
$$\rho {\rm " } \sigma' \circ \sigma^{-1}(a_n)
\in \sigma' \circ \sigma^{-1}(X_n)
{\rm \ \ for \ all \ } n < \omega.$$

But $\sigma' \circ \sigma^{-1}(X_n) = X_n$, as $\mu < \kappa$ and
$\sigma' {\upharpoonright} \sigma^{-1}(\kappa) =
\sigma {\upharpoonright} \sigma^{-1}(\kappa)$, which means that
$$\rho {\rm " } \sigma' \circ \sigma^{-1}(a_n)
\in X_n
{\rm \ \ for \ all \ } n < \omega {\rm , }$$
and $\rho \circ \sigma' \circ \sigma^{-1}$ witnesses that $F$ is
countably complete with respect to $((a_n , X_n) \colon n < \omega)$.

\bigskip
\hfill $\square$ (\ref{goodness})

%
%
%

\section{Generalized bicephali.}
\setcounter{equation}{0}

We let $C_0$ denote the class of all limit cardinals $\kappa$
of $V$ such that ${\cal J}_\kappa^{K^c}$ is universal for coiterable premice of height
$< \kappa$. I.e., if $\kappa \in C_0$, ${\cal M}$ is a premouse with ${\cal M} \cap {\rm OR}
< \kappa$, and ${\cal T}$, ${\cal U}$ are the iteration trees arising from the
successful
comparison of ${\cal J}_\kappa^{K^c}$ with ${\cal M}$, then ${\cal D}^{\cal U} \cap
(0,\infty]_U = \emptyset$, and ${\cal M}^{\cal U}_\infty \triangleleft
{\cal M}^{\cal T}_\infty$.

\begin{lemma}\label{C0-is-club} $( \ \lnot \ \zerohandgrenade \ )$
$C_0$ is closed unbounded in ${\rm OR}$.
\end{lemma}

{\sc Proof.} $C_0$ is trivially be seen to be closed. Suppose that $C_0 \subset \eta$
for some $\eta \in {\rm OR}$. We may then define $F \colon \{ \kappa \ \colon \ \kappa$ is
a limit cardinal $\} \setminus \eta \rightarrow {\rm OR}$
by $F(\kappa) =$ the least $\alpha$ such that there is a coiterable premouse ${\cal
M}$ of height
$\alpha < \kappa$ and ${\cal M}$ wins the coiteration against
${\cal J}_\kappa^{K^c}$. By Fodor, there is an {\it unbounded} class $A \subset {\rm OR}$
such that $F {\rm " } A = \{ \alpha \}$ for some $\alpha \in {\rm OR}$. There are at most
$2^{{\rm Card}(\alpha)}$ many premice of height $\alpha$, so that
by \ref{Kc-is-universal}
there is
some cardinal $\gamma$ such that ${\cal J}^{K^c}_\gamma$
wins the coiteration against all
of them which are coiterable with $K^c$. This gives a contradiction!

\bigskip
\hfill $\square$ (\ref{C0-is-club})

\begin{lemma}\label{overlaps} $( \ \lnot \ \zerohandgrenade \ )$
Let $\kappa \in C_0$, and
let ${\cal M} \trianglerighteq
{\cal J}^{K^c}_\kappa$ be a $0$-iterable premouse. Let $F =
E^{\cal M}_\nu \not= \emptyset$ be
such that $\nu \geq \kappa$ and ${\rm c.p.}(F) < \kappa$. Then $F$ is countably complete.
\end{lemma}

{\sc Proof.} Let $\kappa$, 
${\cal M}$, $F$, and $\nu$ be as in the statement of \ref{overlaps}.
Set $\mu = {\rm c.p.}(F)$.
Let $((a_n,X_n) \colon n<\omega)$ be such that $a_n \in [F(\mu)]^{< \omega}$,
$X_n \in {\cal P}([\mu]^{{\rm Card}(a_n)}) \cap {\cal M}$, and
$a_n \in F(X_n)$ for every
$n<\omega$. We aim to find an order-preserving
$\tau \colon \bigcup_{n<\omega} a_n \rightarrow \mu$
such that $\tau {\rm " } a_n \in X_n$ for every $n<\omega$.
Let $$\sigma \colon {\bar {\cal M}} \rightarrow {\cal M}$$ be fully
elementary such that $\sigma {\upharpoonright} \mu^{+{\cal M}}+1 ={\rm id}$,
${\rm Card}({\bar {\cal M}}) < \mu^{++}$, and
$\{ a_n \colon n<\omega \} \subset {\rm ran}(\sigma)$. Let ${\bar F} = \sigma^{-1}(F)$ in
case $F \in {\cal M}$ (where we then assume without loss of generality that $F \in {\rm ran}(\sigma)$), and let
${\bar F}$ be the top extender of ${\bar {\cal M}}$ otherwise (i.e.,
if $F$ is the top extender of
${\cal M}$).

Let ${\cal T}_0$ and ${\cal T}_1$ be $0$-maximal iteration trees on 
the phalanx $$(({\cal M},{\bar {\cal M}}),{\rm c.p.}(\sigma))$$ and 
on ${\cal M}$, respectively, stemming from the coiteration of that phalanx with
${\cal M}$. 
Also let ${\cal U}$ and ${\cal T}$ be the
$0$-maximal iteration trees
arising from the comparison of ${\bar {\cal M}}$ with ${\cal M}$, respectively.
Let ${\cal M}^*$ be the last model of ${\cal T}_0$. Then
${\bar {\cal M}} T_0 {\cal M}^*$ by
\cite[\S 8 Lemma 1]{NFS}. But this implies that ${\cal U}$ only uses extenders
whose critical points are less than or equal to ${\rm c.p.}(\sigma)$, by 
\ref{phalanx-iterations}.

Obviously, the
coiteration of ${\bar {\cal M}}$ with ${\cal J}^{K^c}_\kappa$
is an initial segment
of the coiteration ${\cal U}$, ${\cal T}$.
Moreover, as $\kappa \in C_0$ and ${\rm Card}({\bar {\cal M}}) < \kappa$,
${\cal J}^{K^c}_\kappa$ wins the coiteration against
${\bar {\cal M}}$.
Hence ${\cal T}$ actually only uses extenders from ${\cal J}^{K^c}_\kappa$ and its images, and
${\cal U}$ does not drop.

This means that the main branch of ${\cal U}$
gives us an embedding $$\pi \colon {\bar {\cal M}} \rightarrow_{\Sigma_0}
{\cal M}_\infty^{\cal U} {\rm , }$$ where ${\cal M}_\infty^{\cal U}$
is an initial segment of ${\cal M}^{\cal
T}_\infty$. Notice that $\pi {\upharpoonright} \mu^{+{\cal M}}+1 = {\rm id}$ as
${\cal U}$ is above ${\rm c.p.}(\sigma)$. Let ${\tilde F} = \pi({\bar F})$ in case ${\bar F}
\in {\bar {\cal M}}$, and let ${\tilde F}$ be the top extender of ${\cal M}^{\cal
U}_\infty$
otherwise.

We claim that ${\tilde F}$
is countably complete.
This trivially follows from
\ref{goodness} if ${\cal T}$ is above $\mu^{+{\cal M}}$.
Otherwise let $E^{{\cal M}_i^{\cal
T}}_\nu$ be the first extender used on ${\cal T}$ with critical point $\leq \mu$.
By Claim 1 in the proof of \ref{normal-linear},
$\nu$ is then strictly greater than the index of
${\tilde F}$. But then \ref{goodness} applied to ${\cal T} {\upharpoonright} i$ yields
that ${\tilde F}$ is countably complete.

Notice that $\pi \circ \sigma^{-1}(a_n) \in {\tilde F}(X_n)$ for all $n<\omega$.
Now let ${\bar \tau} \colon \bigcup_{n<\omega}
{\pi \circ \sigma^{-1}}(a_n) \rightarrow \mu$ be order-preserving and
such that ${\bar \tau} {\rm " }
{\pi \circ \sigma^{-1}}(a_n) \in X_n$ for all $n<\omega$.
Putting $\tau = {\bar \tau} \circ \pi \circ \sigma^{-1}$ hence gives us a
witness as desired.

\bigskip
\hfill $\square$ (\ref{overlaps})

\bigskip
We now have to turn towards our theory of ``bicephali.'' It turned out to be most
convenient to let a (pre-)bicephalus be a {\it pair} of premice rather than
a premouse with two top extenders.

\begin{defn}\label{defn-pb}
An ordered pair ${\cal N} = ({\cal N}^0,{\cal N}^1)$ is called a generalized
prebicephalus provided the following hold.

(a) ${\cal N}^0 = (J_\alpha[{\vec E}]; \in, {\vec E}, F^0)$
is a premouse
with $F^0 \not=
\emptyset$,

(b) ${\cal N}^1 = (J_\beta[{\vec E'}]; \in, {\vec E'},
F^1)$ is a premouse with
$\beta \geq \alpha$ and $F^1 \not= \emptyset$,

(c) $\alpha$ is a cardinal in ${\cal N}^1$ in case $\beta > \alpha$,

(d) ${\vec E'} {\upharpoonright} \alpha = {\vec E}$, i.e.,
$(J_\alpha[{\vec E'} {\upharpoonright} \alpha]; \in ,
{\vec E'} {\upharpoonright} \alpha) = (J_\alpha[{\vec E}]; \in, {\vec E})$, and

(e) ${\rm c.p.}(F^0) \leq {\rm c.p.}(F^1) < F^0({\rm c.p.}(F^0))$.

In this case ${\cal N}^0$ is called the left part of ${\cal N}$,
${\cal N}^1$ is called the right part of ${\cal N}$, and ${\cal N}$ is called the
generalized prebicephalus derived from ${\cal N}^0$, ${\cal N}^1$.
\end{defn}

Notice that if $\beta > \alpha$ then $\alpha$ has to be a successor cardinal 
in ${\cal N}^1$.

We aim to prove that g-prebicephali trivialize, i.e.,
that ${\cal N}^0 = {\cal N}^1$
provided that the generalized prebicephalus ${\cal N} = ({\cal N}^0,{\cal N}^1)$
meets a certain iterability criterion. Here is an immediate trivial
observation:

\begin{lemma}\label{observation}
Let ${\cal N} = ({\cal N}^0,{\cal N}^1)$ be a generalized prebicephalus, and let ${\cal M}$ be a
premouse. Suppose that ${\cal N}^0 \trianglelefteq {\cal M}$ as well as
${\cal N}^1 \trianglelefteq {\cal M}$. Then ${\cal N}^0 = {\cal N}^1$.
\end{lemma}

{\sc Proof.} The hypothesis tells us that ${\cal N}^0 \trianglelefteq {\cal
N}^1$. In particular $F^0$,
the top extender of ${\cal N}^0$, is on the sequence of ${\cal N}^1$.
Assume that ${\cal N}^0 \triangleleft {\cal
N}^1$, i.e., $\alpha = {\cal N}^0 \cap {\rm OR} < {\cal N}^1 \cap {\rm OR}$. Then $F^0 =
E_\alpha^{{\cal N}^1}$, and so
$\rho_1({\cal N}^0) < \alpha$. But $\alpha$ is supposed to
be a cardinal of ${\cal N}^1$ (cf. (c) of \ref{defn-pb}). Contradiction!
Hence $\alpha = {\cal N}^1 \cap {\rm OR}$, and ${\cal N}^0 = {\cal N}^1$.

\bigskip
\hfill $\square$ (\ref{observation})

\begin{defn}\label{defn-iteration}
Let ${\cal N} = ({\cal N}^0,{\cal N}^1)$ be a generalized prebicephalus. Then $${\cal T}
= (({\cal N}^0_\alpha,\pi_{\alpha \beta}^{{\cal N}^0},{\cal N}^1_\alpha,
\pi_{\alpha \beta}^{{\cal N}^1} \colon \alpha T
\beta <\theta),(E_\alpha \colon \alpha+1 < \theta),T)$$
is an unpadded
iteration tree on ${\cal N}$ provided the following hold.

(a) $T$ is a tree order (in the sense of \cite[Def. 5.0.1]{FSIT}),

(b) ${\cal N}^0_0 = {\cal N}^0$ and ${\cal N}^1_0 = {\cal N}^1$,

(c) for all $\alpha < \theta$,
$({\cal N}^0_\alpha,{\cal N}^1_\alpha)$ is a generalized prebicephalus, or else
${\cal N}^0_\alpha = {\cal N}^1_\alpha$,

(d) for all $\alpha+1 < \theta$, $E_\alpha \not= \emptyset$, and
$E_\alpha = E_\nu^{{\cal N}^0_\alpha}$ for some $\nu \leq {\cal N}^0_\alpha
\cap {\rm OR}$, or
else $E_\alpha$ is the top extender of ${\cal N}^1_\alpha$,

%
(e) for all $\alpha+1 < \theta$, $T$-pred$(\alpha+1) =$ the least $\beta \leq \alpha$
such that ${\rm c.p.}(E_\alpha) < {\rm min} \{ E_\beta({\rm c.p.}(E_\beta)), {\cal N}^0_\beta \cap {\rm OR}
\}$,

(f) for all $\alpha+1 < \theta$ and $h = 0,1$,
if $\beta = T$-pred$(\alpha+1)$ then
$$\pi_{\beta \alpha+1}^{{\cal N}^h} \colon {\cal J}_{\eta^h}^{{\cal N}^h_\beta}
\rightarrow_{E_\alpha} {\cal N}^h_{\alpha+1} {\rm , }$$
where $\eta^h \leq {\cal N}^h_\beta$ is maximal such that $E_\alpha$ measures all the
subsets of its critical point in ${\cal J}_{\eta^h}^{{\cal N}^h_\beta}$
(where we understand that ${\rm deg}^{\cal T}(\alpha+1) = 0$ if ${\cal D}^{\cal T} \cap
[0,\alpha+1] = \emptyset$ and ${\rm deg}^{\cal T}(\alpha+1) =$ that $n$ such that
${\rm dom}(\pi_{\beta \alpha+1}^{{\cal N}^h})$ is $n$-sound if ${\cal D}^{\cal T} \cap
[0,\alpha+1] \not= \emptyset$),

(g) if $\alpha<\theta$ is a limit ordinal then $({\cal N}_\alpha^h,(\pi^{{\cal
N}^h}_{\beta \alpha} \colon \beta T \alpha))$ is the transitive direct limit of
$({\cal N}^h_\beta,\pi^{{\cal
N}^h}_{\gamma \beta} \colon \gamma T \beta < \alpha)$, for $h = 0,1$,

(h) for all $\alpha < \theta$, the set $${\cal D}^{\cal T} \cap (0,\alpha]_T =
\{ \beta+1 T \alpha \colon
{\rm dom}(\pi_{T-pred(\beta+1) \beta+1}^{{\cal N}^h})
\not= {\cal N}^h_{T-pred(\beta+1)} \}$$
is finite, for $h=0,1$, and

(i) for all $\alpha T \beta < \theta$,
$\pi_{\alpha \beta}^{{\cal N}^0} = \pi_{\alpha \beta}^{{\cal N}^1} {\upharpoonright}
{\rm dom}(\pi_{\alpha \beta}^{{\cal N}^0})$.
\end{defn}

By a ``putative'' iteration tree on ${\cal N}$ we shall mean a tree ${\cal T}$
of {\it successor} length $\theta = {\bar \theta} +1$ which is as in
\ref{defn-iteration}
except for the fact that possibly (c) fails for $\alpha = {\bar \theta}$ or (h)
fails for $\beta = {\bar \theta}$. It is crucial, but straightforward,
that if ${\cal T}$ is a putative iteration tree on ${\cal N}$ such that
${\cal N}^0_{\bar \theta}$ and ${\cal N}^1_{\bar \theta}$ are both transitive and (h)
holds for ${\cal T}$,
then ${\cal T}$ is in fact an iteration tree (to get (i) one uses (c) of
\ref{defn-pb}).

For our purposes, we shall also consider padded iteration trees on g-prebicephali;
as usual, this just means that the indexing of the models is slowed down by possible
repetition of models. The reader will have no trouble with
modifying \ref{defn-iteration}
accordingly.

Let $${\cal T}
= (({\cal N}^0_\alpha,\pi_{\alpha \beta}^{{\cal N}^0},{\cal N}^1_\alpha,
\pi_{\alpha \beta}^{{\cal N}^1} \colon \alpha T
\beta <\theta),(E_\alpha \colon \alpha+1 < \theta),T)$$ be a(n) (putative)
iteration tree on the generalized prebicephalus ${\cal N}$. Then we'll write ${\cal M}^{\cal T}_\beta$
for $({\cal N}^0_\beta,{\cal N}^1_\beta)$ if ${\cal D}^{\cal T} \cap (0,\beta]_T =
\emptyset$, and we'll write ${\cal M}^{\cal T}_\beta$ for ${\cal N}^0_\beta =
{\cal N}^1_\beta$ if ${\cal D}^{\cal T} \cap (0,\beta]_T \not=
\emptyset$. We shall also use $({\cal M}^{\cal T}_\beta)^0$ for ${\cal
N}^0_\beta$, and $({\cal M}^{\cal T}_\beta)^1$ for ${\cal
N}^1_\beta$. We'll write $\pi^{\cal T}_{\alpha \beta}$ for $\pi_{\alpha \beta}^{{\cal
N}^1}$, and $E_\beta^{\cal T}$ for $E_\beta$.

\begin{defn}\label{it-of-pb}
Let ${\cal N} = ({\cal N}^0,{\cal N}^1)$ be a generalized prebicephalus,
and let ${\cal M}$ be a premouse. We want to
determine a pair $(\partial_0^{{\cal N},{\cal M}},\partial_1^{{\cal N},{\cal M}})$,
the
``least disagreement'' of ${\cal N}$ and
${\cal M}$.

{\em Case 1.} ${\cal N}^0$ and ${\cal M}$ are not lined up. Let $\nu$ be least such
that $$E^{{\cal N}^0}_\nu \not= E^{\cal M}_\nu {\rm , }$$ and let
$$\partial_0^{{\cal N},{\cal M}}
= E^{{\cal N}^0}_\nu {\rm \ and \ } \partial_1^{{\cal N},{\cal M}} = E^{\cal M}_\nu.$$

{\em Case 2.} Not Case 1, but
${\cal N}^1$ and ${\cal M}$ are not lined up. Then let $\partial_0^{{\cal N},{\cal M}}$
be the top extender of ${\cal N}^1$, and let
$\partial_1^{{\cal N},{\cal M}}$ be the top extender
of ${\cal N}^0$ (sic!).

{\em Case 3.} Neither Case 1 nor Case 2. Then we let
$(\partial_0^{{\cal N},{\cal M}},\partial_1^{{\cal N},{\cal M}})
= (\emptyset, \emptyset)$.
\end{defn}

Case 2 needs a brief discussion. We'll have that
${\cal N}^0 \trianglelefteq {\cal M}$, as
${\cal N}^0$ and ${\cal M}$ are lined up but ${\cal N}^1$ and ${\cal M}$
are not. In particular, the top extender of ${\cal N}^0$ appears on the sequence of
${\cal M}$. But the top extender of ${\cal N}^1$ can't appear on the sequence of
${\cal M}$, as ${\cal N}^1$ and ${\cal M}$
are not lined up.

The content of \ref{observation} may now be stated as saying that if ${\cal M}$ is a
premouse with $(\partial_0^{{\cal N},{\cal M}},\partial_1^{{\cal N},{\cal M}}) =
(\emptyset,\emptyset)$ then ${\cal M} \triangleleft {\cal N}^0$ or else ${\cal N}^0 =
{\cal N}^1$.

\begin{defn}\label{defn-bicephalus}
We call a generalized prebicephalus ${\cal N}$
a generalized bicephalus if it is coiterable with $K^c$ in the
following sense.

Let ${\cal T}$, ${\cal U}$ be putative padded
iteration trees of successor length $\theta+1$
on ${\cal N}$, $K^c$ given by the following characterization. If $\beta < \theta$ then
$$(E_\beta^{\cal T},E_\beta^{\cal U}) = (\partial_0^{{\cal M}^{\cal
T}_\beta,{\cal M}^{\cal U}_\beta},\partial_1^{{\cal M}^{\cal
T}_\beta,{\cal M}^{\cal U}_\beta})$$ in case ${\cal M}^{\cal
T}_\beta$ is a generalized prebicephalus, and let otherwise $$(E_\beta^{\cal T},E_\beta^{\cal U})
= (E_\nu^{{\cal M}^{\cal U}_\beta},E_\nu^{{\cal M}^{\cal T}_\beta})$$ where $\nu$ is
least with $E_\nu^{{\cal M}^{\cal U}_\beta} \not= E_\nu^{{\cal M}^{\cal T}_\beta}$.
{\em Then} ${\cal T}$ is an iteration tree on ${\cal N}$.
\end{defn}

We want to emphasize that by \ref{it-of-pb}
we have that
the ${\cal N}$-side of the coiteration of
the generalized prebicephalus ${\cal N}$ with $K^c$, being performed as in
\ref{defn-bicephalus},
will only use extenders which are allowed by
\ref{defn-iteration}
(d).

\begin{lemma}\label{semi-normality}
Let ${\cal N}$, ${\cal T}$, ${\cal U}$ and $\theta$ be as in
\ref{defn-bicephalus}. Then ${\cal U}$ is normal.
Moreover, for all $\beta+1 < \theta$ do we have that
the index of $E_\beta^{\cal T}$ is $\geq$ the index of
$E_\beta^{\cal U}$.
\end{lemma}

{\sc Proof.} The only thing to notice here is that if
$(E_\beta^{\cal T},E_\beta^{\cal U})$ is chosen according to Case 2 in
\ref{it-of-pb} and $\nu$ denotes the index of $E_\beta^{\cal U}$ then
${\cal M}^{\cal U}_{\beta+1}$, $({\cal M}^{\cal T}_{\beta+1})^0$, and
$({\cal M}^{\cal T}_{\beta+1})^1$ pairwise agree up to 
$\nu+1$. This holds because
then $\nu = ({\cal M}^{\cal T}_{\beta})^0 \cap {\rm OR}$, so that either $\nu$ is the
index of $E_\beta^{\cal T}$ or else $\nu$ is a (successor) cardinal in
$({\cal M}^{\cal T}_{\beta})^1$. In any event, $$E_\nu^{({\cal M}^{\cal
T}_{\beta+1})^0} = E_\nu^{({\cal M}^{\cal
T}_{\beta+1})^1} = \emptyset = E^{{\cal M}^{\cal U}_{\beta+1}}_\nu.$$

\bigskip
\hfill $\square$ (\ref{semi-normality})

%
%
\begin{lemma}\label{bicephali-trivialize} $( \ \lnot \ \zerohandgrenade \ )$
Let ${\cal N} = ({\cal N}^0,{\cal N}^1)$ be a generalized bicephalus.
Then ${\cal N}^0 = {\cal N}^1$.
\end{lemma}

{\sc Proof.} Let
${\cal
T}$, ${\cal U}$ denote the (padded) iteration trees
arising from the coiteration of ${\cal N}$ with $K^c$
built as in
\ref{defn-bicephalus}, where
$${\cal T}
= (({\cal N}^0_\alpha,\pi_{\alpha \beta}^{{\cal N}^0},{\cal N}^1_\alpha,
\pi_{\alpha \beta}^{{\cal N}^1} \colon \alpha\leq
\beta <\theta),(E_\alpha \colon \alpha+1 < \theta),T).$$

\bigskip
{\bf Claim.}
${\rm lh}({\cal T}) = {\rm lh}({\cal U}) <
{\rm OR}$.

\bigskip
Given this Claim, the proof of \ref{bicephali-trivialize} can be completed as follows.
By the proof of \ref{W-is-universal}
we shall have that ${\cal M}_\infty^{\cal T}$ is a generalized prebicephalus,
and $({\cal M}_\infty^{\cal T})^0 \trianglelefteq {\cal M}^{\cal U}_\infty$ as
well as $({\cal M}_\infty^{\cal T})^1 \trianglelefteq {\cal M}^{\cal
U}_\infty$. This implies that
$({\cal M}_\infty^{\cal T})^0 = ({\cal M}_\infty^{\cal T})^1$ by
\ref{observation}. Now suppose that ${\cal N}^0 \not= {\cal N}^1$, and let $w \in
{\cal N}^0$ be such that $$w \in F^0
\Leftrightarrow w \notin F^1 {\rm , }$$ where $F^0$ and $F^1$ are the top extenders of
${\cal N}^0$ and ${\cal N}^1$. (Notice that such $w$ would have to exist!)
We have that ${\rm dom}(\pi^{{\cal N}^0}_{0 \infty}) = {\cal N}^0$,
${\rm dom}(\pi^{{\cal N}^1}_{0 \infty}) = {\cal N}^1$, and $\pi^{{\cal N}^0}_{0 \infty}
= \pi^{{\cal N}^1}_{0 \infty} {\upharpoonright} {\cal N}^0$. If we let ${\tilde F}^0$
and ${\tilde F}^1$ denote the top extenders of ${\cal N}^0_\infty$ and
${\cal N}^1_\infty$ then ${\tilde F}^0 =
{\tilde F}^1$, and thus
$$w \in F^0 \Leftrightarrow \pi^{{\cal N}^0}_{0 \infty}(w) \in {\tilde F}^0
\Leftrightarrow \pi^{{\cal N}^1}_{0 \infty}(w) \in {\tilde F}^1
\Leftrightarrow w \in F^1.$$ Contradiction!

\bigskip
{\sc Proof} of the Claim.
Let us assume that ${\rm lh}({\cal T}) = {\rm lh}({\cal U}) = {\rm OR}+1$. The proof of
\ref{W-is-universal} yields that then
$\pi_{0 \infty}^{\cal U} {\rm " } {\rm OR}
\not \subset
{\rm OR}$, so that usual arguments give some $\lambda \in [0,{\rm OR})_T \cap [0,{\rm OR})_U$ so that

\begin{equation}\label{C1}
{\rm c.p.}(\pi_{\lambda \infty}^{\cal T}) =
{\rm c.p.}(\pi_{\lambda \infty}^{\cal U}) = \lambda {\rm \ \ and }
\end{equation}
\begin{equation}\label{C2}
\pi_{\lambda \infty}^{\cal T} {\upharpoonright} {\cal P}(\lambda) \cap {\cal
M}^{\cal T}_{{\lambda}} =
\pi_{\lambda \infty}^{\cal U} {\upharpoonright} {\cal P}(\lambda) \cap {\cal
M}^{\cal U}_{\lambda}.
\end{equation}

\noindent Let $\alpha+1$ be
least in $(\lambda,{\rm OR}]_T$,
and let $\beta+1$ be least in $(\lambda,{\rm OR}]_U$.
Let $E_\beta^{\cal U} = E_{\nu_1}^{{\cal M}_\beta^{\cal U}}$, and
let $E_\alpha^{\cal T} = E^{({\cal M}^{\cal T}_\alpha)^0}_{\nu_0}$ or
$E_\alpha^{\cal T} = E^{({\cal M}^{\cal T}_\alpha)^1}_{\nu_0}$ 
(if the former is not true but the
latter, then ${\cal M}^{\cal T}_\alpha$ is a generalized prebicephalus and $E_\alpha^{\cal T}$
has to be the top extender of $({\cal M}^{\cal T}_\alpha)^1$).
Of course, $\lambda = {\rm c.p.}(E^{\cal U}_\beta) = {\rm c.p.}(E^{\cal T}_\alpha)$.

\bigskip
{\it Case 1.} $E_\alpha^{\cal T} = E^{({\cal M}^{\cal T}_\alpha)^0}_{\nu_0}$.

\bigskip
In this case, the rest of the coiteration is beyond $\nu_0+1$ on both sides, i.e., all
$E_\gamma^{\cal T}$ and $E_\gamma^{\cal U}$ for $\gamma > \alpha$ have index $> \nu_0$.
Moreover, by \ref{defn-iteration} (e) we'll have that $\pi_{\alpha+1 \infty}^{\cal T}
{\upharpoonright} E_\alpha^{\cal T}(\lambda) = {\rm id}$. For this reason, and because
${\cal U}$ is normal by \ref{semi-normality}, (\ref{C1}) and
(\ref{C2}) give that $E_\alpha^{\cal T}$ and $E_\beta^{\cal U}$ are compatible.

\bigskip
{\it Case 1.1.} $\nu_0 < \nu_1$.

\bigskip
Then $E_\alpha^{\cal T} \in {\cal J}^{{\cal M}_\beta^{\cal U}}_{\nu_1}$ by the initial
segment condition for ${\cal M}_\beta^{\cal U}$. But
${\cal U}$ is normal by \ref{semi-normality}, so that
$E_\alpha^{\cal T} \in {\cal M}_\infty^{\cal U}$, i.e., $E_\alpha^{\cal T} \in
{\cal J}^{{\cal M}_\infty^{\cal U}}_{{\rm OR}}$.
On the other hand $E_\alpha^{\cal T} \notin ({\cal M}^{\cal T}_{\alpha+1})^0$,
which by $\pi_{\alpha+1 \infty}^{\cal T}
{\upharpoonright} E_\alpha^{\cal T}(\lambda) = {\rm id}$
implies that $E_\alpha^{\cal T} \notin ({\cal M}^{\cal T}_{\infty})^0$.
But of course we have ${\cal J}^{{\cal M}_\infty^{\cal U}}_{{\rm OR}} =
{\cal J}^{({\cal M}^{\cal T}_{\infty})^0}_{{\rm OR}}$. Contradiction!

\bigskip
{\it Case 1.2.} $\nu_1 < \nu_0$.

\bigskip
We then have $E_\beta^{\cal U} \in {\cal J}^{({\cal M}_\alpha^{\cal T})^0}_{\nu_0}$
by the initial
segment condition for $({\cal M}_\alpha^{\cal T})^0$. So by
$\pi_{\alpha+1 \infty}^{\cal T}
{\upharpoonright} E_\alpha^{\cal T}(\lambda) = {\rm id}$ we know that
$E_\beta^{\cal U} \in ({\cal M}^{\cal T}_\infty)^0$. On the other hand,
$E_\beta^{\cal U} \notin {\cal M}_{\beta+1}^{\cal U}$, which by the normality of
${\cal U}$ implies that $E_\beta^{\cal U} \notin {\cal M}_{\infty}^{\cal U}$.
We thus have a contradiction as in case 1.1!

\bigskip
{\it Case 1.3.} $\nu_1 = \nu_0$.

\bigskip
Then we
get $\alpha \leq \beta$ by (the proof of)
\ref{semi-normality}; but then as $E_\gamma^{\cal U}$ has
index $> \nu_0$ for $\gamma > \alpha$ we must have
$\alpha = \beta$.
Now by the normality of ${\cal U}$ and by
$\pi_{\alpha+1 \infty}^{\cal T}
{\upharpoonright} E_\alpha^{\cal T}(\lambda) = {\rm id}$ we must have
$E_\alpha^{\cal T} = E_\alpha^{\cal U}$. This clearly contradicts the choice of
$(E_\alpha^{\cal T},E_\alpha^{\cal U})$ as $(\partial_0^{{\cal M}^{\cal
T}_\alpha,{\cal M}^{\cal U}_\alpha},\partial_1^{{\cal M}^{\cal
T}_\alpha,{\cal M}^{\cal U}_\alpha})$.

\bigskip
{\it Case 2.} ${\cal M}^{\cal T}_\alpha$ is a generalized prebicephalus and $E_\alpha^{\cal T}$
is the top extender of $({\cal M}^{\cal T}_\alpha)^1$.

\bigskip
In this case, by \ref{it-of-pb}, $E_\alpha^{\cal U}$ is the top extender of
$({\cal M}^{\cal T}_\alpha)^0$.
Let $\nu^*$ be the index of $E_\alpha^{\cal U}$. So
$\nu^* = ({\cal M}_\alpha^{\cal T})^0 \cap {\rm OR}$.

If $\alpha < \beta$ then by \ref{defn-pb} (e), as ${\cal
M}_\alpha^{\cal T}$ is a generalized prebicephalus, and by \ref{semi-normality}
we would get that
$${\rm c.p.}(E_\alpha^{\cal U}) \leq {\rm c.p.}(E_\alpha^{\cal T}) = \lambda = {\rm c.p.}(E^{\cal
U}_\beta) < \nu^* < \nu_1.$$
However $K^c$ is below $\zerohandgrenade$, and ${\cal
T}$ is a normal iteration tree on $K^c$. We thus have a contradiction with Claim 1 in
the proof of \ref{normal-linear}.

Thus $\alpha \geq \beta$. If $\alpha > \beta$ then, because
${\rm c.p.}(E_\alpha^{\cal U}) \leq {\rm c.p.}(E_\alpha^{\cal T}) = \lambda = {\rm c.p.}(E^{\cal
U}_\beta)$, we would get $\beta+1 \notin (0,\infty)_U$, which is nonsense.

We thus have that $\alpha = \beta$. In particular,
$\nu_1 = \nu^* = ({\cal M}_\alpha^{\cal T})^0 \cap {\rm OR}$. Hence by
\ref{defn-iteration} (e) we'll have that

\begin{equation}\label{C3}
\pi_{\alpha+1 \infty}^{\cal T}
{\upharpoonright} {\rm min} \{ E_\alpha^{\cal T}(\lambda) , \nu_1 \} = {\rm id}.
\end{equation}

\noindent Because ${\cal U}$ is normal by \ref{semi-normality} we'll have that

\begin{equation}\label{C4}
\pi_{\alpha+1 \infty}^{\cal U}
{\upharpoonright} E^{\cal U}_\alpha(\lambda) = {\rm id}.
\end{equation}

\noindent
Recall that $E_\alpha^{\cal U}$
is the top extender of $({\cal M}_\alpha^{\cal T})^0$,
and $E_\alpha^{\cal T}$ is the top extender of $({\cal M}_\alpha^{\cal T})^1$.
In particular,

\begin{equation}\label{C5}
E_\alpha^{\cal U}(\lambda) \leq
{\rm min} \{ E_\alpha^{\cal T}(\lambda) , \nu_1 \} {\rm . }
\end{equation}

\noindent But now (\ref{C1}) to
(\ref{C5}) will tell us that

\begin{equation}\label{C6}
E_\alpha^{\cal U} = E_\alpha^{\cal T} |
E_\alpha^{\cal U}(\lambda).
\end{equation}

\noindent The initial segment condition for $({\cal M}_\alpha^{\cal T})^1$ then yields
$E_\alpha^{\cal U} = E_\alpha^{\cal T}$, or else $E_\alpha^{\cal U}
\in ({\cal M}_\alpha^{\cal T})^1$.

\bigskip
{\it Case 2.1.} $E_\alpha^{\cal U}
\in ({\cal M}_\alpha^{\cal T})^1$.

\bigskip
In this case, we'll have that $E_\alpha^{\cal U} \in ({\cal M}_{\alpha+1}^{\cal T})^1$,
too. By (\ref{C3}) and (\ref{C5}) we have that $\pi_{\alpha+1 \infty}^{\cal T}
{\upharpoonright} E_\alpha^{\cal U}(\lambda) = {\rm id}$, and thus $(X,Y) \in E_\alpha^{\cal U}$
if and only if
$$\exists {\tilde Y} \ [ (X,{\tilde Y}) \in
\pi_{\alpha+1 \infty}^{\cal T}(E_\alpha^{\cal U}) \wedge Y = {\tilde Y} \cap
E_\alpha^{\cal U}(\lambda) ].$$
Hence $E_\alpha^{\cal U} \in ({\cal M}_\infty^{\cal T})^1$.

On the other hand, the normality of ${\cal U}$ yields that $E_\alpha^{\cal
U} \notin {\cal M}_\infty^{\cal U}$. However, ${\cal J}_{{\rm OR}}^{{\cal
M}_\infty^{\cal U}} = {\cal J}_{{\rm OR}}^{({\cal M}_\infty^{\cal T})^1}$. This is a
contradiction!

\bigskip
{\it Case 2.2.}
$E_\alpha^{\cal U} = E_\alpha^{\cal T}$.

\bigskip
But
$E_\alpha^{\cal U} = E_\alpha^{\cal T}$ of course
contradicts the choice of
$(E_\alpha^{\cal T},E_\alpha^{\cal U})$ as $(\partial_0^{{\cal M}^{\cal
T}_\alpha,{\cal M}^{\cal U}_\alpha},\partial_1^{{\cal M}^{\cal
T}_\alpha,{\cal M}^{\cal U}_\alpha})$.

\bigskip
\hfill $\square$ (Claim)

\hfill $\square$ (\ref{bicephali-trivialize})

\bigskip
We want to point out that \ref{bicephali-trivialize} 
would no longer be true if we
dropped (e) in the definition \ref{defn-pb} (there are easy counterexamples).

%
%
\begin{lemma}\label{it-of-bicephalus} $( \ \lnot \ \zerohandgrenade \ )$
Let ${\cal N} = ({\cal N}^0,{\cal N}^1)$ be a generalized prebicephalus,
and let $\mu$
be the critical point of
${\cal N}^1$. Let $\kappa > \mu$ be a
cardinal in ${\cal N}^0$ (and hence in ${\cal N}^1$, too).
Suppose that for $h \in \{ 0 , 1 \}$ do we have the following.

For all iteration trees ${\cal V}$
on ${\cal N}^h$ with last model ${\cal
M}^{\cal V}_\infty$ and for all ${\cal
N}^h$-cardinals $\rho \leq \kappa$ such that either

(a) $\rho < \kappa$, and ${\cal V}$ lives on ${\cal J}^{{\cal N}^h}_\kappa$ and
is above $\rho$, but ${\cal V}$ doesn't use extenders which are
total on ${\cal N}^h$, or else

(b) $\rho = \kappa$, and ${\cal V}$ is above $\rho$,

\noindent we have that

($\star$) if $F =
E^{{\cal
M}^{\cal V}_\infty}_\nu \not= \emptyset$ is
such that $\nu > \rho$ and ${\rm c.p.}(F) < \rho$, then $F$ is countably complete.

Then ${\cal N}$ is a generalized bicephalus.
\end{lemma}

Notice that if ${\cal N}$ satisfies the assumption in the statement of
\ref{it-of-bicephalus} then in particular every $E_\nu^{{\cal N}^h} \not= \emptyset$
with ${\rm c.p.}(E_\nu^{{\cal N}^h}) < \kappa$ and which is total on ${{\cal N}^h}$ is
countably complete. (Let ${\cal V} =$ the trivial tree, and $\rho =
{\rm c.p.}(E_\nu^{{\cal N}^h})^{+{{\cal N}^h}}$.)

\bigskip
{\sc Proof} of \ref{it-of-bicephalus}. We have
to show that ${\cal N}$ is iterable. Let
${\cal T}$ be a
putative
iteration tree on ${\cal N}$
arising from the comparison with $K^c$. By the proof of
\ref{bicephali-trivialize}, ${\cal T}$ is a set. Pick an elementary embedding
$\sigma \colon N \rightarrow H_\theta$ (for some large enough regular
$\theta$)
with $N$ countable and transitive, and such that $\{ {\cal T} , \kappa \} \subset
{\rm ran}(\sigma)$. Let
${\bar {\cal T}} = \sigma^{-1}({\cal T})$, and ${\bar \kappa} = \sigma^{-1}(\kappa)$.
It suffices to embed the models of ${\bar {\cal T}}$ into transitive structures.
For notational convenience, we shall assume that ${\bar {\cal T}}$ is 
unpadded.
Say,
$${\bar {\cal T}} =
((({\cal M}_\alpha^{\bar {\cal T}})^0,
\pi_{\alpha \beta}^0,
({\cal M}_\alpha^{\bar {\cal T}})^1,
\pi_{\alpha \beta}^1 \colon \alpha {\bar T} \beta < {\rm lh}({\bar T})),
(E_\alpha \colon \alpha+1 < {\rm lh}({\bar T})), {\bar T}).$$ Put
$\kappa_\alpha = {\rm c.p.}(E_\alpha)$ for $\alpha+1 < {\rm lh}({\bar T})$.
We set $n(0) = 0$, and
for $\alpha>0$ we let $n(\alpha)$ be such that
$$\exists \beta \ ( \ \beta {\bar T} \alpha \wedge \forall \gamma+1 \in
(\beta,\alpha]_{\bar T} \ \exists \eta \ {\cal M}_{\gamma+1}^{\cal T} =
{\rm Ult}_{n(\alpha)}({\cal J}^{{\cal M}^{\cal T}_{T-pred(\gamma+1)}}_{\eta};
E^{\cal T}_\gamma) \ ) \ ).$$
Notice that for all $\alpha<{\rm lh}({\cal T})$ and $h \in \{ 0,1 \}$
do we have that
$({\cal M}^{\bar {\cal T}}_\alpha)^h$ is $n(\alpha)$-sound.

We shall in fact determine
some $\vartheta \leq {\rm lh}({\bar {\cal T}})$, and define, for both $h \in \{ 0 , 1 \}$,
sequences ${\vec {\cal V}}^h = ({\cal V}_i^h
\colon i<\vartheta)$
of iteration trees
on ${\cal N}^h$ such that any of the ${\cal V}_i^h$'s is as in (a) or (b) of
\ref{it-of-bicephalus}; for this purpose we'll simultaneously
define a sequence $(\rho_i \colon i < \vartheta)$ of ordinals.
We shall also
construct maps $\sigma_\alpha^h$
from any $({\cal M}_\alpha^{\bar {\cal T}})^h$ into a model
of some ${\cal V}_i^h$.
For all $i < \vartheta$ we'll have that ${\cal V}_i^0$ and ${\cal V}_i^1$
are given by the very same sequence of extenders; in particular, ${\rm lh}({\cal V}_i^0)
= {\rm lh}({\cal V}_i^1)$.
We shall index the models of ${\cal V}_i^h$
in a non-standard way,
namely, we shall start counting them
with $$\ell(i) = \sum_{j<i} \ {\rm lh}({\cal V}_j^h)$$ (where $\sum$
denotes ordinal summation), so that the models of ${\cal V}_i^h$ will be indexed by
the elements of $[\ell(i),\ell(i)+{\rm lh}({\cal V}_i^h))$, and
we'll conveniently have that
$$\sigma_\alpha^h \colon ({\cal M}_\alpha^{\bar {\cal T}})^h
\rightarrow {\cal M}_\alpha^{{\cal V}_i^h}$$ for that $i < \vartheta$
such that $\ell(i)
\leq \alpha < \ell(i)+{\rm lh}({\cal V}_i^h)$.
It will be clear from the construction that always

\bigskip
{\bf R 0}${}_\alpha$
$\sigma_\alpha^0 = \sigma_\alpha^1
{\upharpoonright} ({\cal M}_\alpha^{\bar {\cal T}})^0$.

\bigskip
Letting $h$ range over $\{ 0 , 1 \}$, we shall inductively
maintain that
the following requirements are met as well.

\bigskip
{\bf R 1}${}_\alpha$ $\ \ \forall \ \beta \in [0,\alpha]_{\bar T}$, if
$i$ is maximal with $\ell(i) \leq \beta$ then
$$\sigma_\beta^h \colon ({\cal M}_\beta^{{\bar {\cal T}}})^h \rightarrow
{\cal M}_\beta^{{\cal V}_i^h}$$ is a weak $n(\beta)$-embedding.

\bigskip
{\bf R 2}${}_\alpha$ $\ \ $ if $i$ is maximal with $\ell(i) \leq \alpha$ then
$\ell(i) {\bar T} \alpha$, and for all $\beta \in [0,\ell(i)]_{\bar T}$ do we
have that $\exists j \leq i \ ( \ \beta = \ell(j) \
\wedge \ {\cal M}_\beta^{{\cal V}_j^h}
= {\cal N}^h \ )$.

\bigskip
{\bf R 3}${}_\alpha$
for all $\beta+1 \in (0,\alpha]_{\bar T}$, setting $\beta^* =
{\bar T}$-pred$(\beta+1)$, we have that
$\sigma_\alpha^h$ agrees with $\sigma_{\beta^*}^h$ up to $\kappa_\beta$; i.e.,
$${\cal J}^{{\cal M}_\alpha^{{\cal V}_i^h}}_{\sigma_\alpha^h(\kappa_\beta)} =
{\cal J}^{{\cal M}_{\beta^*}^{{\cal V}_j^h}}_{\sigma_\alpha^h(\kappa_\beta)}
{\rm , \ and
}$$
$$\sigma_\alpha^h {\upharpoonright} {\cal J}_{\kappa_\beta}^{({\cal M}_\alpha^{\bar
{\cal T}})^h} = \sigma_\beta^h {\upharpoonright}
{\cal J}_{\kappa_{\beta}}^{({\cal M}_{\beta^*}^{\bar
{\cal T}})^h} {\rm , }$$
where $i$ is maximal with $\ell(i) \leq \alpha$ and
$j$ is maximal with $\ell(j) \leq \beta^*$.

\bigskip
{\bf R 4}${}_\alpha$ $\ \ $ if
$i$ is maximal with $\ell(i) \leq \alpha$,
$\beta {\bar T} \gamma \in [0,\alpha)_{\bar T}$
and ${\cal D}^{\bar {\cal T}} \cap (\beta,\gamma]_{\bar T} = \emptyset$ then we have:

(a) $\ell(i) {\bar T} \beta {\bar T} \gamma \Rightarrow
\sigma_\gamma^h \circ \pi^h_{\beta \gamma} = \pi_{\beta \gamma}^{{\cal
V}^h_i} \circ \sigma_\beta^h$, and

(b) $\beta {\bar T} \gamma {\bar T} \ell(i) \Rightarrow
\sigma_\gamma^h \circ \pi^h_{\beta \gamma} = \sigma_\beta^h$.

\bigskip
{\bf R 5}${}_\alpha$ $\ \ $ if
$i$ is maximal with $\ell(i) \leq \alpha$,
then for all $j \leq i$ do we have that ${\cal V}' = {\cal V}_j^h
{\upharpoonright} {\rm min} \{ \ell(j)+{\rm lh}({\cal V}^h_j),\alpha+1 \}$
is an iteration tree on ${\cal N}^h$ above $\rho_j$ as in
(a) or (b) of \ref{it-of-bicephalus}; i.e., either

$\rho_j < \kappa$, and ${\cal V}'$ lives on ${\cal J}^{{\cal N}^h}_\kappa$ and
is above $\rho_j$, but ${\cal V}'$ doesn't use extenders which are
total on ${\cal N}^h$, or else

$\rho_j = \kappa$, and ${\cal V}'$ is above $\rho_j$.

\bigskip
We are now going to run our construction.
In what follows we again let $h$ range over $\{ 0,1 \}$. Put $\ell(0) = 0$.
To commence, set $\sigma_0^h = \sigma {\upharpoonright} ({\cal M}_0^{\bar
{\cal T}})^h$ (notice $({\cal M}_0^{\bar
{\cal T}})^h = \sigma^{-1}({\cal N}^h)$),
and let ${\cal V}_0^h {\upharpoonright} 1$ be trivial.
It is clear that {\bf R 1}${}_0$ to {\bf R 5}${}_0$ hold.

Now suppose that we have defined everything up to $\alpha$,
that is, suppose that for some $i$ are we given
$${\vec {\cal V}}^h {\upharpoonright} i =
({\cal V}^h_j \colon j < i) {\rm , \ } {\cal V}^h_i {\upharpoonright} \alpha+1
{\rm , \ } (\sigma^h_\beta \colon \beta \leq \alpha) {\rm , }
(\rho_j \colon j < i) {\rm , \ and \ } \rho_i {\rm \ if \ } \alpha >
\ell(i)$$
in such a way that {\bf R 1}${}_\alpha$ through {\bf R 5}${}_\alpha$ are satisfied.
Suppose that $\alpha < {\rm lh}({\bar {\cal T}})$, and
set ${\bar F} = E^{\bar {\cal T}}_\alpha$, and ${\bar \mu} = {\rm c.p.}({\bar F})$.
(If $\alpha = {\rm lh}({\bar {\cal T}})$
we're done with our construction.) Then either ${\bar F} =
E_{\nu}^{({\cal M}_\alpha^{\bar {\cal T}})^0}$ for some $\nu \leq
({\cal M}_\alpha^{\bar {\cal T}})^0 \cap {\rm OR}$, or else ${\bar F}$ is the top extender
of $({\cal M}_\alpha^{\bar {\cal T}})^1$.
Depending on whether the former holds, or the latter, we set
$F = \sigma_\alpha^0({\bar F})$ and $h^* = 0$, or
$F = \sigma_\alpha^1({\bar F})$  and $h^* = 1$.
Set $\mu^* = {\rm c.p.}(F)$.

\bigskip
{\it Case 1.} $\alpha = \ell(i) \wedge \mu^* < \kappa \ \wedge$
${\bar F}$ is total on $({\cal M}_\alpha^{{\bar {\cal T}}})^{h^*}$,
or $\alpha > \ell(i) \wedge
\mu^* < \rho_i$.

\bigskip
By {\bf R 5}${}_\alpha$ and our assumptions on ${\cal N}$, we get that
$F$ is countably complete in this case. Let $\tau \colon
\sigma_\alpha^{h^*} {\rm " }
{\bar F}({\bar \mu}) \rightarrow \mu^*$ be order
preserving such that for appropriate $a$, $X \in {\rm ran}(\sigma_\alpha^{h^*})$
we have that $a \in F(X) \Rightarrow \tau {\rm " } a \in X$.
We now declare ${\cal
V}_i^h = {\cal
V}_i^h {\upharpoonright} \alpha+1$, i.e., $\ell(i+1) = \alpha+1$.
If $\alpha = \ell(i)$ we also put $\rho_i = 0$ (we won't be interested in this
value).

By \ref{defn-iteration} and {\bf R 5}${}_\alpha$
it is easy to see that
${\bar T}$-pred$(\alpha+1) \leq \ell(i)$, so that by
{\bf R 2}${}_\alpha$ we know that $\exists j \leq i \ \ell(j) =
{\bar T}$-pred$(\alpha+1)$. By {\bf R 3}${}_\alpha$, $\sigma^h_\alpha$ agrees with
$\sigma^h_{\ell(j)}$ up to $\kappa_\beta$, where $\beta+1$ is least in
$(\ell(j),\alpha]_{\bar T}$. We have $\kappa_\alpha < \kappa_\beta$ by $\lnot \
\zerohandgrenade$
(cf. Claim 2 in the proof of \ref{normal-linear}),
so that $\sigma^h_\alpha$ agrees with
$\sigma^h_{\ell(j)}$ up to $\kappa_\alpha^+$ (calculated in $({\cal M}_\alpha^{\bar
{\cal T}})^h$). By {\bf R 0}${}_\alpha$, hence, $\sigma_\alpha^{h^*}$ agrees with
$\sigma^h_{\ell(j)}$ up to $\kappa_\alpha^+$.

We now define $\sigma_{\alpha+1}^h \colon ({\cal M}^{\bar {\cal
T}}_{\alpha+1})^h \rightarrow {\cal N}^h$ by setting
$$[a,f] \mapsto \sigma_{\ell(j)}^h(f)(\tau {\rm " } \sigma_\alpha^{h^*}(a)).$$
This is well-defined and $\Sigma_0$-elementary, as we may reason as follows.
Let $\Phi$ be a $\Sigma_0$ formula. Then
$$({\cal M}_{\alpha+1}^{{\bar {\cal T}}})^h \models 
\Phi([a_1,f_1],...,[a_k,f_k]) {\rm \ \ }
\Leftrightarrow$$
$$(a_1,...,a_k) \in {\bar F}( \{ (u_1,...,u_k) \ \colon \ 
({\cal M}_{\ell(j)}^{{\bar {\cal T}}})^h
\models \Phi(f_1(u_1),...,f_k(u_k)) \} ) {\rm \ \ } \Leftrightarrow$$

\bigskip
$(\sigma_\alpha^{h^*}(a_1),...,\sigma_\alpha^{h^*}(a_k)) \in$
$$F(\sigma_\alpha^{h^*}( \{ (u_1,...,u_k) \
\colon \ ({\cal M}_{\ell(j)}^{{\bar {\cal T}}})^h
\models \Phi(f_1(u_1),...,f_k(u_k)) \} ) ) {\rm , }$$
which, by the amount of agreement of $\sigma_\alpha^{h^*}$ with
$\sigma_{\ell(j)}^{h}$ and by {\bf R 1}${}_\alpha$,
holds if and only if

\bigskip
$(\sigma_\alpha^{h^*}(a_1),...\sigma_\alpha^{h^*}(a_k)) \in$
$$F( \{ (u_1,...,u_k) \ \colon \ {\cal N}^h \models 
\Phi(\sigma^h_{\ell(j)}(f_1)(u_1),...,\sigma^h_{\ell(j)}(f_k)(u_k)) \} )
{\rm \ \ } \Leftrightarrow$$
$${\cal N}^h \models
\Phi(\sigma^h_{\ell(j)}(f_1)(\tau 
{\rm " })\sigma_\alpha^{h^*}(a_1),...,\sigma^h_{\ell(j)}(f_k)(\tau 
{\rm " })\sigma_\alpha^{h^*}(a_k)).$$

It is straightforward to check that now {\bf R 1}${}_{\alpha+1}$ through
{\bf R 5}${}_{\alpha+1}$ hold.

\bigskip
{\it Case 2.} $\alpha = \ell(i) \wedge \mu^* \geq \kappa$,
or $\alpha = \ell(i) \wedge \mu^* < \kappa \wedge {\bar F}$ is partial on $({\cal
M}_\alpha^{\bar {\cal T}})^{h^*}$,
or $\alpha > \ell(i) \wedge
\mu^* \geq \rho_i$.

\bigskip
In this case, we start or continue copying ${\bar {\cal T}} {\upharpoonright}
[\ell(i),\alpha+2)$ onto ${\cal N}^h$,
getting ${\cal V}^h_i {\upharpoonright} \alpha+2$. We let
${\cal V}^h_i$-pred$(\alpha+1) = {\bar {\cal T}}$-pred$(\alpha+1)$, call it
$\alpha^*$.
We'll have $\ell(i) \leq \alpha^*$.
We then use
the shift lemma \cite[Lemma 5.2]{FSIT}
to get $\pi_{\alpha^* \alpha+1}^{{\cal V}^h_i}$ together with
the copy map $\sigma^h_{\alpha+1}$. It is easy to check that
{\bf R 1}${}_{\alpha+1}$ to
{\bf R 5}${}_{\alpha+1}$ hold.

We put $\rho_i = \kappa$ if $\alpha = \ell(i) \wedge \mu^* \geq \kappa$,
and we put $\rho_i =$ the cardinality of ${\rm c.p.}(F)$ in ${\cal N}^{h^*}$ if
$\alpha = \ell(i) \wedge \mu^* < \kappa \wedge {\bar F}$ is partial on $({\cal
M}_\alpha^{\bar {\cal T}})^{h^*}$.

\bigskip
Now let $\lambda$ be a limit ordinal and suppose that we 
have defined everything
up to $\lambda$. To state this more precisely, we have to consider two cases.

\bigskip
{\it Case 1.} There are cofinally in $\lambda$ many $\alpha < \lambda$ such that there is a $j$ with $\alpha = \ell(j)$.

\bigskip
In this case for some $i$ are we given
$${\vec {\cal V}}^h {\upharpoonright} i =
({\cal V}^h_j \colon j < i)
{\rm , \ } (\sigma^h_\beta \colon \beta < \lambda) {\rm , \ and \ }
(\rho_j \colon j < i)$$
in such a way that {\bf R 1}${}_\alpha$ through {\bf R 5}${}_\alpha$ are satisfied
for every $\alpha < \lambda$. We then declare $\ell(i) = \lambda$, and
we define $$\sigma_\lambda^h \colon ({\cal M}_\lambda^{\bar {\cal T}})^h \rightarrow
{\cal N}^h$$ by setting $$\sigma_\lambda^h(x) = \sigma_\beta^h((\pi_{\beta
\lambda}^h)^{-1}(x)) {\rm , \ where \ } x \in {\rm ran}(\pi_{\beta
\lambda}^h).$$ This works by $\forall \alpha<\lambda \ ( \ $ {\bf R 2}${}_{\alpha}$
and {\bf R 4}${}_{\alpha}$ (b) $\ ) \ $. It is easy to check that
{\bf R 1}${}_{\lambda}$ to
{\bf R 5}${}_{\lambda}$ hold.

\bigskip
{\it Case 2.} Otherwise.

\bigskip
In this case, there is a largest $i$ such that $\ell(i) < \lambda$, and we are given
$${\vec {\cal V}}^h {\upharpoonright} i =
({\cal V}^h_j \colon j < i) {\rm , \ } {\cal V}^h_i {\upharpoonright} \lambda
{\rm , \ } (\sigma^h_\beta \colon \beta < \lambda)  {\rm , \ and \ }
(\rho_j \colon j \leq i)$$
in such a way that {\bf R 1}${}_\alpha$ to {\bf R 5}${}_\alpha$ are satisfied
for every $\alpha < \lambda$. We then continue with copying
${\bar {\cal T}} {\upharpoonright}
[\ell(i),\lambda]$ onto ${\cal N}^h$,
getting ${\cal V}^h_i {\upharpoonright} \lambda+1$ and the copy map $\sigma_\lambda^h$.
This works by
$\forall \alpha<\lambda$
{\bf R 2}${}_{\alpha}$ (a). We leave the easy details to the reader.
It is straighforward to check that
{\bf R 1}${}_{\lambda}$ to
{\bf R 5}${}_{\lambda}$ hold.

\bigskip
\hfill $\square$ (\ref{it-of-bicephalus})

\bigskip
In practice we'll always know that the hypothesis of
\ref{it-of-bicephalus} is satisfied. An example is given in the proof of the
following.

%
%
\begin{lemma}\label{on-the-sequence} $( \ \lnot \ \zerohandgrenade \ )$
Let $\kappa$ be a limit cardinal in $V$
such that for all $\lambda < \kappa$ we have that
$$\lambda {\rm \ is \ } <\kappa{\rm -strong \ in \ } K^c \Rightarrow
\lambda {\rm \ is \ } <{\rm OR}{\rm -strong \ in \ } K^c {\rm , }$$
and $\kappa \in C_0$.
Let ${\cal M} \trianglerighteq {\cal
J}^{K^c}_\kappa$ be a premouse. Suppose that $1 \leq n < \omega$ is such that

(a) $\rho_n({\cal M}) \leq \kappa < \rho_{n-1}({\cal M})$,

(b) ${\cal M}$ is $n$-sound above $\kappa$ (i.e., ${\cal M}$ is $(n-1)$-sound, and
${\cal M}^{n-1}$ is generated by $h^{n-1}_{\cal M}$ from $\kappa \cup \{ p_{{\cal
M},n} \}$), and

(c) ${\cal M}$ is $(n-1)$-iterable.

Then ${\cal M} \triangleleft K^c$, i.e., ${\cal M}$ is an initial segment of
$K^c$.
\end{lemma}

{\sc Proof.} We first verify the following.

\bigskip
{\bf Claim 1.} $\rho_n({\cal M}) \geq \kappa$, and hence ${\cal M}$ is $n$-sound.

\bigskip
{\sc Proof.} Suppose not. Then $\rho_n({\cal M}) < \kappa$, so that
${\frak C}_n({\cal M})$ has size $< \kappa$. Note that ${\frak C}_n({\cal M})$ is
$(n-1)$-iterable and $n$-sound. By $\kappa \in C_0$, there is an $(n-1)$-
maximal tree ${\cal T}$ on ${\frak C}_n({\cal M})$
with ${\cal D}^{\cal T} \cap (0,\infty]_T = \emptyset$ and
such that ${\cal M}_\infty^{\cal
T}$ is a non-simple iterate of ${\cal J}_\kappa^{K^c}$ (it is non-simple by
$\rho_n({\cal M}) < \kappa$). Moreover,
the $(n-1)$-maximal coiteration $({\bar {\cal U}},{\cal U})$ of
${\frak C}_n({\cal M})$, ${\cal M}$ is simple on both sides and produces some
common coiterate $Q = {\cal M}_\infty^{\bar {\cal U}} =
{\cal M}_\infty^{\cal U}$. Using $\pi = \pi^{\cal T}_{0 \infty}$, we
may copy ${\bar {\cal U}}$ onto ${\cal M}_\infty^{\cal T}$, getting an iteration
${\bar {\cal U}}^\pi$ of ${\cal M}_\infty^{\cal T}$ together with a last copy map
$$\sigma \colon Q \rightarrow Q'' = {\cal M}^{{\bar {\cal U}}^\pi}_\infty.$$
Now on the one hand we have
$$\sigma \circ \pi^{\cal U}_{0 \infty} \colon {\cal M}
\rightarrow Q'' {\rm , }$$
and on the other hand, we have that
$Q''$ is a non-simple iterate of ${\cal M}$ (as ${\cal M}_\infty^{\cal T}$ is).
This contradicts the Dodd-Jensen Lemma
(cf. \cite[Lemma 5.3]{FSIT}).

\bigskip
\hfill $\square$ (Claim 1)

\bigskip
We'll only need that $\rho_1({\cal M}) \geq \kappa$ in what follows.

We now use \ref{overlaps} to show that the
phalanx ${\cal P} = ((K^c,{\cal M}),\kappa)$ is normally $(n-1)$-iterable.
For suppose ${\cal U}$ to be a putative normal tree on ${\cal P}$. By
\ref{phalanx-iterations},
we can write $${\cal U} = {\cal U}_0{}^\frown{\cal U}_1$$ where
${\cal U}_0$ is an iteration of
${\cal M}$ above $\kappa$, and ${\cal U}_1$ is an iteration of
$K^c$ except for the fact that the first extender, call it $F$,
used for building ${\cal U}_1$
comes from the last model of ${\cal U}_0$, ${\rm c.p.}(F) < \kappa$, and the index
of $F$ is larger than $\kappa$
(possibly,
${\cal U}_0 = \emptyset$, or ${\cal U}_1 = \emptyset$). But by \ref{overlaps},
$F$ is countably complete, and hence ${\cal U}_1$ is well-behaved by standard
arguments.

Now let ${\cal U}$, ${\cal T}$ be the iteration trees arising from the
comparison of ${\cal P}$ with $K^c$, where we understand
${\cal U}$ to be $(n-1)$-maximal. We have to verify the following.

\bigskip
{\bf Claim 2.} ${\rm root}^{\cal U}(\infty) = 0$, i.e.,
${\cal M}_\infty^{\cal U}$, the last model of ${\cal U}$, sits above ${\cal M}$.

\bigskip
Before turning to its proof, let us first show that Claim 2 implies
\ref{on-the-sequence}. Note that by \ref{phalanx-iterations}, then,
${\cal U}$ is an iteration
of ${\cal M}$.

Suppose that ${\cal U}$ is non-trivial. Then ${\cal M}^{\cal U}_\infty$
is not sound, which by
\ref{Kc-is-universal} gives that we must have ${\cal M}^{\cal U}_\infty =
{\cal M}^{\cal T}_\infty$. This gives a standard contradiction if ${\cal
D}^{\cal U} \cap (0,\infty]_U \not=
\emptyset$ (cf. the proof of Claim 4 in the proof of \cite[Lemma 6.2]{FSIT}).
Hence we have to have ${\cal
D}^{\cal U} \cap (0,\infty]_U =
\emptyset$. As ${\cal T}$ only uses extenders with indices $\geq \kappa$ we'll
clearly have $\rho_\omega({\cal M}^{\cal T}_\infty) \geq \kappa$. But then
$\rho_\omega({\cal M}^{\cal U}_\infty) \geq \kappa$, and so
$${\frak C}_\omega({\cal M}^{\cal U}_\infty) =
{\frak C}_n({\cal M}^{\cal U}_\infty) = {\cal M} {\rm , }$$
and again we get a standard contradiction (as in
the proof of Claim 4 in the proof of \cite[Lemma 6.2]{FSIT}).

We have shown that ${\cal U}$ has to be trivial; that is, ${\cal M} \trianglerighteq
{\cal M}^{\cal T}_\infty$. But it can then easily be verified that
${\cal T}$ has to be trivial, too. This implies
${\cal M} \triangleleft K^c$ as desired.

\bigskip
{\sc Proof} of Claim 2.
Suppose not. That is, suppose that ${\rm root}^{\cal U}(\infty) = -1$, i.e., that
the last model of ${\cal U}$ sits above $K^c$.
By \ref{Kc-is-universal} and the Dodd-Jensen Lemma applied to
$\pi_{-1 \infty}^{\cal U}$ we easily get that
${\cal D}^{\cal U} \cap [-1,\infty]_U = {\cal D}^{\cal T} \cap [0,\infty]_U =
\emptyset$, and
$${\cal M}_\infty^{\cal U} = {\cal M}_\infty^{\cal T} {\rm , \ call \ it \ }
Q {\rm . }$$
In particular,
we have maps $\pi_{-1 \infty}^{\cal U} \colon K^c \rightarrow Q$ and
$\pi_{0 \infty}^{\cal T} \colon K^c \rightarrow Q$.
Let $F^0 = E^{\cal U}_\alpha$ be the first extender used on
$(-1,\infty]_U$.
Set $\mu = {\rm c.p.}(F^0) =
{\rm c.p.}(\pi_{-1 \infty}^{\cal U}) < \kappa$. A simple trick gives the following.

\bigskip
{\bf Subclaim 1.} ${\rm c.p.}(\pi_{0 \infty}^{\cal T}) = \mu$.

\bigskip
{\sc Proof.}
We have to show that $\pi_{0 \infty}^{\cal T} \not= {\rm id}$, and that $\mu$ is the
critical point of $\pi_{0 \infty}^{\cal T}$.
Let
$F^0 = E^{{\cal M}^{\cal U}_\alpha}_{\nu_0}$.
As ${\cal U}$ only uses extenders with indices $\geq \kappa$, and by Claim 1 above,
it is clear that
$\rho_1({\cal J}^{{\cal M}^{\cal U}_\alpha}_{\nu_0}) \geq \kappa$.
Hence \ref{strong-ISC+} immediately gives that $\mu$ is $<\kappa$-strong in
${\cal J}^{{\cal M}^{\cal U}_\alpha}_{\nu_0}$. But
${\cal J}^{{\cal M}^{\cal U}_\alpha}_{\kappa} = {\cal J}^{K^c}_\kappa$, so
$\mu$ is $<\kappa$-strong in $K^c$. By our assumption on $\kappa$, this means that
$\mu$ is $<{\rm OR}$-strong in $K^c$. On the other hand, \ref{not-strong} gives that
$\mu$ is not $<{\rm OR}$-strong in $Q$.

This already implies that $\pi_{0 \infty}^{\cal T} \not= {\rm id}$.
Let $F^1 = E^{\cal T}_\beta = E^{{\cal M}^{\cal
T}_\beta}_{\nu_1}$ be the first extender used on $(0,\infty]_T$.

Let us first assume that $\pi_{0 \infty}^{\cal T} {\upharpoonright} \mu \not= {\rm id}$,
i.e., that ${\rm c.p.}(F^1) < \mu$.
As ${\cal T}$ only uses extenders with indices $\geq \kappa$, we clearly have that
$\rho_1({\cal J}^{{\cal M}^{\cal
T}_\beta}_{\nu_1}) \geq \kappa$. Thus, again by using \ref{strong-ISC+} and
\ref{not-strong},
we end up with getting ${\rm c.p.}(F^1)$ is $<{\rm OR}$-strong in $K^c$, but ${\rm c.p.}(F^1)$ is not
$<{\rm OR}$-strong in $Q$. However, consider now
$\pi^{\cal U}_{-1 \infty} \colon K^c \rightarrow
Q$. As $\pi^{\cal U}_{-1 \infty} {\upharpoonright} {\rm c.p.}(F^1)+1 = {\rm id}$, the fact that
${\rm c.p.}(F^1)$ is $<{\rm OR}$-strong in $K^c$ implies that ${\rm c.p.}(F^1)$ remains $<{\rm OR}$-strong in
$Q$. Contradiction!

A symmetric argument shows that we cannot have $\pi_{0 \infty}^{\cal T}
{\upharpoonright} \mu+1 = {\rm id}$.

\bigskip
\hfill $\square$ (Subclaim 1)

\bigskip
We could also have used the Dodd-Jensen Lemma to show $\pi_{0 \infty}^{\cal T}
{\upharpoonright} \mu = {\rm id}$ in the proof of Subclaim 1. However, we chose not to do so
in order to exhibit the symmetry of the argument.

Now let $F^1 = E^{\cal T}_\beta$ still denote
the first extender used on $(0,\infty]_T$. Subclaim 1
says that ${\rm c.p.}(F^1) = \mu = {\rm c.p.}(F^0)$. Let, again, $$F^0 = E^{{\cal M}^{\cal
U}_\alpha}_{\nu_0} {\rm \ and \ } F^1 = E^{{\cal M}^{\cal
T}_\beta}_{\nu_1}.$$ The rest of this proof is entirely symmetric, so that we may
assume without loss of generality that $\nu_0 \leq \nu_1$. Both $\nu_0$ and $\nu_1$ are cardinals in
$Q$, which implies that $\nu_0 < \nu_1 \Rightarrow \nu_0$ is a cardinal in
${\cal J}^{{\cal M}^{\cal
T}_\beta}_{\nu_1}$. It is hence clear that we may derive from
$${\cal J}^{{\cal M}^{\cal
U}_\alpha}_{\nu_0} {\rm \ , \ } {\cal J}^{{\cal M}^{\cal
T}_\beta}_{\nu_1}$$ a generalized prebicephalus, call it ${\cal N}$. The following
is a straightforward consequence of \ref{goodness}, \ref{overlaps}, and
\ref{it-of-bicephalus}:

\bigskip
{\bf Subclaim 2.} ${\cal N}$ is iterable.

\bigskip
{\sc Proof.} It suffices to verify that ${\cal N}$ satisfies the hypothesis of
\ref{it-of-bicephalus}. Specifically, we want to verify
that our current $\kappa$ may
serve as the $\kappa$ in the statement of
\ref{it-of-bicephalus}. However, this is trivial by virtue of \ref{goodness} and
\ref{overlaps}!

\bigskip
\hfill $\square$ (Subclaim 2)

\bigskip
Now \ref{bicephali-trivialize}
yields that
$E^{\cal U}_\alpha =
F^0 = F^1 = E^{\cal T}_\beta$. This is a contradiction,
as this can't happen in a comparison!

\bigskip
\hfill $\square$ (Claim 2)

\hfill $\square$ (\ref{on-the-sequence})

\bigskip
Using \ref{maximality} of the next section we could have shown \ref{on-the-sequence}
for all $\kappa \in C_0$.
\section{Maximality of $K^c$.}
\setcounter{equation}{0}

\begin{lemma}\label{maximality} $( \ \lnot \ \zerohandgrenade \ )$ (Maximality of $K^c$)
Let ${\cal T}$ be a normal iteration tree on $K^c$ with last model $Q = {\cal M}^{\cal
T}_\infty$. Suppose that $(0,\infty]_T \cap {\cal D}^{\cal T} = \emptyset$.
Let $\mu < \kappa \leq \nu$ be cardinals of $Q$, and suppose that
$\pi^{\cal T}_{0 \infty} {\upharpoonright} \kappa = {\rm id}$.
Then there is no extender $F$ with critical point $\mu$
such that the following hold.

$\bullet \ $ ${\cal N} = ({\cal J}^Q_\nu;F)$ is a premouse,

$\bullet \ $ if ${\cal U}$ is an iteration tree on ${\cal N}$ with last model ${\cal
N}^\star = {\cal
M}_\infty^{\cal U}$, $(0,\infty]_U \cap {\cal D}^{\cal U} = \emptyset$,
and $\pi^{\cal U}_{0 \infty} {\upharpoonright} \kappa = {\rm id}$
then the top
extender of ${\cal
N}^\star$ is countably complete, and

$\bullet \ $ $\rho_1({\cal N}) > \mu$.
\end{lemma}

{\sc Proof.} Suppose not, and let ${\cal T}$, $Q$, $\mu$, $\kappa$, $\nu$, $F$, and
${\cal N}$ be as
in the statement of \ref{maximality}. The proof will actually be by induction on
$\mu$. That is, we assume that $\mu$ is least such that there are
${\cal T}$, $Q$, $\kappa$, $\nu$, $F$, and
${\cal N}$ as
in the statement of \ref{maximality}, and we fix such
${\cal T}$, $Q$, $\kappa$, $\nu$, $F$, and
${\cal N}$.

We first aim to define the {\em pull back} of ${\cal
N}$ via $\pi^{\cal T}_{0 \infty}$.
Set $\pi = \pi^{\cal T}_{0 \infty}$.
Let $\lambda = F(\mu)$, i.e., $\lambda$ is the largest cardinal of ${\cal N}$.
Of course, $\lambda$ is an inaccessible cardinal of both ${\cal N}$ and $Q$.
Let ${\bar \lambda} = \pi^{-1} {\rm " } \lambda$, so that ${\bar \lambda}$ is least with
$\pi({\bar \lambda}) \geq \lambda$.
We define a Dodd-Jensen extender ${\bar F}$
by setting
$$X \in {\bar F}_a \ \Leftrightarrow \ \pi(a) \in F(X)$$ for $a \in [{\bar \lambda}]^{<
\omega}$ and $X \in {\cal P}([\mu]^{{\rm Card}(a)}) \cap K^c$.
I.e., ${\bar F}$ is a $(\mu,{\bar \lambda})$-extender over $K^c$
in the sense of \cite[Definition 1.0.1]{FSIT}. ${\bar F}$ 
inherits the countable completeness
from $F$. Let $$i_{\bar F} \ \colon \
{\cal J}^{K^c}_{\mu^{+K^c}} \rightarrow_{\bar F} {\bar {\cal M}}$$ be the $\Sigma_0$
ultrapower map, and set $${\cal M} = ({\bar {\cal M}},i_{\bar F} {\upharpoonright} {\cal
P}(\mu) \cap K^c) {\rm ; }$$
that is, ${\bar {\cal M}}$ is the model theoretic reduct of ${\cal M}$, where the
latter has the additional top extender $i_{\bar F} {\upharpoonright} {\cal
P}(\mu) \cap K^c$.
We call ${\cal M}$ the pull back of ${\cal N}$ via $\pi$.
We shall also write ${\bar F}$ for $i_{\bar F} {\upharpoonright} {\cal
P}(\mu) \cap K^c$.
We can define a cofinal $\Sigma_\omega$-elementary map
$$\pi' \ \colon \ {\bar {\cal M}}
\rightarrow {\cal J}_\nu^Q {\rm \ \ by \ setting }$$
$$[a,f]^{{\cal J}^{K^c}_{\mu^{+K^c}}}_{\bar F}
\mapsto [\pi(a),f]^{{\cal J}^{K^c}_{\mu^{+K^c}}}_F$$ for
$a \in [{\bar \lambda}]^{<
\omega}$ and $f \colon [\mu]^{{\rm Card}(a)} \rightarrow {\cal J}^{K^c}_{\mu^{+K^c}}$, $f
\in K^c$. We'll in fact have that $$\pi' \ \colon \ {\cal M} \rightarrow {\cal N}$$ is
$\Sigma_0$-elementary (and hence $\Sigma_1$-elementary, as $\pi'$ is cofinal),
because $\pi'({\bar F} \cap x) = F \cap \pi'(x)$ for all $x \in {\cal M}$.
This readily implies that ${\cal M}$ is a premouse: the initial segment condition for
${\cal M}$ is true as $C^{\cal M}_{{\cal M} \cap {\rm OR}} \not= \emptyset \Rightarrow
C^{\cal N}_{{\cal N} \cap {\rm OR}} \not= \emptyset$; but
$C^{\cal N}_{{\cal N} \cap {\rm OR}} = \emptyset$.

The following statements are easy to verify.

\bigskip
$\bullet \ $ ${\bar \lambda}$ is a limit cardinal of $K^c$,

$\bullet \ $ ${\cal J}^{\cal M}_{\bar \lambda} = {\cal J}^{K^c}_{\bar \lambda}$,

$\bullet \ $ $\pi' {\upharpoonright} {\cal J}^{\cal M}_{\bar \lambda} =
\pi {\upharpoonright} {\cal J}^{K^c}_{\bar \lambda}$, and

$\bullet \ $ $\pi'({\bar F}(\mu)) = F(\mu) = \lambda$.

\bigskip
{\bf Claim 1.} There is no $\alpha \geq {\bar \lambda}$ such that $E_\alpha^{K^c}
\not= \emptyset$, and ${\rm c.p.}(E_\alpha^{K^c}) \in [\mu,{\bar \lambda})$.

\bigskip
{\sc Proof.} Suppose otherwise. Set ${\bar \mu} = {\rm c.p.}(E_\alpha^{K^c})$. By
\ref{strong-ISC+}, then, ${\bar \mu}$ is $< {\bar \lambda}$-strong in $K^c$ as
witnessed by ${\vec E}^{K^c}$. By the elementarity of $\pi$ we hence have
$\pi({\bar \mu}) \in [\mu,\pi({\bar \lambda}))$ is $< \pi({\bar \lambda})$-strong
in $Q$ as
witnessed by ${\vec E}^Q$.
By the definition of ${\bar \lambda}$, $\pi({\bar \mu}) < \lambda$.
So, trivially, $\pi({\bar \mu})$ is $< \lambda$-strong
in $Q$ as
witnessed by ${\vec E}^Q$, and thus $\pi({\bar \mu})$ is $< \lambda$-strong
in ${\cal N}$ as
witnessed by ${\vec E}^{\cal N}$. But then ${\cal N}$ is easily seen not to be below
$\zerohandgrenade$.

\bigskip
\hfill $\square$ (Claim 1)

\bigskip
{\bf Claim 2.} There is no $\alpha \in [{\bar \lambda},{\cal M} \cap {\rm OR})$
such that $E_\alpha^{\cal M}
\not= \emptyset$, and ${\rm c.p.}(E_\alpha^{\cal M}) \in [\mu,{\bar \lambda})$.

\bigskip
{\sc Proof.} Suppose otherwise. Set ${\bar \mu} = {\rm c.p.}(E_\alpha^{\cal M})$. By
\ref{strong-ISC+}, then, ${\bar \mu}$ is $< {\bar \lambda}$-strong in ${\cal M}$ as
witnessed by ${\vec E}^{\cal M}$,
and thus
${\bar \mu}$ is $< {\bar \lambda}$-strong in $K^c$ as witnessed by ${\vec E}^{K^c}$.
This then gives a contradiction as in the proof of
Claim 1.

\bigskip
\hfill $\square$ (Claim 2)

\bigskip
{\bf Claim 3.} ${\cal J}^{\cal M}_{{\bar \lambda}^{+{\cal M}}} \triangleleft K^c$.

\bigskip
{\sc Proof.} This is shown by coiterating $K^c$ with
${\cal J}^{\cal M}_{{\bar \lambda}^{+{\cal M}}}$, or rather with the phalanx
${\cal P} = ((K^c,{\cal J}^{\cal M}_{{\bar \lambda}^{+{\cal M}}}),{\bar \lambda})$.

\bigskip
{\bf Subclaim 1.} ${\cal P}$ is iterable.

\bigskip
{\sc Proof.} By standard arguments and \ref{phalanx-iterations},
it is enough to see that if ${\cal V}$ is an iteration of
successor length of
${\cal J}^{\cal M}_{{\bar \lambda}^{+{\cal M}}}$ which only uses extenders with
critical point $\geq {\bar \lambda}$ then any $E^{{\cal M}^{\cal V}_\infty}_\rho
\not= \emptyset$
with $\rho \geq {\bar \lambda}$ and ${\rm c.p.}(E^{{\cal M}^{\cal V}_\infty}_\rho) < {\bar
\lambda}$ is countably complete. However, using $\pi'$ we may copy ${\cal V}$ onto
${\cal N}$, getting an iteration tree ${\cal U}$ on ${\cal N}$.
Let $\pi'_\infty \colon {\cal M}^{\cal V}_\infty \rightarrow {\cal M}^{\cal U}_\infty$
be the last copy map.

${\cal M}^{\cal V}_\infty$ inherits the property expressed by
Claim 2 above, that is, we must have
${\rm c.p.}(E^{{\cal M}^{\cal V}_\infty}_\rho) < \mu$. But then
${\rm c.p.}(E^{{\cal M}^{\cal U}_\infty}_{\pi'_\infty(\rho)}) < \mu$, too, and hence
$E^{{\cal M}^{\cal U}_\infty}_{\pi'_\infty(\rho)}$ is countably complete by
\ref{goodness}. But this implies that $E^{{\cal M}^{\cal V}_\infty}_\rho$
is countably complete.

\bigskip
\hfill $\square$ (Subclaim 1)

\bigskip
Now let ${\cal W}$, ${\cal V}$ denote the iteration trees arising from the comparison
of $K^c$ with ${\cal P}$. By \ref{Kc-is-universal} and standard arguments it suffices
to see that the last model of ${\cal V}$ sits above
${\cal J}^{\cal M}_{{\bar \lambda}^{+{\cal M}}}$. Let us suppose not. We shall derive
a contradiction.

By the Dodd-Jensen lemma and \ref{Kc-is-universal} we'll then have that
${\cal M}^{\cal W}_\infty = {\cal M}^{\cal V}_\infty$, ${\cal D}^{\cal W} \cap
[0,\infty]_W = \emptyset$, ${\cal D}^{\cal V} \cap
[0,\infty]_V = \emptyset$, and $\pi^{\cal W}_{0 \infty}$ is lexicographically $\leq
\pi^{\cal V}_{0 \infty}$. In particular, we must have either
${\rm c.p.}(\pi^{\cal W}_{0 \infty}) = {\rm c.p.}(\pi^{\cal V}_{0 \infty})$ or else
$\pi^{\cal W}_{0 \infty} {\upharpoonright} {\rm c.p.}(\pi^{\cal V}_{0 \infty})
+1 = {\rm id}$.

\bigskip
{\it Case 1.} ${\rm c.p.}(\pi^{\cal W}_{0 \infty}) = {\rm c.p.}(\pi^{\cal V}_{0 \infty})$.

\bigskip
Let $\alpha+1$ be least in $(0,\infty]_W$, and let $\beta+1$
be least in $(0,\infty]_V$. Let $$F^0 = E_{\nu_0}^{{\cal M}^{\cal W}_\alpha} {\rm \
and \ } F^1 = E_{\nu_1}^{{\cal M}^{\cal V}_\beta}.$$ The case assumption says that
${\rm c.p.}(F^0) = {\rm c.p.}(F^1)$. It is easy to see that we may derive from
$${\cal J}_{\nu_0}^{{\cal M}^{\cal W}_\alpha} {\rm \ , \ }
{\cal J}_{\nu_1}^{{\cal M}^{\cal V}_\beta}$$ a generalized prebicephalus, call it ${\cal B}$.

\bigskip
{\bf Subclaim 2.} ${\cal B}$ is iterable.

\bigskip
{\sc Proof.} It suffices to verify that ${\cal B}$ satisfies the hypothesis of
\ref{it-of-bicephalus}. Specifically, we want to verify
that ${\rm c.p.}(\pi^{\cal W}_{0 \infty})^{+K^c}$ may
serve as the $\kappa$ in the statement of
\ref{it-of-bicephalus}. However, this shown by combining \ref{goodness}
with arguments which should be standard by now, and
with a copying argument as in the proof of
Subclaim 1 above.

\bigskip
\hfill $\square$ (Subclaim 2)

\bigskip
Now \ref{bicephali-trivialize}
yields that
$E^{\cal W}_\alpha =
F^0 = F^1 = E^{\cal U}_\beta$. This is a contradiction,
as this can't happen in a comparison!

\bigskip
{\it Case 2.} $\pi^{\cal W}_{0 \infty} {\upharpoonright} {\rm c.p.}(\pi^{\cal V}_{0 \infty})
+1 = {\rm id}$.

\bigskip
In this case we shall use our inductive hypothesis. In fact, luckily, this is easy to
do. Namely, arguments which again should be standard by now
yield that we have with ${\cal W}$, ${\cal M}_\infty^{\cal W}$,
${\rm c.p.}(\pi^{\cal V}_{0 \infty})$, ${\rm sup} (\{ \xi < \nu_1 \colon \pi^{\cal W}_{0
\infty}(\xi) = \xi \})$, $\nu_1$, $F^1$, and
${\cal J}_{\nu_1}^{{\cal M}^{\cal V}_\beta}$ a series of objects which are as ${\cal
T}$, $Q$, $\mu$, $\kappa$, $\nu$, and ${\cal N}$ in the statement of \ref{maximality}.
We leave the straightforward details to the reader.
However, by Claim 2 we'll have that ${\rm c.p.}(\pi^{\cal V}_{0 \infty}) < \mu$, which
contradicts the choice of $\mu$ being minimal!

\bigskip
\hfill $\square$ (Claim 3)

\bigskip
We'll now have to split the argument into three cases.

\bigskip
{\it Case 1.} ${\bar F}(\mu) = {\bar \lambda}$.

\bigskip
In this case Claim 3 immediately gives that ${\bar {\cal M}} \triangleleft K^c$.
(Recall that ${\bar {\cal M}}$ is the model theoretic
reduct of ${\cal M}$ obtained by removing the
top extender ${\bar F}$.) But now we have a contradiction with \ref{weak-max}.

\bigskip
{\it Case 2.} ${\bar F}(\mu) > {\bar \lambda}$ and $\pi({\bar \lambda}) > \lambda$.

\bigskip
In this case Claim 3 gives that we may define the lift up
$${\tilde {\cal M}} = {\rm Ult}_0({\cal M};\pi {\upharpoonright}
{\cal J}^{\cal M}_{{\bar \lambda}^{+{\cal
M}}}).$$

\bigskip
{\bf Claim 4.} ${\tilde {\cal M}}$ is transitive and $0$-iterable.

\bigskip
{\sc Proof.} We first note that $\pi$ is countably complete:

\bigskip
{\bf Fact.} Let $(a_n,X_n \colon n<\omega)$ be such that for all $n<\omega$
do we have
$a_n \in
[{\rm OR}]^{<\omega}$, $X_n \in {\cal P}([\eta_n]^{{\rm Card}(a_n)}) \cap K^c$ for some $\eta_n$,
and $a_n \in
\pi(X_n)$. Then there is an order preserving $\tau \colon
\bigcup_{n<\omega} a_n \rightarrow {\rm OR}$ with $\tau {\rm " } a_n \in X_n$ for all
$n<\omega$.

\bigskip
{\sc Proof.} This is a straightforward consequence of the proof of
\ref{normal-iterability}. For this purpose we may
consider ${\cal T}$ as a tree on ${\cal
J}_\theta^{K^c}$, for some $\theta$. We may and shall assume that $\theta$ is large
enough and regular so that we may pick $\sigma \colon {\bar H} \rightarrow H_\theta$
with ${\bar H}$ countable and transitive, and $\{ {\cal T} \} \cup \{a_n,\pi(X_n)
\colon n<\omega \}
\subset {\rm ran}(\sigma)$. Set ${\bar {\cal T}} =
\sigma^{-1}({\cal T})$. The construction from the proof of
\ref{normal-iterability} will then give us an embedding $$\sigma_\infty \ \colon
\ {\cal M}_\infty^{\bar {\cal T}} \rightarrow {\cal
J}_\theta^{K^c}$$ such that for all $X$ with $\pi(X) \in {\rm ran}(\sigma)$
do we have that $$\sigma_\infty \circ \sigma^{-1} \circ \pi(X) = X.$$
We hence have in $$\tau = \sigma_\infty \circ \sigma^{-1} {\upharpoonright}
\bigcup_{n<\omega} a_n$$ a function as desired.

\bigskip
\hfill $\square$ (Fact)

\bigskip
This implies Claim 4 by standard arguments.

\bigskip
\hfill $\square$ (Claim 4)

\bigskip
Notice that
${\tilde {\cal M}}$ is a premouse (rather than a protomouse, as the top extender
of ${\cal M}$ has critical point $\mu$ and $\pi {\upharpoonright} \mu^{+{\cal M}} =
{\rm id}$). Moreover, $\nu$ (the height of ${\cal
N}$) is a cardinal in ${\tilde {\cal M}}$ because $\pi({\bar \lambda}) > \lambda$ and
so $\pi({\bar \lambda}) > \nu$, and $${\cal B} =
({\cal N},{\tilde {\cal M}})$$ is a generalized prebicephalus.
Notice that we can't have that ${\cal N} = {\tilde {\cal M}}$, just because $\nu =
{\cal N} \cap {\rm OR} < {\tilde {\cal M}}$.
Thus, by \ref{bicephali-trivialize}, in order to reach a contradiction it
suffices to verify the following.

\bigskip
{\bf Claim 5.} ${\cal B}$ is iterable.

\bigskip
{\sc Proof.} We shall apply \ref{it-of-bicephalus}. Specifically, we plan on letting
the current $\kappa$ play the r\^ole of the $\kappa$ in the statement of
\ref{it-of-bicephalus}. We have to show that any iteration tree ${\cal V}$ on ${\cal
N}$, or on ${\tilde {\cal M}}$, meets the hypothesis in the statement of
\ref{it-of-bicephalus}. However, notice that this is clear for ${\cal N}$
by \ref{goodness} and our
assumptions on ${\cal N}$. We are hence left with having to verify the following.

\bigskip
{\bf Subclaim 3.} Let ${\cal V}$ be an iteration tree
on ${\tilde {\cal M}}$ with last model ${\cal
M}^{\cal V}_\infty$ such that for some ${\tilde {\cal M}}$-cardinal
$\rho \leq \kappa$ we have that either

(a) $\rho < \kappa$, and ${\cal V}$ lives on ${\cal J}^{{\tilde {\cal M}}}_\kappa$ and
is above $\rho$, but ${\cal V}$ doesn't use extenders which are
total on ${\tilde {\cal M}}$, or else

(b) $\rho = \kappa$, and ${\cal V}$ is above $\rho$,

\noindent we have that

($\star$) if $F =
E^{{\cal
M}^{\cal V}_\infty}_\nu \not= \emptyset$ is
such that $\nu > \rho$ and ${\rm c.p.}(F) < \rho$, then $F$ is countably complete.

\bigskip
{\sc Proof.}
Let ${\cal V}$ be as in the statement of this Subclaim, with some $\rho \leq \kappa$.
Let $F =
E^{{\cal
M}^{\cal V}_\infty}_\nu \not= \emptyset$ be
such that $\nu > \rho$ and ${\rm c.p.}(F) < \rho$, and let $(a_n,X_n \colon n<\omega)$
be such that $a_n \in F(X_n)$ for all $n<\omega$.
Let $\sigma \colon {\bar H} \rightarrow H_\theta$ for some large enough regular
$\theta$ be such that ${\bar H}$ is countable and transitive, and
$\{ {\cal V} , F \} \cup \{a_n,X_n \colon n<\omega \}
\subset {\rm ran}(\sigma)$. Set ${\bar {\cal V}} =
\sigma^{-1}({\cal V})$.
By the above Fact (in the proof of Claim 4) and standard
arguments there is some $${\bar \sigma} \ \colon \
{\cal M}_0^{\bar {\cal V}} = \sigma^{-1}({\tilde {\cal M}}) \rightarrow {\cal
M}.$$ Recall that we also have a cofinal map
$$\pi' \ \colon \ {\cal M} \rightarrow_{\Sigma_1}
{\cal N}.$$
We may hence copy ${\bar {\cal V}}$ onto ${\cal N}$
using $\pi' \circ {\bar \sigma}$, getting an iteration tree
${\cal V}'$ on ${\cal N}$, together with a last copy map $$\sigma_\infty \ \colon \
{\cal M}_\infty^{\bar {\cal V}} \rightarrow {\cal M}_\infty^{{\cal V}'}.$$ Now by
\ref{goodness} and
the
assumptions on ${\cal N}$ we'll have that $\sigma_\infty \circ \sigma^{-1}(F)$ is
countably complete. We may hence pick some order preserving
$$\tau \colon \bigcup_{n<\omega}
\sigma_\infty \circ \sigma^{-1}(a_n) \rightarrow
{\rm c.p.}(\sigma_\infty \circ \sigma^{-1}(F)) {\rm \ with }$$
$$\tau \circ \sigma_\infty \circ \sigma^{-1} {\rm " } a_n \in
\sigma_\infty \circ \sigma^{-1}(X_n) {\rm \ for \ all \ } n<\omega.$$
However, $$\sigma_\infty \circ \sigma^{-1}(X_n) = X_n {\rm , }$$
as ${\rm c.p.}(F) < \rho$, and ${\cal V}$ is above $\rho$. We therefore have
$$\tau \circ \sigma_\infty \circ \sigma^{-1} {\rm " } a_n \in
X_n {\rm \ for \ all \ } n<\omega {\rm ,}$$ and hence we have in
$\tau \circ \sigma_\infty \circ \sigma^{-1}$ a function as desired.

\bigskip
\hfill $\square$ (Subclaim 3)

\hfill $\square$ (Claim 5)

\bigskip
{\it Case 3.} ${\bar F}(\mu) > {\bar \lambda}$ and $\pi({\bar \lambda}) = \lambda$.

\bigskip
Set ${\tilde \lambda} = {\rm sup} \ \pi
{\rm " } {\bar \lambda}$. It is easy to see that $\pi$
can't be continuous at ${\bar \lambda}$ in this case, i.e., that

\bigskip
$\bullet \ $
${\tilde \lambda} <
\lambda$,

\bigskip
\noindent because otherwise we would have to have
$\pi'({\bar \lambda}) \geq \lambda$, whereas $\pi'({\bar F}(\mu)) = \lambda$.
Moreover, as $\lambda$ is inaccessible in $Q$, we must have that $\pi^{-1}(\lambda) =
{\bar \lambda}$ is inaccessible in $K^c$. This implies, because $\pi$ is an
iteration map which is discontinuous at ${\bar \lambda}$, that
${\bar \lambda}$ is measurable in $K^c$, and that in fact we can
write $${\cal T} = {\cal T}_0{}^\frown{\cal T}_1$$ where
$\pi_{0 \infty}^{{\cal T}_0}$ is continuous at ${\bar \lambda}$,
${\rm c.p.}(\pi^{{\cal T}_1}_{0 \infty}) =
\pi^{{\cal T}_0}_{0 \infty}({\bar \lambda}) = {\tilde \lambda}$,
and $\pi = \pi_{0 \infty}^{{\cal T}_1} \circ \pi_{0 \infty}^{{\cal T}_0}$.
(Let $\rho$ be the least ${\bar \rho}
\in [0,\infty)_T$ with ${\rm c.p.}(\pi_{{\bar \rho} \infty}^{\cal T}) = \pi_{0 {\bar
\rho}}^{\cal T}({\bar \lambda})$. Then $\pi_{0 \infty}^{{\cal T}_0} = \pi_{0
\rho}^{\cal T}$, and $\pi_{0 \infty}^{{\cal T}_1} = \pi_{\rho
\infty}^{\cal T}$.)
Let us write $\pi_0 =
\pi^{{\cal T}_0}_{0 \infty}$ and $\pi_1 =
\pi^{{\cal T}_1}_{0 \infty}$. We have that $\pi_0({\bar \lambda}) = {\tilde \lambda}$,
and $\pi_1({\tilde \lambda}) = \lambda$.

\bigskip
{\bf Claim 6.} ${\cal J}^{\cal M}_{{\bar \lambda}^{+{\cal M}}} =
{\cal J}^{K^c}_{{\bar \lambda}^{+{K^c}}}$.

\bigskip
{\sc Proof.} We already know that ${\cal J}^{\cal M}_{{\bar \lambda}^{+{\cal M}}}
\trianglelefteq
{\cal J}^{K^c}_{{\bar \lambda}^{+{K^c}}}$. It is easy to see that
${\cal M}_\infty^{{\cal T}_0}$ and $Q$ agree up to ${\tilde
\lambda}^{+{{\cal M}_\infty^{{\cal T}_0}}} = {\tilde \lambda}^{+Q}$ and that

\begin{equation}\label{claim6.1}
\pi_0
{\upharpoonright} {\cal J}^{K^c}_{{\bar \lambda}^{+K^c}} \ \colon
\ {\cal J}^{K^c}_{{\bar \lambda}^{+K^c}} \rightarrow
{\cal J}^{Q}_{{\tilde \lambda}^{+Q}}
\end{equation}

\noindent
is exactly the lift up of ${\cal J}^{K^c}_{{\bar \lambda}^{+K^c}}$
by $\pi_0 {\upharpoonright} {\cal J}^{K^c}_{\bar \lambda} =
\pi {\upharpoonright} {\cal J}^{K^c}_{\bar \lambda}$,
which in turn is given by letting ${\cal T}_0$ act on
${\cal J}^{K^c}_{{\bar \lambda}^{+K^c}}$.

We'll now use the interpolation technique. Let $$\sigma_0 \colon
{\cal M} \rightarrow {\cal L} = {\rm Ult}_0({\cal M};\pi {\upharpoonright} {\cal J}_{\bar
\lambda}^{K^c})$$ be the lift up of ${\cal M}$ by
$\pi {\upharpoonright} {\cal J}_{\bar
\lambda}^{K^c}$.
Notice that
$\rho_1({\cal L}) \leq {\tilde \lambda}$ and ${\cal L}$ is ${\tilde \lambda}$-sound,
because every element of ${\cal L}$ can be written in the form $[a,i_{\bar F}(f)(b)]$
where $a \in [{\tilde
\lambda}]^{< \omega}$, $b \in [{\bar
\lambda}]^{< \omega}$, and
$f \colon [\mu]^{{\rm Card}(b)} \rightarrow {\cal J}^{K^c}_{\mu^{+K^c}}$ with
$f \in {\cal J}^{K^c}_{\mu^{+K^c}}$. Of
course, $\sigma_0({\bar \lambda}) = {\tilde \lambda}$.

We may define $$\sigma_1 \colon {\cal L} \rightarrow {\cal N}$$ by
letting $$[a,f]^{\cal M}_{\pi {\upharpoonright}
{\cal J}_{\bar \lambda}^{K^c}} \mapsto \pi'(f)(a)$$ for $a \in [{\tilde
\lambda}]^{< \omega}$,
and $f \colon [\eta]^{{\rm Card}(a)} \rightarrow {\cal M}$ for some $\eta < {\bar \lambda}$
with $f \in {\cal M}$ and $a \in \pi({\rm dom}(f))$.
It
is straightforward to see that $\sigma_1$ is $\Sigma_0$-elementary, that $\sigma_1
{\upharpoonright} {\tilde \lambda} = {\rm id}$,
and that
$$\pi' = \sigma_1 \circ \sigma_0 {\rm . }$$ Of course ${\tilde
\lambda}$ is
inaccessible in both ${\cal L}$ and ${\cal N}$.

Clearly, ${\cal J}^{\cal N}_{{\tilde \delta}^{+{\cal N}}} =
{\cal J}^{\cal Q}_{{\tilde \delta}^{+{\cal Q}}}$. Let us assume that we also have

\begin{equation}\label{claim6.2}
{\cal J}^{\cal L}_{{\tilde \delta}^{+{\cal L}}} =
{\cal J}^{\cal Q}_{{\tilde \delta}^{+{\cal Q}}}.
\end{equation}

\noindent Then we'll have that

\begin{equation}\label{claim6.3}
\sigma_0 {\upharpoonright}
{\cal J}^{\cal M}_{{\bar \lambda}^{+{\cal M}}} \ \colon
\ {\cal J}^{\cal M}_{{\bar \lambda}^{+{\cal M}}} \rightarrow
{\cal J}^{\cal L}_{{\tilde \lambda}^{+{\cal
L}}} =
{\cal J}^{Q}_{{\tilde \lambda}^{+Q}}
\end{equation}

\noindent is the map obtained by
taking the ultrapower of ${\cal J}^{\cal M}_{{\bar \lambda}^{+{\cal M}}}$
by the long extender $\pi {\upharpoonright} {\cal J}^{K^c}_{\bar \lambda}$,
which of course is exactly what is obtained when we let
${\cal T}_0$ act on
${\cal J}^{\cal M}_{{\bar \lambda}^{+{\cal M}}}$.
Hence both $\pi_0
{\upharpoonright} {\cal J}^{K^c}_{{\bar \lambda}^{+K^c}}$
and $\sigma_0 {\upharpoonright}
{\cal J}^{\cal M}_{{\bar \lambda}^{+{\cal M}}}$ are maps induced by ${\cal T}_0$, and
they have the same target model, and ${\cal J}^{\cal M}_{{\bar \lambda}^{+{\cal M}}}
\trianglelefteq
{\cal J}^{K^c}_{{\bar \lambda}^{+{K^c}}}$.
This implies
that we must have ${\cal J}^{\cal M}_{{\bar \lambda}^{+{\cal M}}} =
{\cal J}^{K^c}_{{\bar \lambda}^{+{K^c}}}$.

In order to finish the proof of Claim 6 it hence suffices to verify (\ref{claim6.2}).
We first prove:

\bigskip
{\bf Subclaim 4.} There is no $E_\alpha^{\cal N} \not= \emptyset$ with
${\rm c.p.}(E_\alpha^{\cal N}) \in
[\mu,{\tilde \lambda})$ and $\alpha \in ({\tilde \lambda},\nu)$.

\bigskip
{\sc Proof.} Otherwise by \ref{strong-ISC+} we get that ${\bar \mu} =
{\rm c.p.}(E_\alpha^{\cal N})$ is $< {\tilde \lambda}$-strong in ${\cal N}$
as witnessed by ${\vec E}^{\cal N}$.
Using $\pi_0$, we then get that there is some ${\bar \mu} \in [\mu,{\bar \lambda})$
which is $< {\bar \lambda}$-strong in $K^c$ as witnessed by ${\vec E}^{K^c}$, which is
a contradiction as in the proof of Claim 2 above.

\bigskip
\hfill $\square$ (Subclaim 4)

\bigskip
Now (\ref{claim6.2})
is shown by applying the condensation lemma \cite[\S 8 Lemma 4]{NFS}
to $\sigma_1$.
Recall that $\rho_1({\cal L}) \leq {\tilde \lambda}$, ${\cal L}$ is ${\tilde
\lambda}$-sound, and $\sigma_1 {\upharpoonright} {\tilde \lambda} = {\rm id}$.
We may assume that $\sigma_1 \not= {\rm id}$ and ${\cal L} \not= {\cal N}$, as otherwise
(\ref{claim6.2}) is trivial. Let $\eta = {\rm c.p.}(\sigma_1) \geq {\tilde \lambda}$.
Well, if (b) or (c) in the conclusion of \cite[\S 8 Lemma 4]{NFS} 
were to hold then (in much
the same way as in the proof of \ref{strong-ISC}) we'd get that there is some
$E_\alpha^{\cal N} \not= \emptyset$ with ${\rm c.p.}(E_\alpha^{\cal N}) =
\mu$ and $\alpha \in ({\tilde
\lambda},\nu)$, contradicting Subclaim 4.
We hence must have that
(a) in the conclusion of \cite[\S 8 Lemma 4]{NFS} holds, i.e., that
${\cal L}$ is the $\eta$-core of ${\cal
N}$, and that $\sigma_1$ is the core map.

We are now finally going to use our assumption that $\rho_1({\cal N}) > \mu$
which we make in the statement of
\ref{maximality}. Let ${\cal V}$, ${\cal V}'$ denote the iteration trees arising from
the comparison of ${\cal L}$ with ${\cal N}$. We know that $(0,\infty]_V \cap {\cal
D}^{\cal V} = (0,\infty]_{V'} \cap {\cal
D}^{{\cal V}'} = \emptyset$, and
that ${\cal M}_\infty^{\cal V} = {\cal M}_\infty^{{\cal
V}'}$,
because ${\cal L}$ is the $\eta$-core of ${\cal
N}$ and $\sigma_1$ is the core map. We also know
that $\pi^{\cal V}_{0 \infty} {\upharpoonright} \eta = {\rm id}$:
otherwise we can consider the coiteration of $(({\cal N},{\cal L}),\eta)$
with ${\cal N}$ and by
\ref{phalanx-iterations} we'd get a contradiction as in the solidity proof
(cf. \cite[\S 7]{NFS}).

Suppose that $\pi^{{\cal V}'}_{0 \infty} {\upharpoonright} {\tilde \lambda} \not= {\rm id}$.
We'd then have that ${\rm c.p.}(\pi^{{\cal V}'}_{0 \infty}) \leq \mu$, because
by Subclaim 4 above
${\cal N}$ does not have any extenders $E^{\cal N}_\alpha \not= \emptyset$ with
${\rm c.p.}(E^{\cal N}_\alpha) \in (\mu,{\tilde \lambda})$ and $\alpha > {\tilde \lambda}$,
a fact which is inherited by iterations of ${\cal N}$ above ${\tilde \lambda}$.
But then ${\rm c.p.}(\pi^{{\cal V}'}_{0 \infty}) \leq \mu$ and $\rho_1({\cal N}) > \mu$
imply that
$$\rho_1({\cal M}_\infty^{{\cal
V}'}) >
\pi_{0 \infty}^{{\cal V}'}(\mu) \geq
\pi_{0 \infty}^{{\cal V}'}({\rm c.p.}(\pi_{0 \infty}^{{\cal V}'})) >
{\tilde \lambda}.$$ On the other hand, as
$\pi^{\cal V}_{0 \infty} {\upharpoonright} {\tilde \lambda} = {\rm id}$, we must have
$$\rho_1({\cal M}_\infty^{{\cal
V}}) \leq {\tilde \lambda}.$$ This is a contradiction!

Hence we must have that
$\pi^{{\cal V}'}_{0 \infty} {\upharpoonright} {\tilde \lambda} = {\rm id}$, that is, the
coiteration of ${\cal L}$ with ${\cal N}$ is above ${\tilde \lambda}$ on both sides,
$(0,\infty]_V \cap {\cal D}^{\cal V} =
(0,\infty]_{V'} \cap {\cal D}^{{\cal V}`} = \emptyset$,
and ${\cal M}_\infty^{\cal V} = {\cal M}_\infty^{{\cal
V}'}$. This yields (\ref{claim6.2}).

\bigskip
\hfill $\square$ (Claim 6)

\bigskip
Given Claim 6, we can now easily finish the proof of \ref{maximality}. As $\pi$ is
an iteration map, we must have that $\pi$ is continuous at ${\bar \lambda}^{+K^c}$. By
Claim 6, hence, $\pi {\rm " } {\bar \lambda}^{+{\cal M}}$ is cofinal in
$\pi({\bar \lambda}^{+K^c}) = \nu$. But then, if we define the lift up
$${\tilde {\cal M}} = {\rm Ult}_0({\cal M};\pi {\upharpoonright}
{\cal J}^{\cal M}_{{\bar \lambda}^{+{\cal
M}}})$$ as in Case 2 above, then we shall again have that $\nu$ is a cardinal in
${\tilde {\cal M}}$. We may then continue and derive a contradiction exactly as in
Case 2 above.

\bigskip
\hfill $\square$ (\ref{maximality})

\bigskip
It can be shown that \ref{maximality} also holds when the assumption that
$\rho_1({\cal N}) > \mu$ is dropped.

\section{Weak covering for $K^c$.}

Recall that a cardinal $\kappa$ is called countably closed if ${\bar
\kappa}^{\aleph_0} < \kappa$ for all ${\bar \kappa} < \kappa$.
%
We let $B_0$ denote the class of all countably closed singular cardinals $\kappa$
such that $\kappa \in C_0$ and
for all $\mu < \kappa$ do we have that if $\mu$ is
$<\kappa$-strong in $K^c$ then
$\mu$ is $<{\rm OR}$-strong in $K^c$.

It is easy to see (using \ref{C0-is-club})
that $B_0$ is a stationary class. The significance of
$B_0$ is due to the following fact, the ``weak covering 
lemma for $K^c$.''\footnote{We
want to remark that \ref{weak-covering} is not the strongest covering
lemma which is provable for $K^c$. On the other hand, \ref{weak-covering} 
will suffice for the construction of $K$.}

\begin{lemma}\label{weak-covering} ($\lnot \ \zerohandgrenade$)
Let $\kappa \in B_0$. Then
$\kappa^{+K^c} = \kappa^+$.
\end{lemma}

{\sc Proof.} Fix $\kappa \in B_0$, and set
$\lambda = \kappa^{+K^c}$. Let us assume that $\lambda < \kappa^+$. We
shall eventually derive a contradiction.

Fix $\Omega > \kappa^+$, a regular cardinal. Notice that
for $W = (K^c)^{H_\Omega}$ we have that ${\cal J}^{K^c}_\Omega =
{\cal J}^W_\Omega$. We may pick an elementary embedding
$$\pi \ \colon \ N \rightarrow H_\Omega$$ such that
${}^\omega N \subset N$,
${\rm Card}(N) < \kappa$,
$\{ \kappa {\rm , \ } \lambda \} \subset {\rm ran}(\pi)$, and
${\rm ran}(\pi) \cap \lambda$ is cofinal in $\lambda$.
We may and shall assume that $${\rm Card}(N) =
{\rm cf}(\lambda)^{\aleph_0} {\rm , }$$ and that $N \cap ({\rm Card}(N))^+$ is transitive.

Set ${\bar \kappa} = \pi^{-1}(\kappa)$,
${\bar \lambda} = \pi^{-1}(\lambda)$, and
${\bar K} = {\cal J}^{(K^c)^N}_{\bar \lambda}$.
Let $\delta = {\rm c.p.}(\pi) < {\bar \kappa}$. Notice that $\delta =
N \cap ({\rm Card}(N))^+$, ${\rm Card}(\delta) = {\rm Card}(N)$, and $\pi(\delta) = {\rm Card}(N)^+ =
\delta^{+V}$.We can make an immediate observation.

\bigskip
{\bf Claim 1.} $( {\cal P}(\delta) \cap K^c ) \setminus N \not=
\emptyset$.

\bigskip
{\sc Proof.} Suppose otherwise.
Set ${\tilde F} = \pi {\upharpoonright}
{\cal P}(\delta) \cap K^c$.
${\tilde F}$ is countably
complete by ${}^\omega N \subset N$.
It is easy to see that we have a
contradiction with \ref{weak-max}.

\bigskip
\hfill $\square$ (Claim 1)

\bigskip
Let $\eta$ be least such that $( {\cal P}(\eta) \cap K^c ) \setminus N \not=
\emptyset$. Either $\eta = \delta$ and $\delta$ is inaccessible in ${\bar K}$, or else
$\delta = \eta^{+{\bar K}}$ (and then $\pi(\delta) = \eta^{+K^c}$).

Let ${\cal U}$ and ${\cal T}$ denote the iteration trees arising from the comparison
of ${\bar K}$ with $K^c$.
By \ref{Kc-is-universal}, we'll have that ${\cal M}^{\cal T}_\infty
\triangleright {\cal M}^{\cal U}_\infty$,
and ${\cal D}^{\cal U} \cap [0,\infty]_U = \emptyset$.
In fact, the Dodd-Jensen Lemma can easily be used to see that ${\cal T}$ actually
only uses extenders from ${\cal J}_{\Omega}^{K^c}$ and its images.

\[
\xymatrix{
 &               &         & & & \bullet & \kappa^{+K^c} \\
 & {\bar \kappa}^{+{\bar K}} & \bullet & & & \bullet \ar@{-}[u]  & \kappa \\
 & {\bar \kappa} & \bullet \ar@{-}[u] \ar[rrru] & & &         & \\
 {\cal P}_i & . \ar@{.>}@/_/[ddr] &         & & &         & \\
 &               & \bullet \ar@{-}[uu] \ar@{-}[lu] & \mu_i& &         & \\
 & \rho_\omega({\cal P}_i) & \bullet \ar@{-}[u] & & & \bullet \ar@{-}[uuuu] & \pi(\delta) \\
 &   \delta      & \bullet \ar[rrru]^\pi \ar@{-}[u] & & &         & \\
 &   \eta        & \bullet \ar[rrr] \ar@{-}[u] & & & \bullet \ar@{-}[uu] & \\
 &               & {\bar K} \ar@{-}[u]      & & &  {K^c} \ar@{-}[u]       & \\
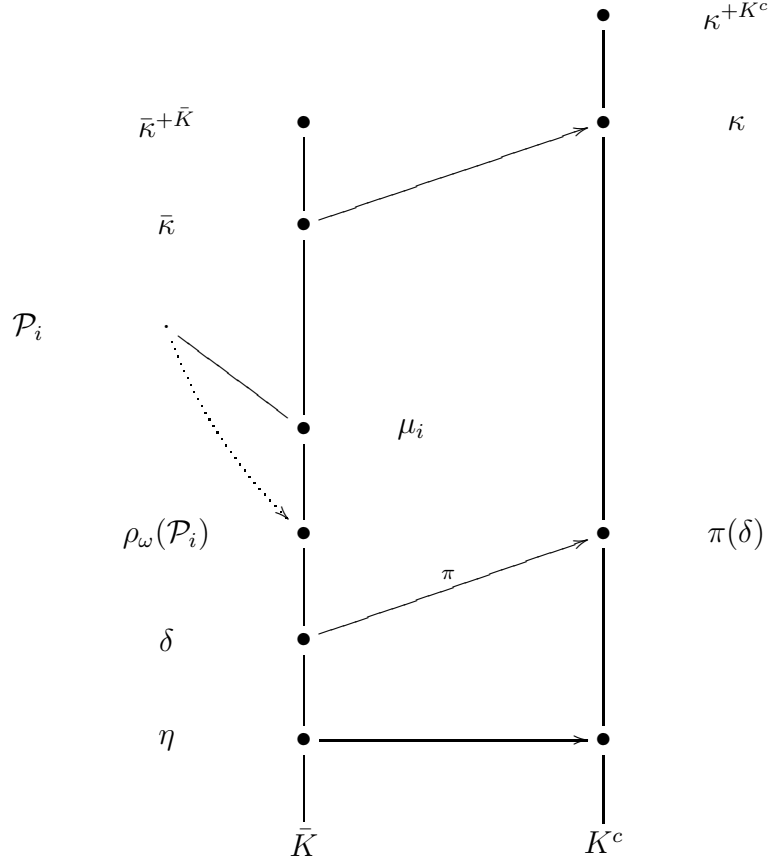
   &  &   & {}\save[]*{\txt{{\sc Figure 3.} The covering argument}} 
\restore & & & \\
}
\]

\bigskip

%
We now aim to prove:

\bigskip
{\bf Claim 2.} ${\cal U}$ is trivial, i.e., ${\cal M}_\infty^{\cal U} = {\bar K}$.

\bigskip
{\sc Proof.}
Suppose that ${\cal U}$
is non-trivial, and let $F = E^{\cal U}_\alpha$ be the first extender used on
$[0,\infty)_U$. Recall that $D^{\cal U} \cap [0,\infty)_U = \emptyset$.
Let $\mu$ be the critical point of $F$. Notice that by \ref{normal-linear}
${\cal U}
{\upharpoonright} \alpha$ only uses extenders with critical point $> \mu^{+{\bar K}}$.

\bigskip
{\bf Subclaim 1.} $F$ is countably complete.

\bigskip
{\sc Proof.} Let $((a_n,X_n) \colon n<\omega)$ be such that $a_n \in F(X_n)$ for all
$n<\omega$. Pick $$\sigma \colon {\bar H} \rightarrow H_{\vartheta}$$ such that
$\vartheta$ is regular and large enough,
${\bar H}$ is transitive, ${\rm Card}({\bar H}) = \aleph_0$,
and $\{ {\cal U} {\upharpoonright} \alpha , F \} \cup \{ a_n ,
X_n \colon n<\omega \} \subset {\rm ran}(\sigma)$. Let ${\bar {\cal U}} =
\sigma^{-1}({\cal U} {\upharpoonright} \alpha)$, and
${\bar {\cal M}} =
\sigma^{-1}({\bar K})$. Notice that ${\bar {\cal U}}$, ${\bar {\cal M}}$, as well as
$\sigma {\upharpoonright} {\bar {\cal M}}$ are all elements of $N$, by ${}^\omega N
\subset N$.

Let us copy ${\bar {\cal U}}$ onto ${\bar K}$, using
$\sigma {\upharpoonright} {\bar {\cal M}}$. The entire copying construction takes
place within $N$, and it gives an iteration tree ${\cal U}'$ on ${\bar K}$
together with copying maps $\sigma_i \colon {\cal M}^{\bar {\cal U}}_i \rightarrow
{\cal M}^{{\cal U}'}_i$ (where $\sigma_0 = \sigma {\upharpoonright} {\bar {\cal M}}$).

As ${\cal U}
{\upharpoonright} \alpha$ only uses extenders with critical point $> \mu^{+{\bar K}}$,
we have that ${\cal U}'$ only uses extenders with critical point
$> \mu^{+{\bar K}}$, too. This implies that $\sigma_0$ and $\sigma_\infty
\colon {\cal M}^{\bar {\cal U}}_\infty \rightarrow
{\cal M}^{{\cal U}'}_\infty$ are such that $\sigma_\infty \circ \sigma_0^{-1}(X_n) =
X_n$
for all $n<\omega$. Moreover, by
\ref{goodness} applied inside $N$ together with ${}^\omega N
\subset N$, we have that $F' = \sigma_\infty \circ \sigma_0^{-1}(F)$ is
(really) countably
complete. Let $a_n' = \sigma_\infty \circ \sigma_0^{-1}(a_n)$ for $n<\omega$, and let
$\tau \colon \bigcup_{n<\omega} \ a_n' \rightarrow \mu$ be order preserving and
such that $\tau {\rm " }
a_n' \in X_n$ for all $n<\omega$. We then have in
$$\tau \circ \sigma_\infty \circ \sigma_0^{-1} {\upharpoonright} \bigcup_{n<\omega}
\ a_n$$ a function as desired.

\bigskip
\hfill $\square$ (Subclaim 1)

\bigskip
Now
notice that \cite[\S 8 Lemma 1]{NFS} 
(a corollary to the solidity proof)
and \ref{phalanx-below-0h} imply that $\mu \geq
\eta$.
Set $\mu_0 = \eta$, and let
$(\mu_i \colon 0<i<\gamma)$ enumerate the successor cardinals of ${\bar K}$
in the half-open interval $(\eta,\mu^{+{\bar K}}]$
(if $\mu = \eta$ then
$\gamma=2$; otherwise, we'll
have $\gamma = \mu+1$).
For every $i < \gamma$ let $\tau(i)$
be the least $\tau$ such that ${\cal J}^{{\cal M}^{\cal T}_\tau}_{\mu_i}
= {\cal J}^{{\bar K}}_{\mu_i}$, and let
${\cal P}_i$ be the
longest initial segment of ${\cal M}^{\cal T}_{\tau(i)}$ which
has the same bounded subsets of
$\mu_i$ as ${\bar K}$ has.

\bigskip
{\bf Subclaim 2.}
${\cal P}_0
= K^c$, and
$i>0 \Rightarrow \rho_\omega({\cal P}_i) < \mu_i$.

\bigskip
{\sc Proof.}
It is enough to verify that ${\cal P}_i$ is a set-sized premouse
for
$i>0$. Let $i>0$, and assume that
${\cal P}_i$ is a weasel. (In particular, ${\cal P}_i = {\cal
M}_{\tau(i)}^{\cal T}$.)
Notice that $\tau(i) > 0$, because $({\cal P}(\eta) \cap K^c) \setminus {\bar K}
\not= \emptyset$.
Let $G = E_j^{\cal T} = E^{{\cal P}_j}_\nu$ be
the first extender used on
$[0,\tau(i))_T$. $G$ is total on $K^c$, and countably complete by \ref{goodness}.
Moreover, of course, $\nu < \mu_i$.
Also, ${\rm c.p.}(G) < \eta$. So
${\rm ran}(G) \subset {\bar K}$ and
${\tilde G} = \pi \circ G$ is well-defined. We have that
$G(\mu)$ is a cardinal of ${\bar
K}$, and hence ${\tilde G}(\mu)$ is a cardinal of $K^c$, using the
elementarity of $\pi$. But as $G$ is countably complete,
${\tilde G}$ is countably complete, too, by the
countable completeness of $\pi$ (i.e., by ${}^\omega N \subset N$).
We thus have a contradiction
with \ref{weak-max}.

\bigskip
\hfill $\square$ (Subclaim 2)

\bigskip
Now set ${\cal P}_\gamma = {\cal
M}_\alpha^{\cal U}$. We shall be interested in the phalanx
$${\vec {\cal P}} =
(({\cal P}_i \colon i<\gamma+1),(\mu_i \colon
i<\gamma)).$$ For $0 < i < \gamma$ let $\mu_i^-$ denote the cardinal predecessor
of $\mu_i$ in
${\cal P}_i$.
Standard arguments give that
for such $i$ is
${\cal P}_i$ sound above $\mu_i^-$, i.e., if $k<\omega$
is such that $\rho_{k+1}({\cal P}_i) < \mu_i \leq \rho_k({\cal P}_i)$ then
${\cal P}_i$ is $k$-sound and $({\cal P}_i)^k$ is the hull of $\mu_i^- \cup \{
p_{({\cal P}_i)^k,1} \}$ generated by $h^k_{{\cal P}_i}$.

\bigskip
{\bf Subclaim 3.} ${\vec {\cal P}}$ is coiterable with $K^c$.

\bigskip
{\sc Proof.} Let us coiterate $K^c$ with
${\vec {\cal P}}$ using $\omega$-maximal trees, getting
iteration trees ${\cal V}$ and ${\cal W}$. Subclaim 3 says that all models of ${\cal
W}$ are transitive, and that we therefore get comparable ${\cal M}^{\cal V}_\infty$ and
${\cal M}^{\cal W}_\infty$.

Let $\infty < {\rm OR}$ denote that ordinal such that either ${\cal M}^{\cal W}_\infty$ is
ill-founded, or else ${\cal M}^{\cal V}_\infty$ and
${\cal M}^{\cal W}_\infty$ are comparable. We aim to show that it is the latter
alternative which holds.
Notice that ${\cal V} {\upharpoonright}
\alpha+1 = {\cal T} {\upharpoonright} \alpha+1$, that ${\cal W} {\upharpoonright}
\alpha$ is trivial, and that $E_\alpha^{\cal W} = F$ which is applied to ${\cal
P}_{\gamma-1}$.

Well, by \ref{phalanx-iterations}, there is $n<\omega$ such
that ${\cal W}$ can be written as $${\cal W}_0 {}^\frown {\cal W}_1 {}^\frown ...
{}^\frown {\cal W}_n$$ where each ${\cal W}_k$ is an iteration of some ${\cal
P}_{i(k)}$
(with $k'>k \Rightarrow i(k')<i(k)$) except for the fact that 
the very first extender
used in ${\cal W}_k$ is equal to $F$ if $k=0$, or else
is taken from the last model of ${\cal W}_{k-1}$ if $k>0$.
Let $F_k$ denote the first extender used in ${\cal W}_k$, and let $\kappa_k =
{\rm c.p.}(F_k)$. Notice that $k'>k
\Rightarrow \kappa_{k'} < \kappa_k$, and that ${\cal W}_k$ is an iteration
which uses only extenders with critical point
$\geq \kappa_k$ (by the rules for iterating a phalanx).

Let us pick an elementary embedding
$$\sigma \colon {\bar H} \rightarrow H_\theta$$ where $\theta$ is regular and large
enough, ${\bar H}$ is countable and transitive, and ${\cal W} \in {\rm ran}(\sigma)$.
Set ${\vec {\cal Q}} = \sigma^{-1}({\vec {\cal P}})$, ${\bar {\cal W}} =
\sigma^{-1}({\cal W})$, ${\bar \kappa}_k = \sigma^{-1}(\kappa_k)$, and ${\bar F} =
\sigma^{-1}(F)$. Then
${\bar {\cal W}}$ is an iteration of ${\vec {\cal Q}}$ which can be written as
$${\bar {\cal W}}_0 {}^\frown {\bar {\cal W}}_1 {}^\frown ...
{}^\frown {\bar {\cal W}}_n$$ where ${\bar {\cal W}}_k = \sigma^{-1}({\cal W}_k)$.

Now as $F = F_0$ is countably complete by Subclaim 1, we may pick some
$\tau \colon (F(\mu) \cap {\rm ran}(\sigma)) \rightarrow \kappa_0 = \mu$
order preserving such that
$a \in F(X) \Rightarrow \tau {\rm " } a \in X$ for appropriate $a$, $X \in
{\rm ran}(\sigma)$.
Set ${\bar \gamma} = \sigma^{-1}(\gamma)$.
Then $${\cal M}^{{\bar
{\cal W}}}_{{\bar \gamma}+1} = {\rm Ult}_k(\sigma^{-1}({\cal P}_{\gamma-1});{\bar F})
{\rm , }$$ where $\rho_{k+1}({\cal P}_{\gamma-1}) \leq
\mu < \rho_{k}({\cal P}_{\gamma-1})$, is
the first model of ${\bar {\cal W}}$ above the list of models ${\vec {\cal Q}}$.
We may define an embedding $${\bar \sigma} \colon {\cal M}^{{\bar
{\cal W}}}_{{\bar \gamma}+1} =
{\rm Ult}_k(\sigma^{-1}({\cal P}_{\gamma-1});{\bar F}) \rightarrow
{\cal P}_{\gamma-1} {\rm , }$$
as
being the extension (via the upward
extensions of embeddings lemma) of the map
$${\bar {\bar \sigma}} \colon ({\cal M}^{{\bar
{\cal W}}}_{{\bar \gamma}+1})^k \rightarrow ({\cal P}_{\gamma-1})^k
{\rm , }$$ where ${\bar {\bar \sigma}}$ is defined by
$$[a,f] \mapsto \sigma(f)(\tau{\rm " }\sigma(a)) {\rm . }$$

Let ${\bar {\cal W}}^*$ be that iteration of the phalanx
$\sigma^{-1}({\vec {\cal P}} {\upharpoonright} \gamma)^\frown {\cal M}^{{\bar
{\cal W}}}_{{\bar \gamma}+1}$
which uses exactly the same extenders (in the same order) as
${\bar {\cal W}}$ does, except for the very first one.
In particular, ${\bar {\cal W}}^*$ and ${\bar {\cal W}}$ have exactly the same models,
except that ${\bar {\cal W}}^*$ misses ${\cal M}_{\bar \gamma}^{{\bar {\cal W}}}$.
Notice that ${\vec {\cal P}} {\upharpoonright} \gamma = ({\cal P}_i \colon i<\gamma)$ is
iterable, as it is generated by an iteration of $K^c$.
It hence suffices to copy the iteration ${\bar {\cal W}}^*$
onto ${\vec {\cal P}} {\upharpoonright} \gamma$, using the maps $\sigma$ and ${\bar
\sigma}$.

However, we have that ${\bar \sigma}$ agrees with $\sigma {\upharpoonright}
\sigma^{-1}({\cal P}_{\gamma-1})$ up to
$\sigma^{-1}(\kappa_1)^+$ (calculated in $\sigma^{-1}({\cal
P}_{\gamma-1})$). Hence the standard copying construction goes
through.

\bigskip
\hfill $\square$ (Subclaim 3)

\bigskip
Let still ${\cal V}$, ${\cal W}$ denote the trees coming from the coiteration of $K^c$
with ${\vec {\cal P}}$.
It is now easy to see that we have to have that $0 = {\rm root}^{\cal W}(\infty)$, i.e.,
that the final model of ${\cal W}$ sits above ${\cal P}_0 = K^c$:
otherwise ${\cal M}^{\cal W}_\infty$ would be non-sound which would give an immediate
contradiction to \ref{Kc-is-universal}.
By applying the
Dodd-Jensen lemma to $\pi_{0 \infty}^{\cal W}$ otherwise, we get
that ${\cal
V}$ has to be simple along its main branch, and ${\cal M}_\infty^{\cal V}
\trianglelefteq {\cal M}^{\cal W}_\infty$. But \ref{Kc-is-universal} implies
that ${\cal W}$ has to be simple along its main branch, and ${\cal M}_\infty^{\cal W}
\trianglelefteq {\cal M}^{\cal V}_\infty$.

Set $Q = {\cal M}_\infty^{\cal V}
= {\cal M}^{\cal W}_\infty$. Let $G$ be the first extender applied on the main
branch of $W$, and
set ${\bar \mu} = {\rm c.p.}(G)$. As the final model of ${\cal W}$ is above
$K^c$, we know that ${\bar \mu} < \eta$.
The Dodd-Jensen lemma tells us that
$\pi^{\cal V}_{0 \infty} \leq \pi^{\cal W}_{0 \infty}$ lexicographically; in
particular, $\pi^{\cal V}_{0 \infty} {\upharpoonright} {\bar \mu} = {\rm id}$.

\bigskip
{\it Case 1.} ${\rm c.p.}(\pi_{0 \infty}^{\cal V}) = {\bar \mu}$.

\bigskip
Let $\beta_0+1$ be least in $(0,\infty]_V$, and let
$\beta_1+1$ be least in $(0,\infty]_W$. We may in this case derive from $${\cal M}^0
= {\cal J}^{{\cal
M}_{\beta_0}^{\cal V}}_{\nu_{\beta_0}} {\rm \ , \ }
{\cal M}^1 = {\cal J}^{{\cal
M}_{\beta_1}^{\cal W}}_{\nu_{\beta_1}}$$
a generalized prebicephalus, call it ${\cal N}$.

\bigskip
{\bf Subclaim 4.}
${\cal N}$ is iterable.

\bigskip
{\sc Proof.} This is shown by
applying \ref{it-of-bicephalus}. Specifically, we plan on letting ${\bar \mu}^{+K^c}$
play the r\^ole of
$\kappa$ in the statement of \ref{it-of-bicephalus}. Set ${\bar \nu} =
{\bar \mu}^{+K^c}$. As ${\bar \nu}$
is a cardinal of $K^c$,
${\cal M}^0 \trianglerighteq {\cal J}^{K^c}_{\bar \nu}$, and
${\cal M}^1 \trianglerighteq {\cal J}^{K^c}_{\bar \nu}$,
by virtue of \ref{goodness} it suffices to
verify the following two statements in order to prove Subclaim 4.

\bigskip
{\bf A}${}_0 \ \ $ if ${\tilde {\cal M}}$ is an iterate of ${\cal M}^0$ above
${\bar \nu}$, and if $H = E^{\tilde {\cal M}}_\nu \not= \emptyset$ is an extender with
${\rm c.p.}(H) < {\bar \nu}$ and $\nu > {\bar \nu}$ then $H$ is countably complete.

\bigskip
{\bf A}${}_1 \ \ $ if ${\tilde {\cal M}}$ is an iterate of ${\cal M}^1$ above
${\bar \nu}$, and if $H = E^{\tilde {\cal M}}_\nu \not= \emptyset$ is an extender with
${\rm c.p.}(H) < {\bar \nu}$ and $\nu > {\bar \nu}$ then $H$ is countably complete.

\bigskip
Now by $\lnot \ \zerohandgrenade$, ${\cal M}^0$ is easily seen to be
an iterate of $K^c$ above
${\bar \nu}$. Hence if ${\tilde {\cal M}}$ and $H$ are as in {\bf A}${}_0$
then ${\tilde {\cal M}}$ is an iterate of $K^c$ above
${\bar \nu}$, too, and $H$ is countably complete by \ref{goodness}. It thus remains to
verify {\bf A}${}_1$.

Fix ${\tilde {\cal M}}$ and $H = E^{\tilde {\cal M}}_\nu$ as in {\bf A}${}_1$, and let
$(a_n,X_n \colon n<\omega)$ be such that $a_n \in H(X_n)$ for every $n<\omega$.
Let ${\cal W}'$ be the iteration of ${\cal M}^1$ which gives ${\tilde {\cal M}}$.
Note that ${\cal M}^1$ in turn is given by
an iteration, call it ${\cal W}^+$, of the phalanx ${\vec {\cal P}}$.
Recall that ${\vec {\cal P}} = ({\cal P}_i \colon i<\gamma+1)$, and that $F$
is the first extender
used in ${\cal W}$. Let ${\vec {\cal P}}^*$ denote the phalanx
$({\cal P}_i \colon i<\gamma)^\frown {\rm Ult}({\cal P}_{\gamma-1};F)$
(with the same exchange
ordinals as ${\vec {\cal P}}$). Then, trivially,
${\cal M}^1$ is given by that iteration, call it
${\cal W}^*$, of ${\vec {\cal P}}^*$ which uses exactly the same extenders (in the same
order) as ${\cal W}^+$ does, except for the very first one, $F$.

Let us pick
an elementary embedding
$$\sigma \colon {\bar H} \rightarrow H_\theta$$ where
$\theta$ is regular and large enough, ${\bar H}$ is countable and transitive,
and $$\{ {\vec {\cal P}}^* , {\cal W}^* , {\cal M}^1 , {\cal W}' , H
\} \cup \{ a_n , X_n
\colon n<\omega \} \subset {\rm ran}(\sigma).$$ Set ${\vec {\cal Q}} =
\sigma^{-1}({\vec {\cal P}}^*)$, ${\bar {\cal W}} = \sigma^{-1}({\cal W}^*)$,
${\bar {\cal M}} = \sigma^{-1}({\cal M}^1)$, and ${\bar {\cal W}}' =
\sigma^{-1}({\cal W}')$. Exactly as in the proof of Subclaim 2 above, we may
first find a map ${\bar \sigma}$
re-embedding
the last model of ${\vec {\cal Q}}$ into ${\cal P}_{\gamma-1}$ (using the countable
completeness of $F$), and we may then use the agreement between ${\bar \sigma}$
and $\sigma$ to copy ${\bar {\cal W}}$ onto ${\vec {\cal P}} {\upharpoonright} \gamma$,
which gives an iteration ${\bar {\cal W}}^c$ of
${\vec {\cal P}} {\upharpoonright} \gamma$ together with the copy maps.
Let $${\bar \sigma}_\infty \colon {\cal M}_\infty^{\bar {\cal W}}
\rightarrow {\cal M}_\infty^{{\bar {\cal W}}^c}$$ be the copy map from the last model
of ${\bar {\cal W}}$ to the last model of ${\bar {\cal W}}^c$.

Now ${\bar {\cal W}}'$ is an iteration of (a truncation of)
${\cal M}_\infty^{\bar {\cal W}}$, and we may hence continue with copying
${\bar {\cal W}}'$ onto ${\cal M}_\infty^{{\bar {\cal W}}^c}$, using
${\bar \sigma}_\infty$, which gives an iteration ${\bar {\cal W}}^{cc}$ of
${\cal M}_\infty^{{\bar {\cal W}}^c}$ together with the copy maps. Let
$$\sigma_\infty \colon {\cal M}_\infty^{{\bar {\cal W}}'}
\rightarrow {\cal M}_\infty^{{\bar {\cal W}}^{cc}}$$
be the copy map from the last model
of ${\bar {\cal W}}'$ to the last model of ${\bar {\cal W}}^{cc}$.

It is now easy to verify that $\sigma$ and $\sigma_\infty$ agree up to
$\sigma^{-1}({\cal J}^{K^c}_{\bar \nu})$.
Moreover, as ${\vec {\cal P}} {\upharpoonright} \gamma$ is given by a normal iteration
of $K^c$, we'll have that ${\bar {\cal W}}^c$ is in fact an iteration of $K^c$ above
${\bar \nu}$. But then
$${\bar {\cal W}}^c {}^\frown {\bar {\cal W}}^{cc}$$ is an iteration
of $K^c$ above
${\bar \nu}$, too. By \ref{goodness}, we thus know that $\sigma_\infty \circ
\sigma^{-1}(H)$ is countably complete. We may hence pick
$\tau \colon \bigcup_{n<\omega} \
\sigma_\infty \circ
\sigma^{-1}(a_n) \rightarrow {\rm c.p.}(H)$ such that $\tau {\rm " } \sigma_\infty \circ
\sigma^{-1}(a_n) \in X_n$ for every $n<\omega$. We have found a function as desired!

\bigskip
\hfill $\square$ (Subclaim 4)

\bigskip
But now by \ref{bicephali-trivialize} we have that Subclaim 4 implies
$$E^{\cal V}_{\beta_0}
= E^{\cal W}_{\beta_1}.$$ This is a contradiction!

\bigskip
{\it Case 2.} $\pi_{0 \infty}^{\cal V} {\upharpoonright} {\bar \mu} +1 = {\rm id}$.

\bigskip
In this case, we get a contradiction with \ref{maximality}. That \ref{maximality} is
applicable here follows from arguments exactly as in Case 1 above. We leave
it to the reader to chase through the obvious details.

\bigskip
\hfill $\square$ (Claim 2)

\bigskip
{\bf Claim 3.} ${\cal T}$ only uses extenders with critical point greater than or equal to $\eta$.

\bigskip
{\sc Proof.} This follows from the argument which was presented in the proof of
Subclaim 2 above.

\bigskip
\hfill $\square$ (Claim 3)

\bigskip
The combination of Claims 1, 2, and 3 now give an initial segment, ${\cal M}$, of ${\cal
M}_\infty^{\cal T}$ with ${\cal M} \triangleright {\bar K}$, ${\bar \lambda}$ is
a cardinal in ${\cal M}$, ${\cal M}$ is ${\bar \kappa}$-sound, and
$\rho_\omega({\cal M}) \leq {\bar \kappa}$.
Let $n<\omega$ be such that $\rho_{n+1}({\cal M}) \leq {\bar \kappa} < \rho_n({\cal
M})$. Let $${\tilde {\cal M}} =
{\rm Ult}_n({\cal M};\pi {\upharpoonright} {\bar K}).$$
Then either ${\tilde {\cal M}}$ is an $n$-iterable
premouse
with ${\tilde {\cal M}} \triangleright {\cal J}_\lambda^{K^c}$,
$\rho_{n+1}({\tilde {\cal M}}) \leq \kappa$, and ${\tilde {\cal M}}$ is sound above
$\kappa$,
or else $n = 0$,
${\cal M}$ has a top extender, ${\bar F}$,
$\pi {\rm " } {\rm c.p.}({\bar F})^{+{\cal M}}$ is
not cofinal in $\pi({\rm c.p.}({\bar F})^{+{\cal M}})$,
and ${\tilde {\cal M}}$ is a protomouse
(see \cite[\S 2.3]{MiSchSt}). In the latter case,
let ${\bar F}'$ be the top extender of ${\tilde {\cal M}}$, and
let ${\cal N}
= {\cal J}_\eta^{K^c}$ be
the longest initial segment of $K^c$ which has only subsets of ${\rm c.p.}({\bar F}')$
which are
measured by ${\bar F}'$. Let $${\tilde {\cal N}} = {\rm Ult}_m({\cal N};{\bar F}') {\rm , }$$
where $\rho_{m+1}({\cal N}) \leq {\rm c.p.}({\bar F}') < \rho_m({\cal N})$.
Using the countable completeness of $\pi$, by standard arguments we'll have
that ${\tilde {\cal N}}$ is an $m$-iterable premouse with ${\tilde {\cal N}}
\triangleright {\cal J}_\lambda^{K^c}$,
$\rho_{\omega}({\tilde {\cal N}}) \leq \kappa$, and ${\tilde {\cal N}}$ is sound above
$\kappa$.

But now, finally, \ref{on-the-sequence} gives a contradiction!

\bigskip
\hfill $\square$ (\ref{weak-covering})

\section{Beavers and the existence of K.}
\setcounter{equation}{0}

In this section we shall isolate $K$, the true core model below $\zerohandgrenade$.
We closely follow \cite[\S 5]{CMIP}; however, 
of course, some revisions are necessary, as
\cite{CMIP} works in the theory ``${\sf ZFC} + \Omega$ is measurable.''

We first need a concept of thick classes, which is originally due to
\cite{mitchell2}. We commence with a simple observation (which is standard).

\begin{lemma}\label{covering-universality} $( \ \lnot \ \zerohandgrenade \ )$
Let $W$ be an iterable weasel. Let $S$ be a stationary class such that for all $\beta
\in S$ we have that $\beta$ is a strong limit cardinal,
${\rm cf}^W(\beta)$ is not measurable in $W$, and
$\beta^{+W} = \beta^+$. Then $W$ is universal.
\end{lemma}

{\sc Proof.} Suppose not.
Then there is a (set- or class-sized) premouse ${\cal M}$ such
that if ${\cal T}$, ${\cal U}$ are the iteration trees arising from the comparison of
${\cal M}$ with $W$ we have that: ${\rm lh}({\cal T}) = {\rm lh}({\cal U}) = {\rm OR}+1$, there is no
drop on $[0,{\rm OR}]_U$, $\pi_{0 \infty}^{\cal U} {\rm " } {\rm OR} \subset {\rm OR}$, and if $\alpha$
is largest in $\{ 0 \} \cup ({\cal D}^{\cal T} \cap (0,\infty)_T)$ then $\pi_{\alpha
\infty}^{\cal T} {\rm " } {\rm OR} \not\subset {\rm OR}$. (We have $\infty = {\rm OR}$ here.)
There are club many $\beta \in [0,\infty)_U$ such that
$\pi_{0 \beta}^{\cal U} {\rm " } \beta \subset \beta$.
Also, there are $\gamma \in [\alpha,\infty)_T$
and $\kappa \in {\cal M}^{\cal T}_\gamma$ such that for
club many $\beta \in [\gamma,\infty)_T$ do we have that
$\pi^{\cal T}_{\gamma
\beta}(\kappa) = \beta$ and $\beta$ is the critical point of $\pi^{\cal T}_{\beta
\infty}$.

Now pick $\beta \in S \setminus ({\cal M}_\gamma^{\cal T} \cap {\rm OR})$
with $\pi_{0 \beta}^{\cal U} {\rm " } \beta \subset \beta$,
$\pi^{\cal T}_{\gamma
\beta}(\kappa) = \beta$, and ${\rm c.p.}(\pi^{\cal T}_{\beta \infty}) = \beta$.
As $\beta$ is a strong limit cardinal and
${\rm cf}^W(\beta)$ is not measurable in $W$, we'll have that $\pi_{0 \infty}^{\cal
U}(\beta) = \beta$, which implies that
$\beta^{+{\cal M}^{\cal U}_\infty} = \beta^+$.
On the other hand, we'll have that $\pi_{\gamma \beta}^{\cal T} {\rm " }
\kappa^{+{\cal M}^{\cal T}_\gamma}$ is cofinal in $\beta^{+{\cal M}^{\cal T}_\beta}
= \beta^{+{\cal M}^{\cal T}_\infty}$,
so that $\beta^{+{\cal M}^{\cal T}_\infty} < \beta^+$. Hence
$\beta^{+{\cal M}^{\cal T}_\infty} < \beta^+ = \beta^{+{\cal M}^{\cal U}_\infty}$.
This is a contradiction!

\bigskip
\hfill $\square$ (\ref{covering-universality})

\begin{defn}\label{defn-thick}
Let $W$ be a weasel, let $S \subset {\rm OR}$, and let $\tau \in {\rm OR}$.
A class $\Gamma \subset {\rm OR}$ is called
$S,\tau$-thick in $W$ provided the following clauses hold.

\noindent (i) $S$ is stationary in ${\rm OR}$,

\noindent (ii) for all but nonstationary many $\beta \in S$ do we have that

(a) $\beta^{+W} = \beta^+$,

(b) $\beta$ is a strong limit cardinal,

(c) ${\rm cf}^W(\xi)$ is not measurable in $W$ as witnessed by ${\vec E}^W$ (i.e., there
is no $E_\nu^W \not= \emptyset$ with ${\rm c.p.}(E_\nu^W) = {\rm cf}^W(\xi)$ and $E_\nu^W$
is total on $W$) for all $\xi \in
\Gamma \cap [\beta,\beta^+)$, and

(d) $\Gamma \cap (\beta,\beta^+)$ is unbounded in $\beta^+$, and
for all regular $\theta > \tau$ we have that
$\Gamma \cap (\beta,\beta^+)$ is $\theta$-closed,
and $\beta \in \Gamma$.

A class $\Gamma \subset {\rm OR}$ is called
$S$-thick in $W$ if there is some $\tau$ such that $\Gamma$ is
$S,\tau$-thick in $W$.
\end{defn}

By the proof of \ref{covering-universality}, if $\Gamma \subset {\rm OR}$ is
$S$-thick in $W$ for some $S \subset {\rm OR}$ then $W$ is universal.

The following four lemmas are easy to prove 
(cf. \cite[Lemmas 3.9 to 3.11]{CMIP}).

\begin{lemma}\label{thick-lemma-1}
Let $W$ be a weasel, let $S \subset {\rm OR}$ be a class, and let $(\Gamma_i \colon
i<\theta)$ be such that $\Gamma_i$ is $S$-thick in $W$ for all $i<\theta$. Then
$\bigcap_{i<\theta} \ \Gamma_i$ is $S$-thick in $W$.
\end{lemma}

\begin{lemma}\label{thick-lemma-2}
Let ${\bar W}$ and $W$ be weasels, and let $S \subset {\rm OR}$ be a class.
Let $\pi \colon {\bar W} \rightarrow W$ be an elementary embedding such that there is
some $\Gamma \subset {\rm ran}(\pi)$ which is $S$-thick in $W$.
Then $\Gamma \cap \{ \xi \colon \pi(\xi) = \xi \}$
is $S$-thick in both ${\bar W}$ and $W$.
\end{lemma}

\begin{lemma}\label{thick-lemma-2a}
Let $W$ be a weasel, let $S \subset {\rm OR}$ be a class, and let $\Gamma$ be $S$-thick in
$W$. Let $F$ be a total extender on $W$ such that ${\rm Ult}(W;F)$ is transitive.
Let $i_F$ denote the ultrapower map. Then
$\Gamma \cap \{ \xi \colon i_F(\xi) = \xi \}$ is $S$-thick in both $W$ and ${\rm Ult}(W;F)$.
\end{lemma}

\begin{lemma}\label{thick-lemma-3}
Let ${\cal T}$ be an iteration tree on the weasel $W$, let $S \subset {\rm OR}$ be a class,
and let $\Gamma$ be $S$-thick in $W$. Let $\alpha \leq {\rm lh}({\cal T})$ (possibly,
$\alpha = {\rm OR}$) be such that ${\cal D}^{\cal T} \cap (0,\alpha]_T = \emptyset$ and
$\pi^{\cal T}_{0 \alpha} {\rm " } {\rm OR} \subset {\rm OR}$. Then
$\Gamma \cap \{ \xi \colon \pi^{\cal T}_{0 \alpha}(\xi) = \xi \}$
is $S$-thick in both $W$ and ${\cal M}_\alpha^{\cal T}$.
\end{lemma}

In order to isolate $K$ we now have to work with some $K^c_\Gamma$ for a small
$\Gamma$ rather than with $K^c$, just because there need not be $S$ and
$\Gamma^*$ such that $\Gamma^*$ is $S$-thick in $K^c$. Let us
consider $K^c_{\{\omega_1\}}$.
We get versions of everything we proved in sections 4 through 7 about $K^c$ also for
$K^c_{\{\omega_1\}}$. (When quoting a result from one of these sections in what
follows, what
we'll in fact intend to quote is the corresponding version for $K^c_{\{\omega_1\}}$.)

In particular, there is a stationary class $B_0^{ \{ \omega_1 \} }$ 
(defined in a similar fashion as
the old $B_0$ was defined for $K^c$) with
$\kappa^{+K^c_{\{\omega_1\}}} = \kappa^+$ for all 
$\kappa \in B_0^{ \{ \omega_1 \} }$. 
With this meaning
of $B_0^{ \{ \omega_1 \} }$ in mind we now 
let $A_0$ denote the class of all elements of $B_0^{ \{ \omega_1 \} }$ 
which are strong limit
cardinals of cofinality $> \omega_1$.
I.e., $A_0$ consists of all singular strong limit cardinals $\kappa$ 
of cofinality $> \omega_1$ such that
${\cal J}^{K^c_{ \{ \omega_1 \} }}_\kappa$ is universal for coiterable 
premice of height $< \kappa$ and such that
if $\mu$ is $< \kappa$-strong in $K^c_{ \{ \omega_1 \} }$ then $\mu$ is 
$< {\rm OR}$-strong in $K^c_{ \{ \omega_1 \} }$.
Of course, $A_0$ is a stationary subclass of $B_0^{ \{ \omega_1 \} }$.

The following is now immediate from \ref{weak-covering}.

\begin{lemma}\label{thick-in-Kc} $( \ \lnot \ \zerohandgrenade \ )$
Let $\Gamma$ be the class of all
ordinals of cofinality $> \omega_1$. Then $\Gamma$ is $A_0$-thick in ${K^c_{ \{ \omega_1 \} }}$.
\end{lemma}

{\sc Proof.} As for (ii) (c) in \ref{defn-thick}, notice that if ${\rm cf}^{{K^c_{ \{ \omega_1 \} }}}(\xi)$ were
measurable in ${K^c_{ \{ \omega_1 \} }}$ as witnessed by ${\vec E}^{{K^c_{ \{ \omega_1 \} }}}$ then we would have to have
${\rm cf}^V(\xi) = {\rm cf}^V({\rm cf}^{{K^c_{ \{ \omega_1 \} }}}(\xi)) = \omega_1$,
because if $E_\nu^{{K^c_{ \{ \omega_1 \} }}} \not= \emptyset$ is total on ${K^c_{ \{ \omega_1 \} }}$ then its critical point
has cofinality $= \omega_1$. Contradiction!

\bigskip
\hfill $\square$ (\ref{thick-in-Kc})
\bigskip

Let $W$ be a weasel, and let ${\vec \alpha} \in {\rm OR}$. For the following purposes
we shall denote by $x = \tau^W[{\vec \alpha}]$ the fact that there is some formula
$\Phi$ such that $x$ is the unique ${\bar x}$ with $W \models
\Phi({\bar x},{\vec \alpha})$,
and $\tau$ is the term given by $\Phi$.
For $\Gamma \subset {\rm OR}$ we shall let $x \in H^W(\Gamma)$ mean that
$x = \tau^W[{\vec \alpha}]$ for some term $\tau$ and ${\vec \alpha} \in \Gamma$.

\begin{defn}\label{defn-hull-prop}
Let $W$ be a weasel, and let $S \subset {\rm OR}$.
Let $\alpha \in {\rm OR}$.
We say that $W$ has the $S$-hull property at $\alpha$ just in case that for all
$S$-thick $\Gamma$ do we have
that ${\cal P}(\alpha) \cap W$ is a subset of the transitive collapse of
$H^W(\alpha \cup \Gamma)$.
\end{defn}

\begin{lemma}\label{hull-prop} $( \ \lnot \ \zerohandgrenade \ )$
There is some club $C \subset {\rm OR}$ such that ${K^c_{ \{ \omega_1 \} }}$ has the $A_0$-hull
property at any $\kappa \in C$.
\end{lemma}

{\sc Proof.} Let $C$ be the class of all limit cardinals $\kappa$ (of $V$)
such that for any $\mu < \kappa$ we have that if $\mu$ is
$< \kappa$-strong in ${K^c_{ \{ \omega_1 \} }}$ then
$\mu$ is $< {\rm OR}$-strong in ${K^c_{ \{ \omega_1 \} }}$. Let $\kappa \in C$. We aim to show
that ${K^c_{ \{ \omega_1 \} }}$ has the $A_0$-hull property at $\kappa$.

Let $\Gamma$ be $A_0$-thick in ${K^c_{ \{ \omega_1 \} }}$, and let
$\sigma \colon W \rightarrow {K^c_{ \{ \omega_1 \} }}$ be elementary with $W$ transitive,
and ${\rm ran}(\sigma) = H^{{K^c_{ \{ \omega_1 \} }}}(\kappa \cup \Gamma)$. Notice that $W$ is universal.
It suffices to
prove that the coiteration of $W$, ${K^c_{ \{ \omega_1 \} }}$ is above $\kappa$ on both sides.
We first show the following easy

\bigskip
{\bf Claim.} Let ${\cal M}$ be an iterate of $W$ above $\kappa$, and let $F = E^{\cal
M}_\nu \not= \emptyset$ be such that $\nu > \kappa$ and ${\rm c.p.}(F) <
\kappa$. Then $F$ is countably complete.

\bigskip
{\sc Proof.} Let ${\cal M} = {\cal M}^{\cal T}_\infty$ where ${\cal T}$ is the
iteration tree on $W$ giving ${\cal M}$. Using
$\sigma \colon W \rightarrow {K^c_{ \{ \omega_1 \} }}$
we may copy ${\cal T}$ onto ${K^c_{ \{ \omega_1 \} }}$, getting an iteration tree ${\cal U}$
on ${K^c_{ \{ \omega_1 \} }}$ together with a last copy
map $\sigma_\infty \colon {\cal M} \rightarrow {\cal M}^{\cal U}_\infty$. As ${\cal
T}$ is above $\kappa$, ${\cal U}$ is above $\kappa$, too. Hence $\sigma_\infty(F)$ is
countably complete by \ref{goodness}. This clearly implies that $F$ is countably
complete, because $\sigma_\infty {\upharpoonright} {\cal P}({\rm c.p.}(F)) \cap W = {\rm id}$.

\bigskip
\hfill $\square$ (Claim)

\bigskip
Now let ${\cal T}$, ${\cal U}$ denote the iteration trees arising from the comparison
of $W$ with ${K^c_{ \{ \omega_1 \} }}$. If $\pi^{\cal T}_{0 \infty} {\upharpoonright} \kappa \not= {\rm id}$ or
$\pi^{\cal U}_{0 \infty} {\upharpoonright} \kappa \not= {\rm id}$ then exactly as in the proof
of \ref{on-the-sequence}
Claim 2 Subclaim 1
we'll get that ${\rm c.p.}(\pi^{\cal T}_{0 \infty}) =
{\rm c.p.}(\pi^{\cal U}_{0 \infty}) < \kappa$. (Notice that any $\mu < \kappa$ is
$<{\rm OR}$-strong in
$W$ if and only if it is $<{\rm OR}$-strong in ${K^c_{ \{ \omega_1 \} }}$!)
Hence if $\alpha+1$ is least in $(0,\infty]_T$ and
$\beta+1$ is least in $(0,\infty]_U$, and if $E_\alpha^{\cal T} = E^{{\cal M}^{\cal
T}_\alpha}_{\nu_0}$ and $E_\beta^{\cal U} = E^{{\cal M}^{\cal
U}_\beta}_{\nu_1}$, then we may derive from
$${\cal J}^{{\cal M}^{\cal T}_\alpha}_{\nu_0} {\rm \ , \ }
{\cal J}^{{\cal M}^{\cal U}_\beta}_{\nu_1}$$ a generalized prebicephalus, call it ${\cal N}$.

But using the Claim above, ${\cal N}$ is iterable by \ref{it-of-bicephalus}
(where we let the current $\kappa$ play the r\^ole of the $\kappa$ in the statement
of \ref{it-of-bicephalus}).
This gives a contradiction as
in the proof of \ref{on-the-sequence}.

\bigskip
\hfill $\square$ (\ref{hull-prop})

\begin{lemma}\label{hull-prop-at-strongs}
Let $W$ be a universal weasel, and let $\Gamma$ be $S$-thick in $W$. Then $W$ has the
hull property at all $\kappa$ such that $\kappa$ is $< {\rm OR}$-strong in $W$.
\end{lemma}

{\sc Proof.} This is shown by induction on $\kappa$.
Let $W$ be a universal weasel, and let $\Gamma$ be $S$-thick in $W$.
Let $$\sigma \colon {\bar W} \cong H^W(\kappa \cup \Gamma) \prec W.$$
Let ${\cal U}$, ${\cal T}$ denote the iteration trees arising from the comparison of
$W$ with ${\bar W}$.
It suffices to
prove that $\pi_{0 \infty}^{\cal U} {\upharpoonright} \kappa =
\pi_{0 \infty}^{\cal T} {\upharpoonright} \kappa = {\rm id}$. Suppose not. Then $\mu =
{\rm min} \{ {\rm c.p.}(\pi_{0 \infty}^{\cal U}),$ ${\rm c.p.}(\pi_{0 \infty}^{\cal T}) \}$ is $<
\kappa$-strong in
$W$, and hence $< {\rm OR}$-strong in both $W$ and
${\bar W}$. By an argument as in the proof of \ref{on-the-sequence} Claim 2 Subclaim 1
we then get that in fact $\mu =
{\rm c.p.}(\pi_{0 \infty}^{\cal U}) = {\rm c.p.}(\pi_{0 \infty}^{\cal T})$. Let $\alpha+1$ be
least in $(0,\infty]_U$, and let
$\beta+1$ be
least in $(0,\infty]_T$.
By inductive hypothesis both $W$ and ${\bar W}$
have the hull-property at $\mu$. Set $Q = {\cal M}_\infty^{\cal
U} = {\cal M}_\infty^{\cal
T}$.

By \ref{thick-lemma-1}, \ref{thick-lemma-2}, and \ref{thick-lemma-3},
$\Gamma' = \{ \xi \in \Gamma
\colon \sigma(\xi) = \pi^{\cal U}_{0 \infty}(\xi) = \pi^{\cal T}_{0 \infty}(\xi) = \xi
\}$ is $A_0$-thick in $W$, ${\bar W}$,
and $Q$. Let $F = E_\alpha^{\cal U}$, and $F^* = E_\beta^{\cal
T}$. Let $a \in [F(\mu) \cap F^*(\mu)]^{< \omega}$, and $X \in {\cal
P}([\mu]^{{\rm Card}(a)}) \cap W = {\cal
P}([\mu]^{{\rm Card}(a)}) \cap {\bar W}$. Let $X =
\tau^{W}[{\vec \xi}] \cap [\mu]^{{\rm Card}(a)}$
for ${\vec \xi} \in \Gamma'$. Then $X = \tau^{{\bar W}}[{\vec \xi}]
\cap [\mu]^{{\rm Card}(a)}$, too. Moreover, $a \in F(X)$ if and only if $a \in
\pi^{\cal U}_{0 \infty}(X) = \tau^Q[{\vec \xi}] \cap
\pi^{\cal U}_{0 \infty}([\mu]^{{\rm Card}(a)})$ if and only if $a \in
\tau^Q[{\vec \xi}] \cap
\pi^{\cal T}_{0 \infty}([\mu]^{{\rm Card}(a)}) = \pi^{\cal T}_{0 \infty}(X)$ if and only if
$a \in F^*(X)$. Hence $F$ and $F^*$ are compatible. Contradiction!

\bigskip
\hfill $\square$ (\ref{hull-prop-at-strongs})

\bigskip
Let $W$ be a weasel, and let $S \subset {\rm OR}$. We shall let $x \in {\rm Def}(W,S)$
mean that whenever $\Gamma$ is $S$-thick in $W$ then $x \in H^W(\Gamma)$.

If $W$, $W'$ are both (universal) weasels, and $S \subset {\rm OR}$,
such that there is some $\Gamma \subset {\rm OR}$
which is $S$-thick in $W$ as well as in $W'$ then it is easy to verify that
${\rm Def}(W,S)$ and ${\rm Def}(W',S)$ have the same transitive collapse (cf. \cite[Corollary 5.7]{CMIP}. Hence for $S$ a stationary class
we denote by $K(S)$ the transitive collapse of ${\rm Def}(W,S)$ for
any (some) universal $W$ such that there is some $\Gamma \subset {\rm OR}$
which is $S$-thick in $W$.

The proof of the following lemma \ref{K-weasel}
needs some care in order to not to go beyond the bounds of
predicative class theory. We shall need the following auxiliary concept,
which was introduced in \cite[Definition 2.7]{kairalf} 
(and generalizes a concept of
\cite{jensen}; see also \cite{schdl}).

\begin{defn}\label{defn-beaver}
Let ${\cal B} = (J_\alpha[{\vec E}];\in,{\vec E},F)$ be a premouse with top extender
$F \not= \emptyset$. Set $\kappa = {\rm c.p.}(F)$.
Then ${\cal B}$ is called a beaver provided there is a universal
weasel $W$ with the following properties.

(a) ${\cal J}^W_\lambda = {\cal J}^{\cal B}_\lambda$ where $\lambda =
\kappa^{+W} = \kappa^{+{\cal B}}$,

(b) $W$ has the definability property at all $\mu < \kappa$ such that $\mu$
is $< \kappa$-strong in $W$, and

(c) ${\rm Ult}_0(W;F)$ is $0$-iterable.
\end{defn}

It is straightforward to check that the proof of 
\cite[Lemma 2.9]{kairalf} works
below $\zerohandgrenade$; this observation establishes:

\begin{lemma}\label{double-beavers} $( \ \lnot \ \zerohandgrenade \ )$
Suppose that ${\cal B} = (J_\alpha[{\vec E}];\in,{\vec E},F)$ and
${\cal B}' = (J_\alpha[{\vec E}];\in,{\vec E},F')$ are beavers. Then $F = F'$.
\end{lemma}
%
%
\begin{lemma}\label{K-weasel}
$( \ \lnot \ \zerohandgrenade \ )$
$K(A_0)$ is a weasel.
\end{lemma}

{\sc Proof.} We aim to show that ${\rm Def}({K^c_{ \{ \omega_1 \} }},A_0)$
is unbounded in ${\rm OR}$. Let $\Gamma$ be the class of all ordinals
of cofinality $> \omega_1$. By \ref{thick-in-Kc}, $\Gamma$ is
$A_0$-thick in ${K^c_{ \{ \omega_1 \} }}$.

\bigskip
{\bf Claim 1.} There is a sequence
$(\Gamma_\kappa \
\colon \ \kappa \in {\rm OR})$ of classes such that

(A) $\Gamma_0 = \Gamma$,

(B) $\Gamma_\kappa \supset \Gamma_{\kappa'}$ for
$\kappa \leq \kappa'$,

(C) if there is some ${\bar \Gamma} \subset \Gamma_\kappa$
which is $A_0$-thick in ${K^c_{ \{ \omega_1 \} }}$ and
such
that $\kappa \notin H^{{K^c_{ \{ \omega_1 \} }}}({\bar \Gamma})$ then
$\kappa \notin H^{{K^c_{ \{ \omega_1 \} }}}(\Gamma_{\kappa+1})$, for all ordinals $\kappa$, and

(D) $\Gamma_\lambda = \bigcap_{\kappa < \lambda} \ \Gamma_\kappa$
for all
limit ordinals $\lambda$.

\bigskip
{\sc Proof.} Of course, Claim 1 is supposed to say that the class $\{
(\kappa,x) \ \colon \ x \in \Gamma_\kappa
\}$ exists. The point of Claim 1 is that we have to be able to {\it define} such a
class.
For this purpose, we need a refinement of the argument in 
\cite[Appendix II]{MOZ}. Let us first indicate how we aim to
choose $\Gamma_{\kappa+1}$, given $\Gamma_\kappa$.

Suppose that
$\Xi$ is $A_0$-thick in ${K^c_{ \{ \omega_1 \} }}$, and $\beta < \alpha
\wedge \beta \notin {\rm Def}({K^c_{ \{ \omega_1 \} }},A_0) \Rightarrow \beta \notin H^{{K^c_{ \{ \omega_1 \} }}}(\Xi)$.
We want to find some canonical
$\Xi'$ which is $A_0$-thick in ${K^c_{ \{ \omega_1 \} }}$, and
$\alpha \notin {\rm Def}({K^c_{ \{ \omega_1 \} }},A_0) \Rightarrow \alpha \notin H^{{K^c_{ \{ \omega_1 \} }}}(\Xi')$.
Let us suppose without loss of generality that $\alpha \in H^{{K^c_{ \{ \omega_1 \} }}}(\Xi) \setminus {\rm Def}({K^c_{ \{ \omega_1 \} }},A_0)$.
Let $\Xi_0
\subset \Xi$ be $A_0$-thick in ${K^c_{ \{ \omega_1 \} }}$ such that $\alpha \notin H^{{K^c_{ \{ \omega_1 \} }}}(\Xi_0)$.
Let
$$\sigma_0 \colon {\bar W}_0 \cong H^{{K^c_{ \{ \omega_1 \} }}}(\Xi_0) \prec {K^c_{ \{ \omega_1 \} }} {\rm , \ and }$$
$$\sigma \colon {\bar W} \cong H^{{K^c_{ \{ \omega_1 \} }}}(\Xi) \prec {K^c_{ \{ \omega_1 \} }} {\rm . }$$
Then $\sigma^{-1}(\alpha)$ is the critical point of $\sigma^{-1} \circ \sigma_0$, and
both ${\bar W}_0$ and ${\bar W}$ have the definability property at
all $\mu < \sigma^{-1}(\alpha)$. Set $\kappa = \sigma^{-1}(\alpha)$.
The proof of \cite[Lemma 1.3]{kairalf} then
works below
$\zerohandgrenade$ and shows
that $\kappa^{+{\bar W}_0} = \kappa^{+{\bar W}}$. Let us write
$\lambda = \kappa^{+{\bar W}}$. Set
$F = \sigma^{-1} \circ \sigma_0 {\upharpoonright} {\cal P}(\kappa) \cap {\bar W}_0$,
and ${\tilde \lambda} = {\rm sup} \ \sigma^{-1} \circ \sigma_0 {\rm " } \lambda$.
Now consider $${\cal B} =
({\cal J}^{\bar W}_{\tilde \lambda},F).$$ Let us suppose without loss of generality that
${\cal B}$ is a
premouse (if not, then the initial segment condition fails for $F$, and we may
replace ${\cal B}$ by a premouse obtained as in the proof of \ref{weak-max}).
Moreover, as ${\rm Ult}_0({\bar W}_0;F)$ can be embedded into ${\bar W}$, we
immediately get that in fact
${\cal B}$
is a beaver.
On the other hand, by \ref{double-beavers},
there is at most one ${\bar F}$ such
that $({\cal J}^{\bar W}_{\tilde \lambda},{\bar F})$ is a beaver.

By \cite[Lemma 1.2]{kairalf} 
we now know that ${\rm Ult}_0({\bar W};F)$ is $0$-iterable,
too. Let $$\pi \colon {\bar W} \rightarrow_F {\tilde W}.$$ By coiterating ${\bar W}$,
${\tilde W}$ we get a common coiterate $Q$ together
with iteration maps $\pi_{{\bar W},Q}$ and $\pi_{{\tilde W},Q}$.
Set $$\Xi' = \{ \xi \in \Xi \ \colon \ \sigma(\xi) = \pi_{{\bar W},Q}(\xi) =
\pi_{{\tilde W},Q} \circ \pi(\xi) = \xi \}.$$ Then
$\Xi'$ is $A_0$-thick in ${K^c_{ \{ \omega_1 \} }}$
by \ref{thick-lemma-1}, \ref{thick-lemma-2}, \ref{thick-lemma-2a},
and
\ref{thick-lemma-3}.
Notice that $x \in H^{{K^c_{ \{ \omega_1 \} }}}(\Xi')$ of course implies that
$\pi_{{\bar W},Q} \circ \sigma^{-1}(x) \in {\rm ran}(\pi_{{\tilde W},Q} \circ \pi)$.

\bigskip
{\bf Subclaim 1.} $\alpha \notin H^{{K^c_{ \{ \omega_1 \} }}}(\Xi')$.

\bigskip
{\sc Proof.} Suppose otherwise, and let $\alpha = \tau^{{K^c_{ \{ \omega_1 \} }}}[{\vec \xi}]$ where
${\vec \xi} \in \Xi'$. Then $\kappa = \tau^{\bar W}[{\vec \xi}]$. Notice that both
${\bar W}$ and ${\tilde W}$ have the definability property at all $\mu < \kappa$. By
\cite[Corollary 1.5]{kairalf}, hence, 
$\pi_{{\bar W},Q} {\upharpoonright} \kappa =
\pi_{{\tilde W},Q} {\upharpoonright} \kappa = {\rm id}$.

\bigskip
{\it Case 1.} $\pi_{{\bar W},Q} {\upharpoonright} \kappa+1 = {\rm id}$.

\bigskip
In this case we get that $\kappa = \tau^Q[{\vec \xi}] \in {\rm ran}(\pi_{{\tilde W},Q} \circ
\pi)$. However, $\kappa$ is the critical point of $\pi$, and
$\pi_{{\tilde W},Q} {\upharpoonright} \kappa = {\rm id}$. Contradiction!

\bigskip
{\it Case 2.} $\kappa$ is the critical point of $\pi_{{\bar W},Q}$.

\bigskip
Because ${\cal B}$ is a premouse, we have that ${\cal J}_{\tilde \lambda}^{\bar W} =
{\cal J}_{\tilde \lambda}^{\tilde W}$. By $\lnot \ \zerohandgrenade$
(i.e., by \ref{not-strong}), and because
$\pi_{{\tilde W},Q} {\upharpoonright} \kappa = {\rm id}$, we hence get that in fact
$\pi_{{\tilde W},Q} {\upharpoonright} F(\kappa) = {\rm id}$.
By \cite[Lemma 1.3]{kairalf} 
both ${\bar W}$ and ${\tilde W}$ have the hull property at
$\kappa$.
Let $F^*$ be the first extender used on the main branch giving $\pi_{{\bar W},Q}$.

Fix $a \in [F(\kappa) \cap F^*(\kappa)]^{< \omega}$,
and $X \in {\cal P}([\kappa]^{{\rm Card}(a)}) \cap {\bar
W} = {\cal P}([\kappa]^{{\rm Card}(a)}) \cap {\tilde
W}$. Pick ${\bar \tau}$ and ${\vec \zeta} \in \Xi'$ such that
$X = {\bar \tau}^{\bar W}[{\vec \zeta}] \cap [\kappa]^{{\rm Card}(a)}$.
Then we get that $a \in F(X)$ if and only if $a \in \pi(X)$ if and only if $a \in \pi_{{\tilde W},Q}
\circ \pi(X)$ (as $\pi_{{\tilde W},Q} {\upharpoonright} F(\kappa) = {\rm id}$) if and only if
$a \in \pi_{{\bar W},Q}(X)$ (as $\pi_{{\tilde W},Q}
\circ \pi(X) = {\bar \tau}^Q[{\vec \zeta}] \cap [\pi_{{\tilde W},Q}
\circ \pi(\kappa)]^{{\rm Card}(a)}$ and $\pi_{{\bar W},Q}(X) =
{\bar \tau}^Q[{\vec \zeta}] \cap [\pi_{{\bar
W},Q}(\kappa)]^{{\rm Card}(a)}$) if and only if $a \in F^*(X)$.
Hence
$F$ and $F^*$ are compatible.

But we can say more. Using
the initial segment for ${\cal B}$ or for the the model where
$F^*$ comes from we can easily deduce that in fact $F^* = F$. However, then,
$\rho_1({\cal J}^{\bar W}_{\tilde \lambda}) < F(\kappa)$. On the other hand,
$F(\kappa)$ is a cardinal in ${\bar W}$, by how $F$ was obtained. Contradiction!

\bigskip
\hfill $\square$ (Subclaim 1)

\bigskip
We are now going to define $\{ (\kappa,x) \colon x \in \Gamma_\kappa \}$ in such a way
that $\Gamma_\kappa = \Xi \Rightarrow
\Gamma_{\kappa+1} = \Xi'$ for all ordinals $\kappa$,
where $\Xi \mapsto \Xi'$ is as above.

Let $A_0'$ denote the class of limit points of $A_0$. We closely follow
\cite[Appendix II]{MOZ}. 
We'll construct $\Gamma_i^\delta$, $Y_i^\delta$ for all $i \in {\rm OR}$
and certain $\delta \in A_0'$, by recursion on $i \in {\rm OR}$. We shall inductively
maintain that the following statements
are true (whenever the objects referred to are defined).

\begin{equation}\label{TH1}
\Gamma^\delta_i \subset \delta {\rm , \ and \ } \Gamma^\delta_i \subset Y^\delta_i
\prec_{\Sigma_1} {\cal J}_\delta^{{K^c_{ \{ \omega_1 \} }}} {\rm , }
\end{equation}
\begin{equation}\label{TH2}
i \geq j \Rightarrow \Gamma^\delta_i \subset \Gamma^\delta_j \wedge
Y^\delta_i \subset Y^\delta_j {\rm , }
\end{equation}
\begin{equation}\label{TH3}
\delta' \geq \delta \Rightarrow \Gamma^\delta_i = \Gamma^{\delta'}_i \cap \delta
\wedge
Y^\delta_i = Y^{\delta'}_i \cap {\cal J}_\delta^{{K^c_{ \{ \omega_1 \} }}}
\end{equation}

To commence, we let $\Gamma_0^\delta$ be all ordinals $< \delta$ of cofinality
$\omega$, and we let $Y_0^\delta = {\cal J}^{{K^c_{ \{ \omega_1 \} }}}_\delta$, for all $\delta \in A_0'$.
If $\lambda$ is a limit ordinal then we let $\Gamma_\lambda^\delta =
\bigcap_{i<\lambda} \ \Gamma_i^\delta$ and $Y_\lambda^\delta =
\bigcap_{i<\lambda} \ Y_i^\delta$, for all $\delta$ such that $\Gamma_i^\delta$
and $Y_i^\delta$ are defined whenever $i<\lambda$ (otherwise $\Gamma_\lambda^\delta$
and $Y_\lambda^\delta$ will be undefined).
Notice that (\ref{TH1}) will then be true for $\lambda$, as any ${\cal
J}^{{K^c_{ \{ \omega_1 \} }}}_\delta$ has a $\Sigma_1$-definable $\Sigma_1$ Skolem function.

Now suppose that $\Gamma_i^\delta$
and $Y_i^\delta$ are defined. If $i \notin Y_i^\delta$ then we put
$\Gamma_{i+1}^\delta = \Gamma_i^\delta$
and $Y_{i+1}^\delta = Y_i^\delta$. Suppose now that $i \in Y_i^\delta$. Consider
$$\sigma \colon W = W^\delta_i \cong Y_i^\delta \prec_{\Sigma_1} {\cal
J}^{{K^c_{ \{ \omega_1 \} }}}_\delta.$$ Suppose that $W \cap \delta = \delta$
(otherwise $\Gamma_{i+1}^\delta$
and $Y_{i+1}^\delta$ will be undefined). If there are no $\alpha < \delta$ and
$F$ such that
$${\cal B} = ({\cal J}_\alpha^W,F)$$ is a beaver with ${\rm c.p.}(F) = \sigma^{-1}(i)$ then
we let $\Gamma_{i+1}^\delta = \Gamma_i^\delta$ and $Y_{i+1}^\delta = Y_i^\delta$.
Otherwise let $\alpha < \lambda$ be least such that $F$ is the unique (by
\ref{double-beavers}) $F$ so that ${\cal B}$ as above is a beaver. Let $$i_F \colon
W \rightarrow_F {\tilde W} =
{\tilde W}^\delta_i = {\rm Ult}_0(W;F).$$
Notice that ${\tilde W}$ must be $0$-iterable (as ${\cal B}$ is a beaver),
and ${\tilde W} \cap {\rm OR} = \delta$.
Let ${\cal U}$ and ${\cal T}$ denote the
iteration trees arising from the coiteration of $W$ with ${\tilde W}$. Suppose that:

\bigskip
$\bullet \ \ $ ${\cal D}^{\cal U} \cap (0,\infty]_U =
{\cal D}^{\cal T} \cap (0,\infty]_T = \emptyset$, and

$\bullet \ \ $ ${\cal M}^{\cal U}_\infty \cap {\rm OR} = {\cal M}^{\cal T}_\infty \cap {\rm OR} =
\delta$.

\bigskip
(Otherwise $\Gamma_{i+1}^\delta$
and $Y_{i+1}^\delta$ will be undefined.)
Now put
$$\Gamma_{i+1}^\delta = \{ \xi \in \Gamma_i^\delta \colon \xi = \sigma(\xi) =
\pi^{\cal U}_{0 \infty}(\xi) = \pi^{\cal T}_{0 \infty} \circ i_F(\xi) \} {\rm , \
and}$$
$$Y_{i+1}^\delta = \{ \sigma(x) \colon x \in W \wedge \pi^{\cal U}_{0 \infty}(x) \in
{\rm ran}(\pi^{\cal T}_{0 \infty} \circ i_F) \}.$$

This finishes the recursive definition of $\Gamma_i^\delta$ and $Y_i^\delta$.
Now set $D_i = \{ \delta \colon \Gamma_i^\delta$ and $Y^\delta_i$ are defined$\}$,
and let $$\Gamma_i = \bigcup_{\delta \in D_i} \ \Gamma_i^\delta {\rm , \ and}$$
$$Y_i = \bigcup_{\delta \in D_i} \ Y_i^\delta.$$
Besides (\ref{TH1}), (\ref{TH2}), (\ref{TH3}), we also want to verify, inductively,
that:

\begin{equation}\label{TH4}
\forall i \in {\rm OR} \ \exists \eta \ \ ( \ D_i \setminus \eta = A_0' \setminus
\eta \ ) {\rm , }
\end{equation}
\begin{equation}\label{TH5}
H^{{K^c_{ \{ \omega_1 \} }}}(\Gamma_i) \subset Y_i \prec {K^c_{ \{ \omega_1 \} }}  {\rm , }
\end{equation}
\begin{equation}\label{TH6}
(\Gamma_i \colon i \in {\rm OR}) {\rm \ is \ as \ in \ the \ statement \ of \ Claim \ 1.}
\end{equation}

It is now straightforward to see that in order to inductively prove that
(\ref{TH1}) through (\ref{TH6}) hold for all $i \in {\rm OR}$
it suffices to show the following (cf. \cite[Appendix II]{MOZ}).

\begin{proposition}\label{proposit}
Suppose that
(\ref{TH1}) through (\ref{TH6}) hold for some $i \in {\rm OR}$. Let $\delta \in A_0'$ be large
enough. Suppose that
$Y_{i+1}^\delta \not= Y_i^\delta$, and let ${\cal B} = ({\cal
J}_\alpha^{W_i^\delta},F)$ be the beaver used to define $\Gamma_{i+1}^\delta$ and
$Y_{i+1}^\delta$.
Let ${\cal U}$, ${\cal T}$ denote the coiteration of $W_i^\delta$ with
${\tilde W}_i^\delta = {\rm Ult}_0(W^\delta_i;F)$. Let $\theta = {\rm lh}({\cal U}) = {\rm lh}({\cal
T})$.
Let $$\sigma' \colon W_i \cong Y_i \prec {K^c_{ \{ \omega_1 \} }}.$$ Let
${\cal U}'$, ${\cal T}'$ denote the coiteration of $W_i$ with ${\rm Ult}_0(W_i;F)$.
Then $\delta \in D_{i+1}$, $\theta \in [0,\infty]_{U'} \cap [0,\infty]_{T'}$, and
$\pi_{\theta \infty}^{{\cal U}'} {\upharpoonright} \delta =
\pi_{\theta \infty}^{{\cal T}'} {\upharpoonright} \delta = {\rm id}.$
\end{proposition}

It is clear that ${\cal U}$ and ${\cal T}$ are initial segments of
${\cal U}'$ and ${\cal T}'$, respectively. What \ref{proposit}
says is that the rest of the coiteration is above $\delta$ (and that
$\delta \in D_{i+1}$).

\bigskip
{\sc Proof} of \ref{proposit}.
It is easy to see that $\delta \in D_{i+1}$. Suppose
that $\theta \notin [0,\infty]_{U'} \cap [0,\infty]_{T'}$. Let $\alpha+1$ be least
in $(0,\infty]_{U'} \setminus \theta$, and let $\alpha^* = U'$-${\rm pred}(\alpha+1)$.
Let $\beta+1$ be least
in $(0,\infty]_{T'} \setminus \theta$, and let $\beta^* = T'$-${\rm pred}(\alpha+1)$.
Then $\mu = {\rm min} \{ {\rm c.p.}(\pi_{\alpha^* \infty}^{{\cal U}'}) ,$
${\rm c.p.}(\pi_{\beta^* \infty}^{{\cal T}'}) \} < \delta$.
As $\delta \in A_0'$, we in fact get
essentially as in the proof of \ref{on-the-sequence} Claim 2 Subclaim 1 that
$\mu = {\rm c.p.}(\pi_{\alpha^* \infty}^{{\cal U}'}) =
{\rm c.p.}(\pi_{\beta^* \infty}^{{\cal T}'})$. But then \ref{hull-prop-at-strongs} gives a
standard contradiction (see the proof of \ref{hull-prop-at-strongs})!
This yields \ref{proposit}.

\bigskip
Now the previous construction together with \ref{double-beavers} as well as the
set theoretical definability of ``beaver-hood'' (see \cite[\S 2]{kairalf})
finishes the proof of Claim 1.

\bigskip
\hfill $\square$ (Claim 1)

\bigskip
Now let us assume that ${\rm Def}({K^c_{ \{ \omega_1 \} }},A_0)$ is bounded,
$\beta = {\rm sup} ({\rm Def}({K^c_{ \{ \omega_1 \} }},A_0))$, say.
We aim to derive a contradiction. Fix $(\Gamma_\kappa \colon \kappa \in {\rm OR})$ as given
by Claim 1. 

Let $b_\kappa$ denote the least ordinal in $H^{{K^c_{ \{ \omega_1 \} }}}(\Gamma_\kappa)
\setminus
\beta$ for $\kappa \geq \beta$ (hence $b_\kappa
\geq \kappa$). By further thinning out the $\Gamma_\kappa$'s if
necessary we may assume without loss of generality that $b_\kappa < b_{\kappa+1}$ for
$\kappa \geq \beta$.

\bigskip
{\bf Claim 2.}
There is some $\nu > \beta$, a limit of $b_\kappa$'s, such
that $\nu \in H^{{K^c_{ \{ \omega_1 \} }}}(\nu \cup \Gamma_{\nu +1})$.

\bigskip
Given Claim 2, the proof of \ref{K-weasel} can be completed as
follows (cf. \cite[p. 38]{CMIP}).
Fix $\nu$ as in Claim 2.
Let $\nu = \tau^{{K^c_{ \{ \omega_1 \} }}}[{\vec \alpha},{\vec
\rho}]$ where ${\vec \alpha} < \nu$ and ${\vec \rho} \in
\Gamma_{\nu+1}$. Then ${\vec \alpha} < b_\kappa$ for some
$b_\kappa < \nu$.
Hence
$${{K^c_{ \{ \omega_1 \} }}} \models \ \exists {\vec \alpha} < b_\kappa ( b_\kappa <
\tau^{{K^c_{ \{ \omega_1 \} }}}[{\vec \alpha},{\vec
\rho}] < b_{\nu + 1}).$$
But $b_\kappa$,
${\vec
\rho}$, $b_{\nu+1}
\in \Gamma_\kappa$, and $H^{{K^c_{ \{ \omega_1 \} }}}(\Gamma_\kappa) \prec {{K^c_{ \{ \omega_1 \} }}}$, so
that there is some ${\vec \alpha}^\star \in b_\kappa \cap
H^{{K^c_{ \{ \omega_1 \} }}}(\Gamma_\kappa) \subset \beta$ with $$b_\kappa <
\tau^{{K^c_{ \{ \omega_1 \} }}}[{\vec \alpha}^\star,{\vec
\rho}] < b_{\nu + 1}.$$
By ${\vec \alpha}^\star \in \beta \cap H^{{K^c_{ \{ \omega_1 \} }}}(\Gamma_\kappa)$,
we have that ${\vec \alpha}^\star \in {\rm Def}({{K^c_{ \{ \omega_1 \} }}},A_0)$.
Hence $\tau^{{K^c_{ \{ \omega_1 \} }}}[{\vec \alpha}^\star,{\vec
\rho}] \in H^{{K^c_{ \{ \omega_1 \} }}}(\Gamma_{\nu+1})$ and $\beta \leq b_\kappa \leq
\tau^{{K^c_{ \{ \omega_1 \} }}}[{\vec \alpha}^\star,{\vec
\rho}] < b_{\nu+1}$. This
contradicts the definition of $b_{\nu+1}$!

\bigskip
{\sc Proof} of Claim 2. Let $C$ denote the class of all limit
points of $\{ b_\kappa \ \colon \ \beta \leq \kappa < {\rm OR} \}$. Of
course, $C$ is club in ${\rm OR}$.

Let us suppose Claim 2 to be false, so that for all $\kappa \in
C$ do we have
$$\sigma_\kappa \ \colon \ {{K^c_{ \{ \omega_1 \} }}}_\kappa \cong H^{{K^c_{ \{ \omega_1 \} }}}(\kappa
\cup \Gamma_{\kappa+1}) \prec {{K^c_{ \{ \omega_1 \} }}}$$
with ${\rm c.p.}(\sigma_\kappa) = \kappa$ and $\sigma_\kappa(\kappa) \leq
b_{\kappa+1}$.
Let $F_\kappa = \sigma_\kappa {\upharpoonright} {\cal P}(\kappa) \cap {{K^c_{ \{ \omega_1 \} }}}_\kappa$
for $\kappa \in C$.

The following is due to John Steel and is included here with his permission.

\bigskip
{\bf Subclaim 2 (Steel).} The class $\{ \kappa \in C \ \colon \
F_\kappa$ is countably complete $\}$ contains a club.

\bigskip
{\sc Proof.} Suppose not, so that $$S_0 = \{ \kappa \in C \ \colon \
F_\kappa {\rm \ is \ not \ countably \ complete \ } \}$$ is stationary.
We may then pick $((a^n_\kappa,X^n_\kappa) \
\colon \ n < \omega \wedge \kappa \in C)$ such
that $a^n_\kappa \in F_\kappa(X^n_\kappa)$, but there is no order
preserving $\tau \colon \cup_{n < \omega} \ a^n_\kappa 
\rightarrow \kappa$ with $\tau {\rm " } a^n_\kappa \in X^n_\kappa$
for any $\kappa \in C$.

We have to define a function $G \colon S_0 \rightarrow V$, saying how
the $a^n_\kappa$'s sit inside $\cup_n \ a^n_\kappa$. Let
$\tau_\kappa \ \colon \ \cup_n \ a^n_\kappa \rightarrow {\rm otp}
(\cup_n \ a^n_\kappa) < \omega_1$ be order preserving, and let $G(\kappa) =
({\rm otp}(\cup_n \ a^n_\kappa),(\tau_\kappa {\rm " } a^n_\kappa \colon n <
\omega))$.

Now by \ref{hull-prop}, there is some stationary $S_1 \subset S_0$
such that ${{K^c_{ \{ \omega_1 \} }}}$ has the hull property at every
$\kappa \in S_1$.
Hence
for any $n < \omega$ and $\kappa \in S$
may we pick a term $\tau^n_\kappa$ and ${\vec \alpha}^n_\kappa <
\kappa$ and ${\vec \gamma}^n_\kappa \in \Gamma_{\kappa+1}$ such that
$$\sigma_\kappa(X^n_\kappa) =
(\tau^n_\kappa)^{{K^c_{ \{ \omega_1 \} }}}[{\vec \alpha}^n_\kappa,{\vec
\gamma}^n_\kappa].$$
We now consider
$$F(\kappa) = (G(\kappa),(\tau^n_\kappa,{\vec
\alpha}^n_\kappa,((\tau^n_\kappa)^{{K^c_{ \{ \omega_1 \} }}}[{\vec \alpha},{\vec
\gamma}^n_\kappa] \cap \beta \colon {\vec \alpha} < \beta)
\colon n < \omega)) {\rm , }$$
being essentially a regressive function on $S$.
By Fodor, there is an unbounded $D \subset S$ such that $F$ is
constant on $D$. Fix $\kappa < \kappa' \in D$
with $\sigma_\kappa(\kappa) \leq b_{\kappa+1} \leq \kappa'$.

Notice first that we have an
order preserving map $$\tau \
\colon \cup_n \ a^n_{\kappa'} {\rightarrow} \cup_n \
a^n_\kappa$$ with $\tau {\rm " } a^n_{\kappa'} =
a^n_\kappa$ for all $n < \omega$, due to the fact that $G(\kappa) =
G(\kappa')$. Hence we would have a contradiction to the choice of
$(a^n_{\kappa'},X^n_{\kappa'} \ \colon \ n < \omega)$ if we were able
to show that $a^n_\kappa \in X^n_{\kappa'}$ for all $n < \omega$.

In order to do this, as $a^n_\kappa < \sigma_\kappa(\kappa) \leq
b_{\kappa+1}$ and $a^n_\kappa \in \sigma_\kappa(X^n_\kappa)$ by the
choice of $(a^n_\kappa , X^n_\kappa \ \colon \ n < \omega)$ (and
because $X^n_{\kappa'} \cap b_{\kappa+1} =
\sigma_{\kappa'}(X^n_{\kappa'}) \cap b_{\kappa+1}$ by $\kappa' \geq
b_{\kappa+1}$),
it suffices to establish that
$$\sigma_\kappa(X^n_\kappa) \cap b_{\kappa+1} =
\sigma_{\kappa'}(X^n_{\kappa'}) \cap b_{\kappa+1}.$$

We have that $\tau^n_\kappa = \tau^n_{\kappa'} = \tau$ and
${\vec \alpha}^n_\kappa = {\vec \alpha}^n_{\kappa'} = {\vec \alpha}$
for some $\tau$, ${\vec \alpha}$, and $\sigma_\kappa(A^n_\kappa)
= \tau^{{K^c_{ \{ \omega_1 \} }}}({\vec \alpha},{\vec \gamma}^n_\kappa)$ and
$\sigma_{\kappa'}(X^n_{\kappa'}) =
\tau^{{K^c_{ \{ \omega_1 \} }}}({\vec \alpha},{\vec \gamma}^n_{\kappa'})$.
So if $\sigma_\kappa(A^n_\kappa) \cap b_{\kappa+1} \not=
\sigma_{\kappa'}(X^n_{\kappa'}) \cap b_{\kappa+1}$
then
$${{K^c_{ \{ \omega_1 \} }}} \models \ \exists {\vec \alpha} < b_{\kappa+1} \
\tau^{{K^c_{ \{ \omega_1 \} }}}[{\vec \alpha},{\vec \gamma}^n_\kappa] \not=
\tau^{{K^c_{ \{ \omega_1 \} }}}[{\vec \alpha},{\vec \gamma}^n_{\kappa'}].$$
As $b_{\kappa+1}$, ${\vec \gamma}^n_\kappa$, ${\vec
\gamma}^n_{\kappa'} \in \Gamma_{\kappa+1}$, there is some witness
${\vec \alpha}^\star \in h^{{K^c_{ \{ \omega_1 \} }}}(\Gamma_{\kappa+1})$.
But ${\vec \alpha}^\star \in h^{{K^c_{ \{ \omega_1 \} }}}(\Gamma_{\kappa+1}) \cap b_{\kappa+1}$,
so that ${\vec \alpha}^\star \in {\rm Def}({{K^c_{ \{ \omega_1 \} }}},A_0) \subset \beta$.

Now we have that
$${{K^c_{ \{ \omega_1 \} }}} \models \ \exists d < b_{\kappa+1} \ ( d \in
\tau^{{K^c_{ \{ \omega_1 \} }}}[{\vec \alpha}^\star,{\vec
\gamma}^n_\kappa] \Leftrightarrow
d \notin
\tau^{{K^c_{ \{ \omega_1 \} }}}[{\vec \alpha}^\star,{\vec
\gamma}^n_{\kappa'}]).$$
Here we have that ${\vec \alpha}^\star$,
${\vec \gamma}^n_\kappa$, ${\vec
\gamma}^n_{\kappa'} \in \Gamma_{\kappa+1}$, so that there is some
witness $d \in h^{{K^c_{ \{ \omega_1 \} }}}(\Gamma_{\kappa+1})$.
But $d < b_{\kappa+1}$ and $b_{\kappa+1} \cap h^{{K^c_{ \{ \omega_1 \} }}}(\Gamma_{\kappa+1}) =
{\rm Def}({{K^c_{ \{ \omega_1 \} }}},A_0) \subset \beta$, so that we may conclude that $d < \beta$.
I.e., $\tau^{{K^c_{ \{ \omega_1 \} }}}[{\vec \alpha}^\star,{\vec \gamma}^n_\kappa] \cap \beta
\not= \tau^{{K^c_{ \{ \omega_1 \} }}}[{\vec \alpha}^\star,{\vec \gamma}^n_{\kappa'}] \cap
\beta$.

However, as ${\vec \alpha}^\star < \beta$, $\kappa < \kappa' \in D$
gives us that
$\tau^{{K^c_{ \{ \omega_1 \} }}}[{\vec \alpha}^\star,{\vec \gamma}^n_\kappa] \cap \beta
= \tau^{{K^c_{ \{ \omega_1 \} }}}[{\vec \alpha}^\star,{\vec \gamma}^n_{\kappa'}] \cap
\beta$.
Contradiction!

\bigskip
\hfill $\square$ (Subclaim 2)

\bigskip
We may hence pick $\kappa \in C$ such that $F_\kappa$ is countably complete and ${{K^c_{ \{ \omega_1 \} }}}$ has
the hull property at $\kappa$. We then have an
immediate contradiction with \ref{weak-max}.

\bigskip
\hfill $\square$ (Claim 2)

\hfill $\square$ (\ref{K-weasel})

\begin{defn}\label{defn-K} $( \ \lnot \ \zerohandgrenade \ )$
We shall write $K$ for $K(A_0)$. $K$ is called the core model below $\zerohandgrenade$.
\end{defn}

\begin{lemma}\label{K-is-full} $( \ \lnot \ \zerohandgrenade \ )$
$K$ is full, i.e., whenever $F$ is a countably complete extender such
that $({\cal J}^K_\alpha,F)$ is a premouse then $F = E^K_\alpha$.
\end{lemma}

{\sc Proof.} Let $W$ be such that $\alpha \subset {\rm Def}(W,A_0)$.
As $F$ is countably complete, ${\rm Ult}(W,F)$ is
iterable. In particular, ${\cal J}^W_\alpha = {\cal J}^K_\alpha$.
Then \cite[Lemma 1.6]{kairalf} gives the desired result.

\bigskip
\hfill $\square$

Combined with \ref{W-is-universal} this gives at once:

\begin{cor}\label{K-is-universal}
$( \ \lnot \ \zerohandgrenade \ )$
$K$ is universal.
\end{cor}

It is now clear that the proofs in
\cite{MiSchSt} and \cite{Mi-Sch} give the following.

\begin{thm} $( \ \lnot \ \zerohandgrenade \ )$
Let $\beta \geq \omega_2$. Then ${\rm cf}^V(\beta^{+K}) \geq {\rm Card}^V(\beta)$.
\end{thm}

As in \cite[Theorem 1.1]{schdl}, we'll also have the following.

\begin{thm}\label{KcapHC} $( \ \lnot \ \zerohandgrenade \ )$
Let $1 \leq \alpha \leq \omega_1$, and suppose that $${\cal J}^K_{\omega_1} \models
{\rm \ there \ are \ } < \omega \alpha {\rm \ many \ strong \ cardinals.}$$
Then the set of reals coding ${\cal J}^K_{\omega_1}$ is an element of
$J_{1+\alpha}({\mathbb R})$.
\end{thm}

For the projective level, \cite[Theorems 3.4 and 3.6]{kairalf} 
(which are due to the
present author)
give the following
refinements.

\begin{thm}\label{KcapHC2} $( \ \lnot \ \zerohandgrenade \ )$
Let $n < \omega$.
Suppose that $${\cal J}^K_{\omega_1} \models
{\rm \ there \ are \ exactly \ } n {\rm \ strong \ cardinals.}$$
If $\omega_1$ is inaccessible in $K$ then ${\cal J}^K_{\omega_1}$ is (lightface)
$\Delta^1_{n+5}$ in the codes; and if  $\omega_1$ is a successor cardinal in $K$ then
${\cal J}^K_{\omega_1}$ is $\Delta^1_{n+4}(x)$ in the codes for some $x \in {\mathbb
R}$ coding an initial segment of ${\cal J}^K_{\omega_1}$.
\end{thm}

Jensen has shown in \cite[\S 5.3 Lemma 5]{jensen} 
that if $0^\P$ does not exist then
every universal weasel is an iterate of $K$ (via a normal iteration).
His proof in fact straightforwardly generalizes to the situation where
$K^c$ does not have a cardinal which is strong up to a measurable cardinal, 
i.e., a
measurable cardinal $\kappa$ such that $${\cal
J}^{K^c}_\kappa \models {\rm \ " there \ is \ a \ strong \ cardinal."
}$$ On the other hand, we are now going to construct
an example of a universal
weasel not being an iterate of $K$, assuming that
$K^c$ does have a "strong up to
a measurable."
The example
is due to John Steel and is included here with his permission.
(Recall that \cite[p. 86]{CMIP} had already
shown indirectly that such a
weasel has to exist if $K \cap {\rm HC}$ is not $\Delta^1_5$.)

\begin{lemma}\label{example} {\bf (Steel, } $ \ \lnot \ \zerohandgrenade \ $ {\bf ) }
Suppose that
$\mu < \kappa$ are such that $\kappa$
is measurable in $K$ and ${\cal J}^K_\kappa \models$ ``$\mu$
is a strong cardinal.'' Then if $H$ is
$Col(\omega,\kappa^{+K})$-generic over $K$, inside $K[H]$ there exists
a universal
weasel $W$ which is not an iterate of $K$. In fact, $W$ may be
chosen such that $W \triangleright ({\cal J}_{\kappa^{+K}}^K,F)$
for some extender $F$.

Moreover, if $\mu$ is the only ${\bar \mu}$ with ${\cal J}^K_\kappa \models$
``${\bar \mu}$
is a strong cardinal,'' then no universal
$W' \triangleright ({\cal J}_{\kappa^{+K}}^K,F)$ has
the definability property at $\mu$.
\end{lemma}

{\sc Proof.} First fix $g \in K[H]$ being $Col(\omega,\mu^{++K})$-generic over
$K$. We shall use a theorem of Woodin (unpublished) which tells us
that in $K[g]$ there is a tree $T_3$ projecting to a universal
$\Pi^1_3$-set of reals in all further small extensions of $K[g]$,
where "small" means that they are
obtained by forcing with some $P \in {\cal J}^K_\alpha[g]$
where $\alpha = E^K_\beta(\mu)$ for some $E^K_\beta$ with critical
point $\mu$.
This implies that in $K[g]$ there is a tree $T$ projecting to the set
of all real codes for mice
${\cal M}$ such that $K^c({\cal M})$ exists,
where $T$ works in all such extensions; the reason is that this is a
$\Pi^1_3$-set of reals (cf. \cite[Corollary 2.18 (a)]{kairalf}).

Let $E = E^K_\nu$ witness that $\kappa$ is measurable in $K$, and
let $\pi \colon K \rightarrow_E M$ with $M$ being transitive. We may
extend $\pi$ to ${\tilde \pi} \colon K[g] \rightarrow M[g]$, and
by elementarity, ${\tilde \pi}(T) \in M[g]$ is a tree projecting to
the set of all real codes for mice ${\cal M}$ such that $K^c({\cal
M})$ exists, with ${\tilde \pi}(T)$ working in
all further small extensions of $M[g]$.

Now let $G$ be $Col(\omega,\nu)$-generic over $K[g]$, so that
$G$ is $Col(\omega,\nu)$-generic over $M[g]$,
too. Notice that
${\tilde \pi}(T)$ still "works" in $M[g][G]$. In $M[g][G]$, let $x
\in {\mathbb R}$ code ${\cal J}^M_\nu = {\cal J}^K_\nu$. We may build a
tree $T^\star \in M[g][G]$ searching for a pair $(y,f)$ such that
$(x \oplus y , f) \in [{\tilde \pi}(T)]$. Here, by $a \oplus b$ we mean a
canonical code for $({\cal N},F)$ obtained from $(a,b)$ where $a$
codes ${\cal N}$ and $b$ codes $F$.

We claim that $T^\star$ is ill-founded (in $K[g][G]$, and hence in
$M[g][G]$). Let $T_n$ denote $T$ up to the $n^{{\rm th}}$ level.
Well, if $y \in {\mathbb R} \cap K[g][G]$ codes $E_\nu^K$
then $x \oplus y \in p[T]$, and hence if
$f \in K[g][G]$ is such that $(x \oplus y,f) \in [T]$ then
$$\forall n \ (x \oplus y {\upharpoonright} n, f {\upharpoonright} n) \in
T_n {\rm , \ hence}$$
$$\forall n \ (x \oplus y {\upharpoonright} n, \pi(f {\upharpoonright} n))
\in
T_n {\rm , }$$
so that $(y, \cup_n \pi(f {\upharpoonright} n)) \in T^\star$.

Thus in $M[g]$, the following holds true:
$$\forces_{Col(\omega,\nu)} {\rm \ ``there \ is \ an \ extender \ }
F {\rm \ such \ that \ } K^c(({\cal J}^M_\nu,F)) {\rm \ exists,''}$$
and by elementarity of ${\tilde \pi}$,
in $K[g]$ we have that
$$\forces_{Col(\omega,\kappa^+)} {\rm \ ``there \ is \ an \ extender \ }
F {\rm \ such \ that \ } K^c(({\cal J}^K_{\kappa^+},F))
{\rm \ exists.''}$$

Now let ${\bar H} \in K[H]$ be $Col(\omega,\kappa^+)$-generic over $K[g]$.
By what we have shown, in $K[g][{\bar H}]$, and hence in
$K[H]$ we may pick some $F$ such that $({\cal J}^K_{\kappa^+},F)$ is a
premouse and $W = K^c(({\cal J}^K_{\kappa^+},F))$ exists.
Notice that $\rho_1(({\cal J}^K_{\kappa^+},F)) < \kappa$.

It is now easy to see that $W$ cannot be an iterate of $K$, as
$\kappa^{+K}$ is a cardinal in every such iterate, whereas it is not a
cardinal in $W$.

Now suppose that the universal weasel $W' \triangleright ({\cal J}^K_{\kappa^+},F)$
would have the definability property at $\mu$, and that
$\mu$ is the only ${\bar \mu}$ with ${\cal J}^K_\kappa \models$
``${\bar \mu}$
is a strong cardinal.'' Then the coiteration of $W'$ and $K$ would be above $\kappa$
on both sides by \cite[Corollary 1.5]{kairalf}. 
We'd thus get that $\kappa^{+W'} =
\kappa^{+K}$. Contradiction!

\bigskip
\hfill $\square$ (\ref{example})

\section{An application.}

By classical results, Projective Determinacy (${\sf PD}$, for
short) implies that every projective set of reals is Lebesgue measurable and
has the property of Baire, and that every projective subset of ${\mathbb R}
\times {\mathbb R}$ has a uniformizing function with a projective
graph (cf. for example \cite[6A.16]{Mo}, \cite[6A.18]{Mo}, 
and \cite[6C.5]{Mo}). In \cite{CSPU}, Woodin asked whether
${\sf PD}$ is actually equivalent with 
the fact that
every projective set of reals is Lebesgue measurable and
has the property of Baire, and that every projective subset of ${\mathbb R}
\times {\mathbb R}$ has a uniformizing function with a projective
graph. This question became the $12^{{\rm th}}$ Delfino problem (see \cite{KeMaSt}), and
it was widely believed to have an affirmative answer, -- a belief
which
was supported by a theorem of Woodin's according to which at least \boldmath
$\Pi$\unboldmath$^1_1$ determinacy does follow from these ``analytic''
consequences of ${\sf PD}$.

However, Steel in the
fall of 1997 proved that the answer to Woodin's question is negative.
This left open the question of the exact consistency strength of the
assumption in the statement of the $12^{{\rm th}}$ Delfino problem, that is, of
above-mentioned ``analytic'' consequences of ${\sf PD}$.
Our application, which yields the following theorem \ref{appl1},
determines that strength in terms of large cardinal 
assumptions.\footnote{It was actually the wish to complete the proof of
this theorem which motivated our work reported in the
previous sections.}

\begin{thm}\label{appl1}
The following theories are equiconsistent.

(1) ${\sf ZFC} +$ ``every projective set of reals is Lebesgue measurable and
has the property of Baire" $+$ ``every projective subset of ${\mathbb R}
\times {\mathbb R}$ has a uniformizing function with a projective
graph,"

(2) ${\sf ZFC} +$ ``every projective set of reals is Lebesgue measurable and
has the property of Baire" $+$ ``the pointclass consisting of the
projective sets of reals has the scale property," and

(3) ${\sf ZFC} +$ ``there are infinitely many
cardinals $\kappa_0 <
\kappa_1 < \kappa_2 < ...$ each of which is
$({\rm sup}_{n<\omega} \ \kappa_n)^+$-strong, i.e.,
for every $n<\omega$ and for every $X \subset {\rm sup}_{n<\omega} \ \kappa_n$
there is an elementary embedding
$\pi \ \colon \ V \rightarrow M$ with $M$ being transitive, ${\rm c.p.}(\pi)
= \kappa_n$, and $X \in M$."
\end{thm}

{\sc Proof.} Con(3) $\Rightarrow$ Con(2)
was shown by Steel (building on
work of Woodin, cf. \cite{steel}). (2) $\Rightarrow$ (1) is 
classical (cf. for example \cite[4E.3]{Mo}).
Con(1) $\Rightarrow$ Con(3) was shown in \cite{kairalf} 
(building on
Woodin's \cite{CSPU} and the author's work on the complexity of
$K \cap {\rm HC}$, cf. \ref{KcapHC} and \ref{KcapHC2} above; 
cf. \cite[Theorem 4.2]{kairalf}) with 
${\sf ZFC}$ replaced
by ${\sf ZFC} +$ ``${\mathbb R}^\#$ exists" (the authors of
\cite{kairalf} used the core model theory of \cite{CMIP} which needs
the existence of a measurable cardinal, or some substitute for it).
By our work done here in
the previous sections, ${\mathbb R}^\#$
is now no longer needed for
running the arguments of \cite{kairalf}.

\bigskip
\hfill $\square$ (\ref{appl1})

\bigskip
The reader will have no problems with finding the appropriate reformulations of
\cite[Theorems 1.1, 1.4, and 1.5]{schdl} in light of our work done here.

We finally remark that the following variants 
of the $12^{{\rm th}}$ Delfino problem remain wide open.

\bigskip
{\bf Question.} Suppose that every projective set of reals is Lebesgue 
measurable and has the property of Baire.
Suppose that either

(a) for each $n<\omega$ does 
every \boldmath $\Pi$\unboldmath$^1_{2n+1}$ subset of ${\mathbb R} \times 
{\mathbb R}$ have a uniformizing function with a graph in  
\boldmath $\Pi$\unboldmath$^1_{2n+1}$, or else

(b) every lightface projective subset of 
${\mathbb R} \times 
{\mathbb R}$ has a uniformizing function with a graph which is lightface projective.

Does then Projective Determinacy hold?

\bigskip
Steel has shown (cf. \cite[Corollary 7.14]{CMIP}) that variant (a) implies 
\boldmath $\Delta$\unboldmath$^1_2$ determinacy. It might be reasonable to
conjecture that variant (a) in fact implies ${\sf PD}$.
On the other hand it is conceivable that a refinement of coding arguments as
being developed in \cite{syralf} will establish that variant (b) does not
imply \boldmath $\Delta$\unboldmath$^1_2$ determinacy.      



\begin{thebibliography}{99}

\bibitem{dodd} T. Dodd, {\it Strong cardinals}, handwritten notes, 1987.

\bibitem{dodd-jensen} T. Dodd, R. B. Jensen,
{\it The core model}, Ann. Math. Logic {\bf 20} (1981), pp. 43
- 75.

\bibitem{syralf} S. D. Friedman and R.-D. Schindler, {\it Universally Baire
sets and definable well-orderings of the reals}, Journal Symb. 
Logic, submitted.

\bibitem{kairalf} K. Hauser and
R.-D. Schindler, {\it Projective uniformization
revisited},
Ann. Pure Appl. Logic {\bf 103} (2000), pp. 109 - 153.

\bibitem{classic} R. B. Jensen, {\it The fine structure of the constructible
hierarchy}, Ann. Math. Logic {\bf 4} (1972), pp. 229 - 308.

\bibitem{jensen2} R. B. Jensen, {\it Inner models and large cardinals}, Bulletin of
Symb. Logic {\bf 1} (1995), pp. 393 - 407.

\bibitem{MOZ} R. B. Jensen, {\it Measures of order zero},
handwritten notes.

\bibitem{jensen} R. B. Jensen, {\it The core model for
non-overlapping extender sequences}, handwritten notes.

\bibitem{NFS} R. B. Jensen, {\it A new fine structure for higher
core models}, handwritten
notes.

\bibitem{CAR} R. B. Jensen, {\it Corrections and remarks}, handwritten
notes.

\bibitem{ronald} R. B. Jensen, {\it T-mice}, handwritten notes.

\bibitem{KeMaSt} A.S. Kechris, D.A. Martin and J.R. Steel (eds.), {\em
Appendix: Victoria Delfino
Problems II,} in {\it Cabal Seminar 81-85}, Springer
Lecture Notes in Math. {\bf 1333} (1988), pp.221-224.

\bibitem{koepke} P. Koepke, {\it The short core model and its application},
Ph. D. Thesis, Freiburg i. Br.

\bibitem{st-loewe} B. L\"owe and J. Steel, {\it An introduction to core model
theory}, in: ``Sets and Proofs'' (Cooper and Truss, eds.), 
Cambridge University Press
1999, pp. 103 - 157.

\bibitem{mitchell} W.J. Mitchell, {\it The core model for sequences of measures I},
Math. Proc. Cambridge Philos. Soc. {\bf 95} (1984), pp. 228 - 260.

\bibitem{mitchell2} W.J. Mitchell, {\it The core model for sequences of measures II},
typescript.

\bibitem{mitchell3} W.J. Mitchell, {\it The core model up to a Woodin cardinal}, in:
``Logic, Methodology, and Philosophy of Science IX'' (Prawitz et al., eds.).

\bibitem{Mi-Sch} W.J. Mitchell and E. Schimmerling, {\it Covering without countable
closure}, Math Research Letters {\bf 2} (1995), pp. 595 - 609.

\bibitem{MiSchSt} W.J. Mitchell, E. Schimmerling and J.R. Steel,
{\it The covering lemma up to a Woodin cardinal}, Ann. Pure
Appl. Logic {\bf 84} (1997), pp. 219 - 255.

\bibitem{FSIT}  W. J. Mitchell and J. R. Steel, {\it Fine structure
and iteration trees}, Lecture Notes in Logic {\bf 3},
Springer-Verlag (1994), 130 pages.

\bibitem{Mo} Y. Moschovakis, {\it Descriptive set theory},
North-Holland (1980), xii + 637 pages.

\bibitem{maximality-paper} E. Schimmerling, J. Steel, {\it The maximality of the core
model}, to appear.

\bibitem{thesis} R.-D. Schindler, {\it The core model up to one strong cardinal}, Ph.
D. Thesis, Bonner Mathematische Schriften {\bf 295} (1997).

\bibitem{schdl} R.-D. Schindler, {\it Strong cardinals
and sets of reals in
$L_{\omega_1}({\mathbb R})$}, Math. Logic Quarterly, {\bf 45} (1999), pp. 361 - 369.

\bibitem{d-premice} R.-D. Schindler, {\it 
The core model for d-premice}, in preparation.

\bibitem{ssz} R.-D. Schindler, J. Steel, and M. Zeman, {\it Deconstructing inner model
theory}, Journal of Symbolic Logic, submitted.

\bibitem{CMIP} J. R. Steel, {\it The core model
iterability problem},
Lecture Notes in Logic {\bf 8}, Springer Verlag (1996), 112 pages.

\bibitem{CMMW} J. R. Steel, {\it Core models with more Woodin cardinals},
preprint.

\bibitem{steel} J. R. Steel, in preparation.

\bibitem{john-philip} J. R. Steel and P. D. Welch, {\it $\Sigma^1_3$ 
absoluteness and the second uniform indiscernible}, Israel
Journal of Mathematics {\bf 104} (1998), pp. 157 - 190.

\bibitem{CSPU}  W. H. Woodin, {\it On the
consistency strength of projective
uniformization}, in {\it Logic Colloquium '81}, J. Stern (ed.),
North-Holland (1982), pp. 365-383.

%
\bibitem{zeman}
M. Zeman, {\it Inner models and
strong cardinals}, de Gruyter Series in Logic and its Applications Vol. 4,
in preparation.

\end{thebibliography}
\end{document}